\documentclass[preprint,sort&compress]{elsarticle}

\usepackage[margin=1in]{geometry}

\usepackage{amsmath,amstext,amsbsy,amssymb,mathdots,amsfonts,amsthm}
\usepackage{relsize}
\usepackage{algorithm,algorithmicx,algpseudocode}
\usepackage{color}
\usepackage{subfig}
\usepackage{url}
\usepackage{caption}

\usepackage{tikz}
\usepackage{overpic}
\usetikzlibrary{positioning}
\usetikzlibrary{shapes,arrows}
\usepackage{relsize}
\usepackage{psfrag}
\usepackage{enumerate}
\usepackage{amsopn}
\usepackage{gensymb}
\usepackage{multirow}
\usepackage{colortbl}

\definecolor{Gray}{gray}{0.9}

\def \vu{ \boldsymbol{u} }

\def \vx{ \mathbf{x} }
\def \vp{ \boldsymbol{p} }

\def \vf{ \boldsymbol{f} }
\def \vxi{\boldsymbol{\xi}}
\def \vom{\boldsymbol{\omega}}

\def \vc { \boldsymbol{c} }

\def\calI{{ \mathcal I}}

\def \vq {\mathbf{q} }

\newcommand{\lf}{\left}
\newcommand{\rt}{\right}
\newcommand{\lp}{\lf (}
\newcommand{\rp}{\rt )}
\newcommand{\ep}{\varepsilon}

\def \vu{ \boldsymbol{u} }

\def \vx{ \boldsymbol{x} }

\def \vy{ \boldsymbol{y} }

\def \vf{ \mathbf{f} }
\def \vc{ \mathbf{c} }
\def \vd{ \mathbf{d} }

\def \vp{ \boldsymbol{p} }

\def \vc { \boldsymbol{c} }

\def \vu{ \boldsymbol{u} }

\def \vx{ \boldsymbol{x} }

\def \vf{ \mathbf{f} }

\def \vc{ \mathbf{c} }

\newcommand{\ds}{\displaystyle}
\newcommand{\Sph}{\mathbb{S}^2}

\def \vc { \mathbf{c} }

\newcommand{\Asp}{\widetilde{A}}
\newcommand{\csp}{\widetilde{\vc}}
\newcommand{\fsp}{\widetilde{\vf}}

\newdefinition{rmk}{Remark}

\newcommand{\revone}[1]{{\color{black} #1 }}
\newcommand{\revtwo}[1]{{\color{black}{#1}}}
\newcommand{\revthree}[1]{{\color{black}{#1}}}
\newcommand{\revfour}[1]{{\color{black}{#1}}}

\journal{JCP}

\begin{document}

\begin{frontmatter}

\title{Mesh-free Semi-Lagrangian Methods for Transport on a Sphere Using Radial Basis Functions}

\author[addr1]{Varun Shankar\corref{corresp}}
\address[addr1]{Department of Mathematics, University of Utah, UT, USA}
\ead{vshankar@math.utah.edu}
\cortext[corresp]{Corresponding Author}

\author[addr2]{Grady B. Wright}
\address[addr2]{Department of Mathematics, Boise State University, ID, USA}
\ead{gradywright@boisestate.edu}

\begin{abstract}
We present three new semi-Lagrangian methods based on radial basis function (RBF) interpolation for numerically simulating transport on a sphere.  The methods are mesh-free and are formulated entirely in Cartesian coordinates, thus avoiding any irregular clustering of nodes at artificial boundaries on the sphere and naturally bypassing any apparent artificial singularities associated with surface-based coordinate systems. \revone{For problems involving tracer transport in a given velocity field,} the semi-Lagrangian framework allows these new methods to avoid the use of any stabilization terms (such as hyperviscosity) during time-integration, thus reducing the number of parameters that have to be tuned. The three new methods are based on interpolation using 1) global RBFs, 2) local RBF stencils, and 3) RBF partition of unity.  For the latter two of these methods, we find that it is crucial to include some low degree spherical harmonics in the interpolants. Standard test cases consisting of solid body rotation and deformational flow are used to compare and contrast the methods in terms of their accuracy, efficiency, conservation properties, and dissipation/dispersion errors.  For global RBFs, spectral spatial convergence is observed for smooth solutions on quasi-uniform nodes, while high-order accuracy is observed for the local RBF stencil and partition of unity approaches.

\end{abstract}


\end{frontmatter}

\section{Introduction}
\label{sec:intro}

Radial basis function (RBFs) methods have been used for over a decade to solve partial differential equations (PDEs) on spheres. These methods can broadly be classified into \emph{global} RBF collocation methods~\cite{FlyerWright:2007,FornbergPiret:2008, FlyerWright:2009, Gia:2005, WFY}, RBF-generated finite difference (RBF-FD) methods~\cite{FoL11,FlyerLehto2012,Bollig12,Tillenius2015406,FlyerWrightFornberg}, and more recently RBF-partition of unity (RBF-PU) collocation methods~\cite{KevinThesis}.  Global RBF methods when used with infinitely-smooth RBFs show spectral convergence on smooth problems at the cost of \emph{dense} differentiation matrices; in contrast, RBF-FD and RBF-PU methods produce sparse differentiation matrices and high-order algebraic convergence rates.  All of these methods can use ``scattered'' nodes in their discretizations of a sphere, and have the benefit of being independent of any surface-based coordinate system. They thus avoid any unnatural grid clustering and do not suffer from any coordinate singularities.

In this paper, we present three new RBF methods with similar benefits for numerically solving the transport equation on the surface of a sphere in an incompressible velocity field.  For the unit sphere $\Sph$, this PDE is given 
\begin{align}
%\frac{\partial q}{\partial t} + \vu \cdot \nabla_{\Sph} q = 0,\quad \text{or}\quad \frac{Dq}{Dt} &= 0, \label{eq:advection}
\frac{Dq}{Dt} = 0,\quad \frac{D}{Dt}:=\frac{\partial }{\partial t} + \vu \cdot \nabla_{\Sph}, \label{eq:advection}
\end{align}
%\begin{align}
%\frac{\partial q}{\partial t} + \vu \cdot \nabla_{\Sph} q = 0,
%\end{align}
where $q$ is the scalar quantity being transported, $\vu$ is a surface divergence-free vector field that is tangent to $\Sph$, and $\nabla_{\Sph}$ denotes the surface gradient operator on $\Sph$.  Since global atmospheric flows are dominated by the horizontal advection process, the numerical solution to the transport problem is a fundamental part of any solver  for these flows.

Currently, all RBF discretizations for the transport equation (and more general hyperbolic equations like the shallow water equations) on a sphere suffer from the same drawback: the eigenvalues of the differentiation matrices corresponding to the surface gradient operator may in general have positive real parts, leading to either a slow or rapid onset of instability during a numerical simulation~\cite{FlyerWright:2007}. As of this paper, the only known approach to rectify this problem is to add an artificial hyperviscosity of the form $(-1)^{(k+1)}\gamma \Delta^k$ in the right hand side of the PDE, where $k\geq 2$ is an integer and $\gamma>0$ is some small real number that scales \emph{inversely} with the total number of nodes $N$. The intuition here is that higher powers of the Laplacian $\Delta$ will damp out the eigenvectors associated with the rogue eigenvalues of the discretized surface gradient, while leaving the others essentially untouched~\cite{FoL11,FlyerLehto2012,Safdari-Vaighani2015,FlyerNS}. With global RBFs, the hyperviscosity operator typically takes the form of $\gamma A^{-1}$, where $A$ is the global RBF interpolation matrix whose inverse mimics the properties of high powers of the Laplacian~\cite{FoL11}; a similar approach can be employed for the RBF-PU method~\cite{KevinThesis}.  Unfortunately, for a given PDE and node set on a sphere, the precise values of $\gamma$ and $k$ required to stabilize the numerical solution may need to be determined by trial and error, which can add to the computational expense of the method.   

A common way to naturally stabilize local Eulerian methods for transport is to use ``upwinding'', which uses dynamic direction-dependent information about the flow field. However, this form of upwinding typically requires an underlying mesh and so is impractical for truly mesh-free local methods like RBF-FD and RBF-PU collocation.  \revfour{An important class of methods that naturally possess upwinding are semi-Lagrangian (SL) schemes, widely acknowledged for their formal independence from the CFL stability condition~\cite{staniforth1991SL}}. \revthree{SL methods have successfully been used for simulating various problems in fluid dynamics, especially in numerical weather prediction and climate modeling, from simple tracer transport to more complex problems involving wide ranges of spatiotemporal scales and intricate forcing terms, e.g.~\cite{staniforth1991SL,StaniforthWoodJCP2008,SmolPudy92,SmolarkiewiczMargolin97,xiu2001SL}.  For problems on a sphere, SL methods have generally used latitude-longitude grids and spherical coordinate systems (e.g.~\cite{rasch1990SL,LaytonSpotz2003,TOLSTYKH2002180}), or other regular surface-based grids and local surface-based coordinate systems (e.g.~\cite{CascadeInterpolation,Lauritzen20101401,Carfora2007}), which can lead to \revfour{a loss of accuracy} because of singularities that arise in the mappings from the physical sphere to the surface-based coordinate systems.}

In this paper, we present three new high-order SL methods for transport on a sphere based on interpolation with global RBFs, local RBFs, and RBF-PU methods. These methods are mesh-free, allowing for scattered node discretizations, and are formulated entirely in Cartesian coordinates so as to avoid any surface based coordinate singularities.  The local RBF and RBF-PU methods also allow for a type of ``$p$-refinement'' for increasing the accuracy for a given fixed set of discretization nodes.  We demonstrate that the SL framework lends our new methods both accuracy and intrinsic stability, thereby eliminating the need for a hyperviscosity term. For the local RBF and RBF-PU methods, we propose using ``scale-free'' RBFs appended with spherical harmonics (an idea related to that of Flyer et al.\ \cite{FlyerPHS} for planar domains) to further reduce the number of tuning parameters (\emph{i.e.},\ the shape-parameter) and to bypass so-called error stagnation.  We compare and contrast all three methods using three standard test cases from the literature---solid-body rotation of a cosine bell from~\cite{Wil92} and deformational flow of two bells from~\cite{NairLauritzen2010}.  The focus of these comparisons is on the overall accuracy, dissipation and dispersion properties, mass conservation, and computational cost.  We find that the computational costs of these new methods are comparable to those of existing RBF collocation and finite-difference techniques. In particular, we find that the local RBF and RBF-PU methods are highly scalable to large node sets. We note that RBFs have previously been used for SL advection in~\cite{TELA:TELA0009}, but the focus there was on planar domains and global RBF methods.  Additionally, they have been used in a conservative SL advection method in~\cite{iske2002SL} for planar domains using Voronoi cells and local thin plate splines.  This is the first application of RBFs to SL transport on a sphere with accurate and scalable numerical methods and in a entirely mesh-free formulation.

We note that, while SL methods are applicable to a wide class of advection dominated problems, they have some limitations, the primary one being a lack of local conservation unless explicitly formulated in a conservative manner. However, the \revfour{large-time-step conservative SL methods are computationally expensive}~\cite{Zerroukat1JCP2010}. 

The paper is organized as follows. In the next section, we review the global RBF, local RBF and RBF-PU methods in the context of interpolation. In Section 3, we review the SL advection technique and discuss how to use use the three RBF methods within the SL framework in an efficient fashion. In Section 4, we compare and contrast the new SL methods on three standard test cases for transport on the unit sphere. Finally, we conclude with a summary of our results and future research directions in Section 5.
\section{Global, local, and partition of unity RBF interpolation on $\Sph$} \label{sec:rbf_review}
RBFs are a well-established method for interpolating/approximating data over a set of ``scattered'' nodes $X=\{\vx_j\}_{j=1}^N \subset\Omega\subseteq\mathbb{R}^d$.  The standard method uses linear combinations of shifts of a kernel $\phi :\Omega \times\Omega \to \mathbb{R}$ with the property that $\phi(\vx,\vy) := \phi(\|\vx-\vy\|)$ for $\vx,\vy\in\Omega$, where $\|\cdot\|$ is the Euclidean norm in $\mathbb{R}^d$.  Kernels with this property are referred to as \emph{radial kernels} or simply \emph{radial functions}.  In the case where $\Omega=\Sph$, these kernels are sometimes referred to as \emph{spherical basis functions} since $\phi(\|\vx-\vy\|) = \phi(2 - 2\vx^{T}\vy)$ when $\vx,\vy\in \Sph$, i.e. $\phi$ will only depend on the cosine of the angle between $\vx$ and $\vy$.  We will use the term RBFs and not make use of this simplification since our numerical schemes use extrinsic coordinates.

Below we review the three interpolation methods used in this study. We give cursory details on the first two methods, as they have appeared in many other places in the literature; see the recent book by Fornberg \& Flyer~\cite{FFBook} for further details and applications of these two methods.

%Since the method is meshfree, the nodes $X$ are not required to live to be on a grid or mesh and can be chosen however we wish for our application. In this study, we use node sets that provide near optimal resolution over $\Sph^2$.  Since $N=20$ nodes is the maximum number that can be exactly equally distributed over $\Sph^2$, one is resigned to using node sets that are only quasi-uniformly distributed over $\Sph^2$.  These node sets, which can be generated from a variety of algorithms~\cite{HardinSaff}, have the property that the average spacing between nodes, $h$,  roughly satisfies $h\sim 1/\sqrt{N}$.  In the results presented here we have used both minimum energy nodes~\cite{WomerSloanSpherePts,Go2003} and icosahedral nodes~\cite{}.  All of the node sets used here can be downloaded from~\cite{SpherePts}.

\subsection{Global RBF interpolation} \label{sec:global_rbfs}
%While our interest in this paper is interpolation on the surface of the unit sphere, $\Sph$, we introduce global RBF interpolation for the more general problem of interpolating on some subdomain of $\mathbb{R}^d$ since the direct method of constructing the interpolants is impervious to the domain~\cite{FuselierWright2012}.  At the end of this section, we comment on the specific case of interpolation on $\mathbb{S}^2$.
%Let $\Omega \subseteq \mathbb{R}^d$, and $\phi :\Omega \times\Omega \to \mathbb{R}$ be a kernel with the property $\phi(\vx,\vy) := \phi(\|\vx-\vy\|)$ for $\vx,\vy\in\Omega$, where $\|\cdot\|$ is the standard Euclidean norm in $\mathbb{R}^d$.  We refer to kernels with this property as \emph{radial kernels} or \emph{radial functions}. 
Given a set of nodes $X = \{\vx_k\}_{k = 1}^N \subset \Sph$ and a continuous target function $f:\Sph \to \mathbb{R}$ sampled at the nodes in $X$, the standard \emph{global} RBF interpolant to the data has the form
\begin{align}
s(\vx)= \sum_{k=1}^N c_k \phi(\|\vx - \vx_k\|).
\label{eq:srbfinterp}
\end{align}
The expansion coefficients $\{c_k\}_{k=1}^N$ are determined by enforcing $\lf.s\rt|_{X} = \lf.f\rt|_{X}$, which can be expressed by the following linear system:
\begin{align}
\underbrace{
\begin{bmatrix}
\phi(r_{1,1}) & \phi(r_{1,2}) & \hdots & \phi(r_{1,N}) \\
\phi(r_{2,1}) & \phi(r_{2,2}) & \hdots & \phi(r_{2,N})\\
\vdots & \vdots & \ddots & \vdots\\
\phi(r_{N,1}) & \phi(r_{N,2}) & \hdots & \phi(r_{N,N})
\end{bmatrix}}_{\ds A_X}
\underbrace{
\begin{bmatrix}
c_1 \\
c_2 \\
\vdots \\
c_N
\end{bmatrix}}_{\ds \vc_X}
=
\underbrace{
\begin{bmatrix}
f_1 \\
f_2 \\
\vdots \\
f_N
\end{bmatrix}}_{\ds \vf_X},
\label{eq:rbf_linsys}
\end{align}
where $r_{i,j} = ||\vx_i - \vx_j||$. If $\phi$ is, for example, a positive-definite radial kernel on $\mathbb{R}^3$, and all nodes in $X$ are distinct, then the matrix $A_{X}$ above is guaranteed to be positive definite, so that \eqref{eq:srbfinterp} is well-posed.  Examples of various choices for $\phi$, including relaxed conditions to guarantee the well-posedness of \eqref{eq:srbfinterp} can be found in~\cite[Ch.\ 3--12]{Fasshauer:2007}.  Results on the approximation properties of RBF interpolants on the sphere for target functions of various smoothness can be found in, for example,~\cite{HubbertMorton:2004,NarcSunWard:2007,JetterStocklerWard:1999}.  A result from~\cite{JetterStocklerWard:1999} is that for infinitely-smooth target functions, convergence rates that are faster than any polynomial order can be realized for various infinitely-smooth $\phi$.  Additionally, global RBF interpolants based on various infinitely-smooth kernels have been shown to converge to spherical harmonic interpolants as the kernels are scaled to become flat (the so-called flat limit)~\cite{FornbergPiret:2007}.  

One issue with global RBF interpolation is that the linear system for determining the interpolation coefficients~\eqref{eq:rbf_linsys} is dense, thereby leading to an $O(N^3)$ computational cost to solve it using a direct method. Another issue is that the matrices can become ill-conditioned.  While the RBF-QR method~\cite{FornbergPiret:2007} can be used to bypass the ill-conditioning associated with small shape parameters (i.e.\ flat radial kernels),  algorithms for bypassing the $O(N^3)$ computational cost \emph{and} ill-conditioning--without sacrificing high-order accuracy--have not yet been developed. The next two subsections discuss techniques for addressing these issues.  We note that an alternative to these methods, that is not pursued here, is to use compactly-supported kernels in a multilevel type framework~\cite{GiaSloanWendland:2010}.

\subsection{Local RBF interpolation}
%One issue with global RBF interpolation is that the linear system to solve for the interpolation coefficients~\eqref{eq:rbf_linsys} is dense so that the computational cost to solve it, when using a direct method, is $O(N^3)$.  Another issue is that, especially for infinitely smooth kernels, the matrices can become ill-conditioned without scaling the basis functions appropriately.  One approach to alleviate this issue is to use compactly supported kernels in a multilevel type framework~\cite{GiaSloanWendland:2010}.  

\begin{figure}[h]
\centering
\includegraphics[width=0.25\textwidth]{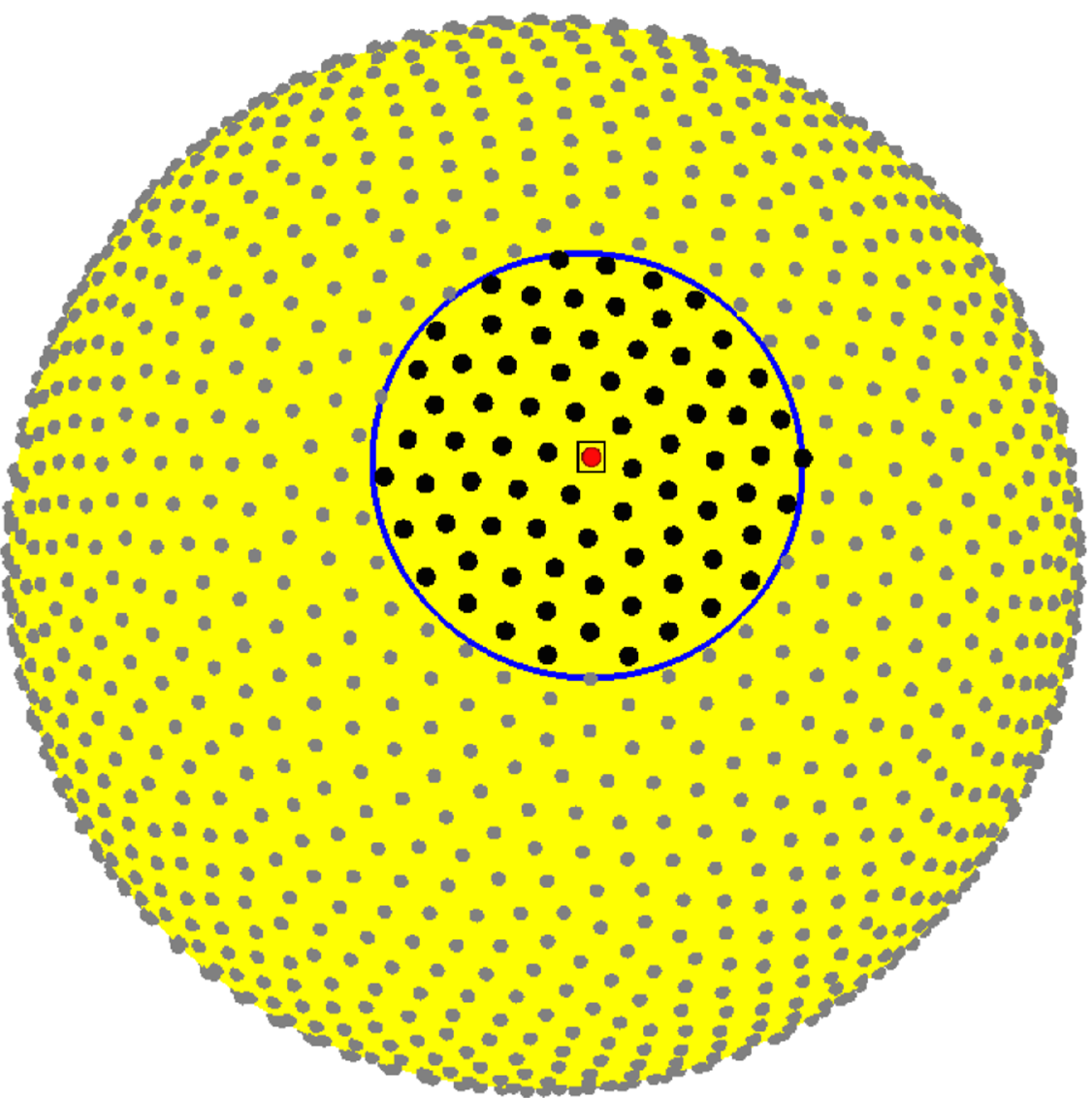}
\caption{Illustration of the nodes in an example stencil for the local RBF method on the sphere.  The global nodes $X$ are marked as small greyed balls, while the nodes that make up the stencil $X_k$ are marked as black balls, with the center point $\vx_k$ marked in red an with a black square around it.\label{fig:RBFFDIllustration}}
\end{figure}

In this method, RBF interpolants are used locally over a small collection of points surrounding each node in the set $X$ (this technique forms the foundation for RBF-FD methods~\cite{Wright200699}).  The method proceeds as follows.  For each node $\vx_k$, $k=1,\ldots,N$, in the set $X$, we select subsets of $X$ that consist of $\vx_k$ and its $n-1$ nearest neighbors, where $n << N$.  We refer to these subsets as \emph{stencils} and denote them by $X_k$, and refer to node $\vx_k$ as the \emph{center point} of the stencil $X_k$; see Figure \ref{fig:RBFFDIllustration} for an illustration.
%One approach for reducing the computational cost and potential ill-conditioning of the global RBF method is to use \emph{local} RBF interpolation.  
%The local method proceeds as follows.  For each selecting subsets (henceforth referred to as \emph{stencils}) $\{P_k\}_{k=1}^N$ of the set of nodes $X$, where each stencil $P_k$ consists of the $k^{th}$ node $\vx_k$ (henceforth referred to as \emph{center point} of the stencil) and its $n-1$ nearest neighbors, where $n << N$.  
%Additionally, an index set $\calI_k = \{\calI^k_1,\hdots,\calI^k_n\}$ that contains the \emph{global} indices of its nodes in the set $X$ is tracked with each stencil $P_k$.  
The nearest neighbors are typically determined in a preprocessing step using a data structure such as a kd-tree, which typically requires $O(N\log N)$ operations to construct.  On each stencil, an RBF interpolant is constructed to the nodes in $X_k$ and the corresponding function values.  It is then used to reconstruct the underlying target function globally at all points within some neighborhood of the center point of the stencil.  While this approach does not produce a globally smooth interpolant to the underlying target function, we find that it works well in the SL setting.  

For the local RBF interpolants, rather than use~\eqref{eq:srbfinterp} on each stencil, we follow an approach related to that proposed by Flyer et al.\ \cite{FlyerPHS} for RBF-FD methods to avoid stagnation errors for interpolation in $\mathbb{R}^d$.  Their approach augments the standard form of the interpolant~\eqref{eq:srbfinterp} with a linear combination of polynomials up to a degree corresponding to roughly $n/2$ polynomial basis functions, together with a set of ``moment conditions'' on the RBF coefficients.  A natural extension of their approach to the sphere is to augment the standard interpolant with \emph{spherical harmonics} up to a degree $L$ that grows with $n$ (see Section \ref{sec:kernel_choice}), since these are polynomials in $\mathbb{R}^3$ that are linearly independent when restricted the sphere.  The exact form of the augmented local RBF interpolants on each stencil $X_k$ that we use is then given by
\begin{align}
s_k(\vx)= \sum\limits_{j={\calI^k_1}}^{\calI^k_n} c_j^k \phi(\|\vx - \vx_j\|) + \sum\limits_{i=1}^{(L+1)^2} d_i^k p_i(\vx),
\label{eq:rbf_poly}
\end{align}
where $\calI^k_1,\ldots,\calI^k_n$ are the indices into the global node set $X$ of the nodes contained $X_k$ and $\{p_i\}_{i=1}^{(L+1)^2}$ are a basis for the space of spherical harmonics of degree $L$.  To determine the unknown coefficients, we impose that $s_k$ interpolates the data associated with the stencil $X_k$ and that the coefficients $c_j^k$ satisfy the following $(L+1)^2$ moment conditions:
\begin{align*}
\sum\limits_{j={\calI^k_1}}^{\calI^k_n} c_j^k p_i(\vx_j) = 0,\quad i=1,\ldots,(L+1)^2,
\end{align*}
which are the standard conditions imposed in the literature for polynomials in $\mathbb{R}^d$~\cite{FlyerPHS,Fasshauer:2007}.  The interpolation and moment conditions lead to the following linear system:
\begin{align}
\underbrace{
\begin{bmatrix}
A_{k} & P_{k} \\
P_{k}^T & 0
\end{bmatrix}}_{\ds \Asp_k}
\underbrace{
\begin{bmatrix}
\vc_k \\
\vd_k
\end{bmatrix}}_{\ds \csp_k}
=
\underbrace{
\begin{bmatrix}
\vf_k\\
{\bf 0}
\end{bmatrix}}_{\ds \fsp_k},
\label{eq:rbf_poly_linsys}
\end{align}
where $A_{k}$ is the RBF interpolation matrix for the nodes in stencil $X_k$ (see \eqref{eq:rbf_linsys}) and $P_k$ is a $n$-by-$(L+1)^2$ matrix with column $i$ containing entries $p_i(\vx_j)$, $j=\calI^k_1,\ldots,\calI^k_n$.  These local interpolants can be used to approximate functions to high-order algebraic accuracy determined by the stencil size $n$ and the degree of the spherical harmonics.  
%Such local interpolants when differentiated can be used to generate scattered-node finite difference (FD) formulas, known in the literature as RBF-FD~\cite{WrightFornberg}.  
Much like in~\cite{FlyerPHS}, we find that augmenting the standard RBF interpolant with spherical harmonics significantly enhances the accuracy of the local RBF interpolants. We will henceforth use tilde's above variables that correspond to the augmented RBF system \eqref{eq:rbf_poly_linsys}. 
%When used on a vector of function values, it means that the vector is padded with zeros at the end to match the number of columns in $\Psi_k$.

\subsection{RBF-PU interpolation}
The RBF-PU method is similar to the local RBF interpolation method. However, instead of forming and using different local interpolants for each node $\vx_k$,  local interpolants are constructed over a collection of \emph{patches} that cover the sphere.  These interpolants are then combined using compactly-supported weight functions on each patch that, all together, form a partition of unity. In this way, the method results in a smooth interpolant over the sphere, whereas the local approach does not. The RBF-PU method was first introduced for problems in the plane by Wendland~\cite{Wendland02fastevaluation} (see also~\cite[Ch.\ 29]{Fasshauer:2007}) and for interpolation problems on a sphere by Cavoretto \& De Rossi in~\cite{CaDeRo10}. It has since been extended to other domains and applications~\cite{Safdari-Vaighani2015,CaDeRoPe:2015,KevinThesis}. Specific details on the RBF-PU construction are provided below in the context of our application.

\begin{figure}[t]
\centering
\begin{tabular}{cc}
\includegraphics[width=0.3\textwidth]{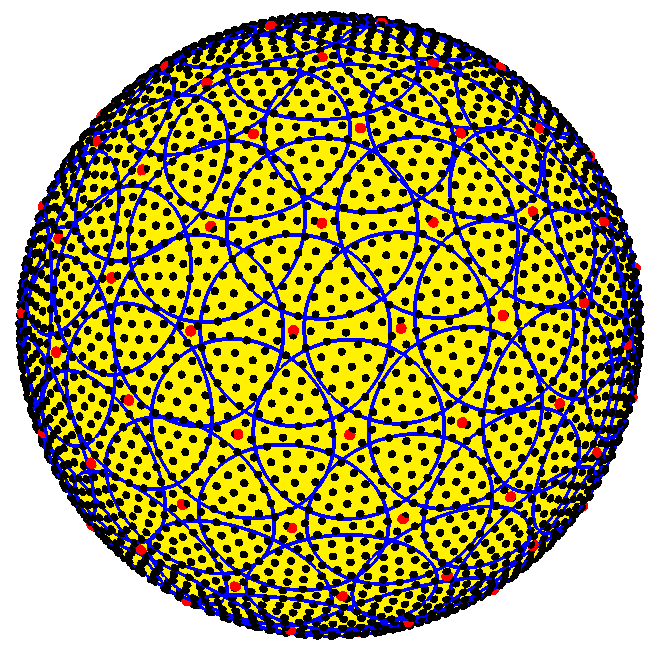} &
\raisebox{-0.04in}{\includegraphics[width=0.335\textwidth]{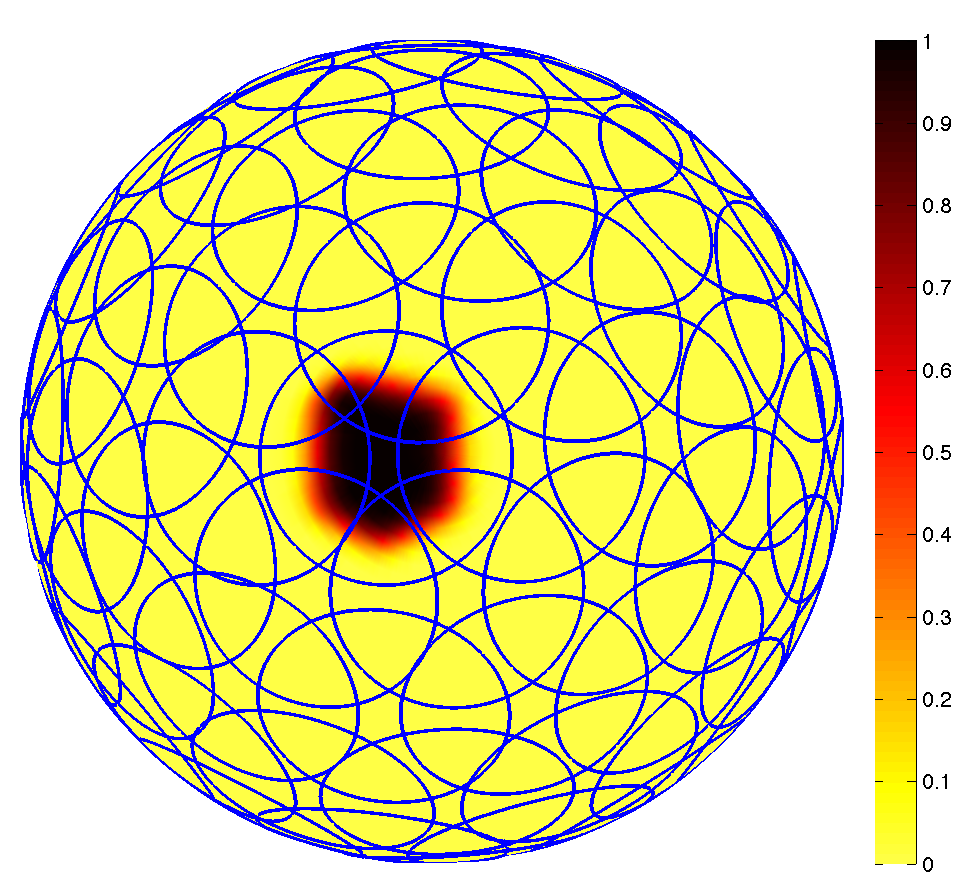}} \\
(a) & (b)
\end{tabular}
\caption{(a) Illustration of a typical collection of patches (outlined in blue) for a quasi-uniformly distributed node set of $N=4096$ nodes (marked in black) on $\Sph$.  The centers of the patches are marked in red and there are roughly 100 nodes per patch. (b) Illustration of one of the PU weight functions \eqref{eq:pu_func} for the patches from part (a).\label{fig:patches}}
\end{figure}

Let $\Omega_1,\ldots,\Omega_M$ be a collection of distinct spherical caps on $\Sph$ with the properties that
the set of all patches provides an open covering of  $\Sph$, i.e.\ $\ds \cup_{\ell=1}^{M} \Omega_\ell = \Sph$, and each $\Omega_\ell$ contains at least one node from $X$.
These caps are the patches in the PU method.  Figure \ref{fig:patches} (a) illustrates a typical collection of patches for a quasi-uniformly distributed node set $X$ on $\Sph$, which is what we use in this study.  Let $\vom_\ell\in\Sph$, $\ell=1,\ldots,M$, denote the center of patch $\Omega_\ell$ and $R_\ell > 0$ denote the radius of the patch, measured as the Euclidean distance from $\vom_\ell$.  For each patch, we define the following compactly-supported weight function:
\begin{align}
\varphi_\ell(\vx) = \varphi\lp\frac{\|\vx-\vom_\ell\|}{R_\ell}\rp,
\label{eq:wght_func}
\end{align}
where $\varphi:\mathbb{R}^{+}\rightarrow\mathbb{R}$ has compact support over the interval $[0,1)$ .  In this study, we use the radially-symmetric cubic B-spline
\begin{align}
\varphi(r) =
\begin{cases}
\frac{2}{3} + 4 \lp r -1 \rp r^2 & \text{if $0\leq r < \frac{1}{2}$,} \\
-\frac{4}{3} {\lp r - 1 \rp}^3 & \text{if $\frac{1}{2} \leq r < 1$,} \\
0 & \text{if $1 \leq r$}.
\end{cases}
\label{eq:bspline}
\end{align}
Using \eqref{eq:wght_func}, we define the PU weight functions for the collection of patches $\Omega_1,\ldots,\Omega_M$ as
\begin{align}
w_\ell(\vx) = \dfrac{\varphi_\ell(\vx)}{\ds \sum_{j=1}^{M} \varphi_j(\vx)},\; \ell=1,\ldots,M.
\label{eq:pu_func}
\end{align}
Note that each $w_\ell$ is only supported over $\Omega_\ell$ and that the summation on the bottom only involves terms that are non-zero over patch $\Omega_\ell$, which should be much smaller than $M$.  Figure \ref{fig:patches} (b) displays one of these weight functions for a quasi-uniform distribution of overlapping patches.

We use these PU weight functions to define the RBF-PU interpolant as follows.  Let $X_\ell$ denote the nodes from $X$ that belong to patch $\Omega_\ell$ and let $s_{\ell}$ denote the RBF interpolant of the form \eqref{eq:rbf_poly} to the target function $f$ over $X_\ell$.  Then the RBF-PU interpolant is given by
\begin{align}
s(\vx) = \sum_{\ell=1}^{M} w_\ell(\vx)s_{\ell}(\vx).
\label{eq:PU_interp}
\end{align}
Since $\{w_\ell\}_{\ell=1}^M$ form a partition of unity, it is straightforward to show that $s$ interpolates the target function at all nodes in $X$.  The interpolant is also smooth over the entire sphere.  Here we are using the augmented RBF interpolants \eqref{eq:rbf_poly} over each patch (which differs from~\cite{CaDeRo10}).   We have found that this gives better accuracy than the standard interpolant \eqref{eq:srbfinterp}.  In Section \ref{sec:RBF-PU-SLA}, we discuss our strategy for choosing the patches.

\subsection{Choosing the radial kernels}\label{sec:kernel_choice}

In the past, infinitely-smooth radial kernels, such as the inverse multiquadric (IMQ) kernel $\phi(r) =(1 + (\ep r)^2)^{-1/2}$, have been primarily used for RBF approximations to advection-dominated PDEs on a sphere, e.g.\ \cite{FoL11,FlyerLehto2012,FlyerWright:2007,FornbergPiret:2008, FlyerWright:2009, WFY}. When used with global RBFs, this has been demonstrated to result in convergence rates higher than any polynomial order for smooth solutions.  These kernels have also been shown to be better than finitely-smooth kernels even when the solutions are not smooth~\cite{FornbergPiret:2008}.  However, the advantages of using infinitely-smooth kernels over finitely-smooth ones with RBF-FD was recently called into question~\cite{FlyerPHS,FlyerNS}. Indeed, recent work on Euclidean geometries~\cite{FlyerNS,Shankar2016} has demonstrated that it is more beneficial to use the finitely-smooth polyharmonic spline (PHS) kernels in the context of RBF-FD methods, augmented with a moderate-degree polynomial. These kernels take the form $\phi(r) = r^{2k} \log r$, $k\geq 1$, or $\phi(r) = r^{2k+1}$, $k\geq  0$.  The choice between these two kernels is often made depending on the dimension of the interpolation problem for theoretical reasons,  while the choice of $k$ controls the smoothness of the kernels, and hence their accuracy~\cite{Fasshauer:2007}.
%Furthermore, in~\cite{FlyerPHS}, $k$ was shown to have little to no impact on approximation rates. 

In this work, we use the IMQ kernel for the global RBF method in order to benefit from the higher-than-polynomial convergence rates for smooth problems.   
% and set the shape parameter is set to $\ep=4.9$, which corresponds to a condition number of $10^{14}$ on the finest node set used for global RBFs.  
For the local RBF and RBF-PU methods we follow the recent work of~\revtwo{\cite{FlyerPHS,FlyerNS,Shankar2016,BayonaEtAl2017}} and use the PHS kernel $\phi(r) = r^{2k+1}$ in conjunction with spherical harmonics, as given in \eqref{eq:rbf_poly}.  We choose the degree of the spherical harmonics as
\begin{align*}
L=\left\lfloor\frac12(\sqrt{n} -1)\right\rfloor,
\end{align*}
and the order of the PHS as $k=L$.  For the many test problems we considered, these parameters produced excellent results.  Choosing $k=L$ also coincides with theoretical results on solvability of the PHS interpolation problem~\cite{Fasshauer:2007}. In contrast, recent work~\revtwo{\cite{FlyerPHS,FlyerNS,Shankar2016,BayonaEtAl2017}} has used polynomials of a much higher-degree than that required to prove solvability. For problems posed in general domains in $\mathbb{R}^2$ and $\mathbb{R}^3$, it was demonstrated in these studies (see, for example, Figures 7 and 8 from~\cite{FlyerPHS}) that using this approach for doing local interpolations/approximating derivatives one typically obtains better accuracy.  In our extensive experiments, we did not find this to carry over to the case for either the local or RBF-PU methods when the domain is $\mathbb{S}^2$ and spherical harmonics are used.  %It is possible that using the approach developed
In~\cite{ReegerFornberg2016}, the authors use a different approach for constructing local interpolants over the sphere (for doing quadrature), where the local stencil (or patch) is projected onto the tangent plane of $\mathbb{S}^2$ at the center of the stencil and the interpolations are done in the plane.  This approach, which we did not test here, may give similar results for the sphere as the observations from~\cite{FlyerPHS,FlyerNS,Shankar2016,BayonaEtAl2017}.

%An added benefit of using the PHS kernel and the formulas for $L$ and $k$ is that a shape parameter does not have to be selected.  
%We also tried the kernel  $\phi(r) = r^{2k} \log r$ with these same parameters and obtained similar results, but opted for the simpler PHS kernel.

\section{RBFs for SL advection on a sphere}
\label{sec:rbfsl}
As discussed in the introduction, SL advection is a technique for solving the transport equation \eqref{eq:advection}.
%  For the advection of a  scalar quantity $q$ on $\Sph$ in an incompressible velocity field $\vu$, the transport equation is given by 
%\begin{align}
%%\frac{\partial q}{\partial t} + \vu \cdot \nabla_{\Sph} q = 0,\quad \text{or}\quad \frac{Dq}{Dt} &= 0, \label{eq:advection}
%\frac{Dq}{Dt} = 0,\quad \frac{D}{Dt}:=\frac{\partial }{\partial t} + \vu \cdot \nabla_{\Sph}, \label{eq:advection}
%\end{align}
%where $\vu$ is tangent to $\Sph$ and $\nabla_{\Sph}$ is the surface gradient operator on $\Sph$.  
The primary idea of SL advection is to use the Lagrangian frame to find upwind directions for which $q$ moves, and use an Eulerian set of nodes to perform interpolations of $q$. \revthree{The upwinding step ensures that the numerical domain of dependence matches the physical domain of dependence}. The presence of an Eulerian frame avoids spatial resolution issues and allows for high-order accuracy, a feature that is difficult to achieve in the Lagrangian frame without remeshing, or redistributing nodes.

\revthree{The upwinding is done as follows: assume that a Lagrangian parcel arrived at each of the Eulerian nodes carrying with it some amount of the scalar $q$. The amount of $q$ at this node must have been brought forward from the \emph{departure point} of the parcel. Thus, we find the departure point by tracing the parcel along the flow field $\vu$ and determine $q$ by interpolating the scalar field to the departure point. This procedure is laid out explicitly in Algorithm \ref{alg:sl_alg}. More details on the general SL advection procedure can be found, for example, in~\cite{FalconeFerrettiBook}. In the following subsections, we discuss the different aspects of the algorithm as it pertains to advection on $\mathbb{S}^2$.}

\begin{algorithm}
\caption{Semi-Lagrangian advection on $\Sph$}
\label{alg:sl_alg}
\begin{algorithmic}
        \State \textbf{Input:} Velocity field $\vu(\vx,t)$ tangent to $\Sph$; 
        initial scalar field $q(\vx,0)$; 
        $X = \{\vx_j\}_{j=1}^N$, $X \subset \mathbb{S}^2$;
	final time $t_f$; time-step $\Delta t$.
	\State Set $\vq_{X}^{0} = \{q(\vx_j,0)\}_{j=1}^N$, $t=0$, and $m=0$.
	\While{$t \leq t_f$}
%	  \State \textbf{Interpolate} $Q^m$ and store the interpolation coefficients. \comment{What does the superscript $m$ denote?}	
		\State  For $j=1,\ldots,N$, trace back $\vx_j$ to time $t$ to find departure point $\vxi^m_j$.
		\State Interpolate $\vq^{m}_{X}$ to $\Xi^{m}=\{\vxi_j^{m}\}_{j=1}^N$ to obtain $\vq_{\Xi}^m$.	
		\State Set $\vq^{m+1}_X =  \vq_{\Xi}^m$, $m=m+1$, and $t =  m\Delta t$.
	\EndWhile
\end{algorithmic}
\end{algorithm}

\subsection{Eulerian node sets}\label{sec:eulerian_nodes}
Central to all SL advection schemes is a fixed set of Eulerian nodes $X=\{\vx_j\}_{j=1}^N$ over the given domain that is used to interpolate the advected quantity to the \revthree{Lagrangian parcels}. The approximate solution to the PDE will ultimately be computed only at the set of Eulerian nodes $X$.  Since the interpolation schemes used in this study are all based on RBFs, the Eulerian nodes are not required to live on a grid or mesh, and can be chosen freely for our application. In this study, we use node sets that provide near optimal resolution over $\Sph$.  Since $N=20$ nodes is the maximum number that can be exactly equally distributed over $\Sph$, one is resigned to using node sets that are only quasi-uniformly distributed over $\Sph$.  These node sets, which can be generated from a variety of algorithms~\cite{HardinSaff:2004}, have the property that the average spacing between nodes, $h$, satisfies $h\sim 1/\sqrt{N}$.  In the results presented here we have used maximum determinant nodes~\cite{WomersleySloan:2001} and icosahedral nodes~\cite{icos}.  While the latter of these sets forms a natural grid, we do not use this fact in our algorithms.  All of these node sets, including many others, are available in the \texttt{SpherePts} software package~\cite{SpherePts}.

\subsection{Trajectory reconstruction}
The process of tracing \revthree{Lagrangian parcels} back along the velocity field $\vu$ is known as \emph{trajectory reconstruction}. %In the Lagrangian frame, if $\vxi$ is the position of the particle, we find the departure point $\vxi^d$ by solving the ODE
Let $\vxi_j(t)\in\Sph$ denote the position of a \revthree{Lagrangian parcel} as a function of time such that $\vxi_j(t_{m+1}) = \vx_j$, for some time $t=t_{m+1} > 0$. Here $\vx_j$ is the $j^{\text{th}}$ node in Eulerian node set $X=\{\vx_i\}_{i=1}^{N}\subset\Sph$.  Then we find the departure point of this \revthree{parcel} at time $t=t_{m}=t_{m+1}-\Delta t$ by solving the simple ODE
\begin{align}
\frac{d\vxi_j}{dt} = \vu,\; \vxi_j(t_{m+1}) = \vx_j,
\label{eq:sl_ode}
\end{align}
\emph{backward} in time over the interval $[t_m,t_{m+1}]$.  Here $\Delta t$ defines the interval over which we trace back the \revthree{parcel} from $\vx_j$, and it also defines the time-step for solving the advection equation \eqref{eq:advection}.  

When the velocity field $\vu$ is known for all time and space the trajectory ODE \eqref{eq:sl_ode} can be solved in a straightforward manner using a variety of techniques.  However, since we use the Cartesian coordinate representation of the velocity field $\vu$ in the trajectory reconstruction, the \emph{Lagrangian parcels} may not reside on the sphere at time $t=t_{m}$.  Similar to~\cite{BOSLER2017639}, which also uses Cartesian coordinates (albeit in a fully Lagrangian scheme), we have found that the \revthree{parcels} remain very close to the sphere when a high-order integrator is used.  In our experiments, we found that both the standard fourth-order Runge-Kutta (RK4) method and Fehlberg's fifth-order Runge-Kutta (RK5) scheme~\cite{fehlberg68}, worked well.  We ultimately decided on the RK5 scheme since we found that it allowed for larger time-steps than RK4.  We also found that orthogonally projecting the \revthree{Lagrangian parcels} exactly back to the surface of the sphere at each stage of the RK scheme improved the overall accuracy of the method. 

The trajectory reconstruction procedure is repeated for $j=1,\ldots,N$ giving a set of $N$ departure points at time $t=t_m$, which we denote by $\Xi^{m}=\{\vxi_j^{m}\}_{j=1}^N$ and refer to as the \emph{Lagrangian point set}.  It is at these points that we need to compute $q$.  However, $q$ is only (approximately) known at the nodes in $X$ and in general $X\neq\Xi^{m}$, for any $m>0$. Thus, we need to interpolate $q$ from $X$ to $\Xi^m$.

\revthree{
\begin{rmk}
SL methods are not subject to the advective CFL condition for the purpose of numerical stability (on tracer transport problems). So, in principle, if the velocity field is known for all time and space, then one can use integrators such as MATLAB's ode45 method to take one time-step over the whole simulation time.  However, for problems that feature forcing terms or that require solving for the velocity field during the simulation~\cite{xiu2001SL,LaytonSpotz2003,TOLSTYKH2002180,Weizhu2003}, one may be forced to take smaller intermediate time-steps.  Additionally, for complex flows with significant variability in small scales, time-step restrictions arise from conditions on the inverse of the flow Jacobian in order to produce meaningful results~\cite{CossetteJCP2014}.  This can result in time-step restrictions similar to the advective CFL condition.  In problems where the flow fields are given and are smooth and slowly varying (as considered in this study), it is common to use a fixed time-step and (high-order) time integrator such that spatial errors dominate (see, for example,~\cite{LauritzenEtAl2014}). This is the approach we follow.
\end{rmk}
}

\subsection{Interpolation}

Let $\vq_X^{m}$ denote the vector containing the approximate solution to \eqref{eq:advection} at the Eulerian node set $X$ (with the same ordering as the nodes in $X$) at time $t=t_{m}$.  Below we discuss the details of interpolating $\vq_X^{m}$ to the Lagrangian point set $\Xi^{m}$ using the three different RBF interpolation methods presented in Section \ref{sec:rbf_review}.  These methods can all be used in the interpolation phase of Algorithm 1.
%Once the Lagrangian node set $X_d$ is found, we only need to interpolate the numerical solution $Q$ \comment{What is $Q$? VS: The numerical analog to $q$.  GW: We need to say this somewhere.} to $X_d$ to obtain the advected quantity $Q_d$. There are many possible ways to do so. In this section, we discuss three high-order methods based on the types of RBF interpolation outlined in Section \ref{sec:rbf_review}. The explicit algorithms for these three methods are detailed in \ref{sec:appendix0}. 

\subsubsection{Global RBF interpolation}

The most straightforward interpolation scheme is built on the global RBF interpolant \eqref{eq:srbfinterp}.  To compute this interpolant, we first solve the linear system \eqref{eq:rbf_linsys}, with the right hand side set equal to $\vq_X^{m}$, to determine the interpolation coefficients $\vc_X$.  We then evaluate the interpolant at the points in $\Xi^m$, resulting in a vector of approximate values, which we denote by $\vq_{\Xi}^m$.   Finally, we set $\vq_{X}^{m+1}=\vq_{\Xi}^m$ to complete a time step.

The coefficient matrix $A_X$ in \eqref{eq:rbf_linsys} does not change throughout the simulation (since the nodes $X$ are fixed for all time), thus allowing us to factorize $A_X$ once as a preprocessing step.  We use the inverse multiquadric (IMQ) radial kernel, $\phi(r) = (1+(\ep r)^2)^{-1/2}$, in this study, which is positive-definite, so that the matrix $A_X$ can be factorized using Cholesky factorization.  Since $A_X$ is dense, this requires $O(N^3)$ operations initially.  Each time step then requires $O(N^2)$ operations to solve \eqref{eq:rbf_linsys} for each $\vq_X^{m}$.  The evaluation of \eqref{eq:srbfinterp} at $\Xi^m$ requires an additional $O(N^2)$ computations.  While this cost becomes prohibitive for large $N$, the global RBF method remains competitive for smooth solutions because of its high accuracy---faster than any algebraic order when the target function is sufficiently smooth~\cite{JetterStocklerWard:1999}---as illustrated in the numerical results section.

\subsubsection{Local RBF interpolation}

The local RBF interpolation scheme is more complicated to implement than the global case. The first step of this scheme requires determining the $n$-point local stencil from $X$ that will be used for interpolating to each of the points in $\Xi^m$.  The approach we take is to determine the nearest neighbor in $X$ for each $\vxi_j^m\in\Xi^m$, $j=1,\ldots,N$.  The nearest neighbor determines the center point for the stencil that will be used for the interpolation.  Letting $\mathcal{K}_j$, denote the index of the nearest neighbor in $X$ for $\vxi_j$, the interpolation stencil is then denoted as $X_{\mathcal{K}_j}$.  After these stencils have been determined, we perform the interpolations to the points in $\Xi^m$.  These interpolations are done, for each $\vxi_j^m$, $j=1,\ldots,N$, using \eqref{eq:rbf_poly}, with $k=\mathcal{K}_j$, which results in the vector of approximate values $\vq_{\Xi}^m$.  As in the global case, we then set $\vq_{X}^{m+1}=\vq_{\Xi}^m$ to complete the time step. 

The computation of the local interpolant for $\vxi_j$ requires solving the linear system \eqref{eq:rbf_poly_linsys} with the right hand side set equal to $\tilde{\vq}_{X_{\mathcal{K}_j}}^m$.  While it may be possible for some points in $\Xi^m$ to share the same stencil, in general there will be $O(N)$ of these linear systems to solve.  The matrices for these linear systems can all be pre-computed and factorized as a pre-processing step since the points in $X$ do not change and all the stencils can be determined \emph{prior} to the simulation.  Unlike the global case, we use $LU$-factorizations of the matrices, since they are no longer positive-definite due to the inclusion of the polynomial terms and the use of the PHS kernel. The method initially requires $O(n^3 N)$ operations to factorize the linear systems for each stencil, then the cost is $O(n^2 N)$ operations to determine the interpolation coefficients of the local interpolants for all the points in $\Xi^m$ and to evaluate the local interpolants to obtain $\vq_{\Xi}^m$.  To determine the nearest neighbors for each $\vxi_j^m\in\Xi^m$ in $X$, we use a kd-tree of the nodes $X$.  This requires $O(N\log N)$ operations to construct initially, and $O(\log N)$ operations for each nearest neighbor search on average. Since $n$ is chosen independent of $N$ and typically $n << N$, the method has an asymptotic computational cost of $O(N\log N)$ per time-step, with a one-time initial cost of $O(N\log N)$.

Even when including the cost of kd-tree searches for nearest neigbhors, the local method is much more computationally efficient than the global RBF method. With this approach, it is possible to obtain high-order methods with algebraic convergence rates that depend on the stencil size $n$ and the degree of the appended spherical harmonics.

%Further, to obtain close to uniform accuracy over all stencils, we ensure that the condition numbers of the unaugmented RBF interpolation matrices $A_{P_k}$ are all on the same order of magnitude. To do this, we solve for the shape parameter $\ep_k$ that produces the target condition number for each stencil as a preprocessing step, and store this as well. \comment{We need to explain or give a reference to why this gives stencils with approximately uniform accuracy.  Note that we have said nothing yet about the shape parameter. VS: I am not sure how to do this. What would you suggest?}. Once the departure point is found, we interpolate $Q^m_{P_k}$ on the appropriate stencil by using the triangular factors in backsolves. This incurs a total cost of $O(n^2 N)$ per time-step, compared to the $O(N^2)$ cost of the global RBF method. Even including the cost of kd-tree searches for nearest neigbhors, this is much cheaper than the global RBF case. With this approach, it is possible to obtain high-order methods with algebraic convergence rates that depend on the stencil size $n$ and the degree of the appended spherical harmonics.

\subsubsection{RBF-PU interpolation}\label{sec:RBF-PU-SLA}

%The implementation RBF-PU interpolation scheme has similarities to both the global and local methods.  The 
%The first step of the RBF-PU interpolation scheme 
The first step of the RBF-PU interpolation scheme is to construct the local interpolants $s_{\ell}$ in \eqref{eq:PU_interp} over the patches (spherical caps) $\Omega_\ell$, $\ell=1,\ldots,M$.  This requires solving $M$ linear systems of the form \eqref{eq:rbf_poly_linsys}, where the node sets and right hand sides used are determined by which nodes from $X$ are included in patch $\Omega_\ell$.  Once all the local interpolants are constructed, they are combined into the globally-smooth interpolant \eqref{eq:PU_interp}, which is then evaluated at all the points in $\Xi^{m}$.  Note that for a given $\vxi_j^m\in\Xi^m$, the sum in \eqref{eq:PU_interp} only needs to be taken over terms corresponding to the patches $\Omega_{\ell}$ containing $\vxi_j$.  As with the global and local methods, the result is a vector of approximate values $\vq_{\Xi}^m$, which are assigned to $\vq_{X}^{m+1}$ to complete the time step.

The computational complexity of the RBF-PU scheme is determined by how the patches $\Omega_{\ell}$, $\ell=1,\ldots,M$, are distributed; this also directly effects the accuracy of the interpolants.  Since the node sets $X$ are assumed to be quasi-uniformly distributed, it makes sense to distribute the patches in a quasi-uniform manner to control the computational cost.  Distributing the patches in a quasi-uniform manner is equivalent to distributing the set of patch centers $\{\vom_{\ell}\}_{\ell=1}^M$ in a quasi-uniform manner.  For these sets we use minimum energy (ME) points, which are computed by arranging the points in the set such that the Reisz energy (with a power of 2) of the set on $\Sph$ attains a minimum~\cite{HardinSaff:2004}.  We use the pre-computed quasi-ME point sets from~\cite{SpherePts}, which are available for $2\leq M \leq 5000$.

What remains to determine the distribution of patches is $M$, the total number of patches (centers), and $R_{\ell}$, $\ell=1,\ldots,M$, the radii of the patches.  We determine these quantities by setting 1) the number nodes, $n$, that each patch is to approximately contain and 2) the average number of patches, $a$, that a given node is to belong to.  The first of these controls the computational cost for each patch and effects the accuracy of the local interpolants $s_{\ell}$, as each interpolant will be based on approximately $n$ nodes.  The second value $a$ controls the computational cost of evaluating the global RBF-PU interpolant \eqref{eq:PU_interp}, as it determines how much overlap there is amongst the patches.  This value also controls the locality of the interpolant.

%The patches are determined by their centers and radii.  For the centers, we choose quasi-uniformly distributed points We use two criteria to determine the distribution of the patches:  1) the number nodes, $n$, that each patch is to approximately contain and 2) the average number of patches that a given node is to belong to.   The first criteria controls the computational cost for each patch and effects the accuracy of the local interpolants $s_{\ell}$, as each interpolant will be based on approximately $n$ nodes.  The second criteria controls the computational cost of evaluating the global RBF-PU interpolant \eqref{eq:PU_interp}, as it determines how much overlap there is amongst the patches.  Since the node sets $X$ are assumed to be quasi-uniformly distributed it makes sense to distribute the patches in a quasi-uniform manner.  We thus need to determine the number of patches we npatches are determined by their centers and radii

We use $n$ to determine the radii of the patches as follows.  If $X$ contains $N$ quasi-uniformly distributed nodes and there are to be approximately $n$ nodes per patch, then the area per node over the entire sphere should approximately equal the area per node over patch $\Omega_\ell$, i.e.\ $4\pi/N \approx \pi R_{\ell}^2/n$, where $R_{\ell}$ is the radius of $\Omega_{\ell}$.  This gives the estimate $R_{\ell} \approx 2 \sqrt{n/N}$.
Since the centers of the patches are also assumed to be quasi-uniformly distributed, the radii can all be chosen in the same way.  We thus set one radius,
\begin{align}
R=2\sqrt{n/N},\label{eq:rad_patches}
\end{align}
for every patch.  We use $a$ to determine the number of patches $M$ as follows.  Since the patch centers $\{\vom_{\ell}\}_{\ell=1}^M$ are quasi-uniformly distributed, the same area arguments as above can be used to arrive at the estimate $4\pi/M \approx \pi R^2/a$.  Solving for $M$ in this equation and using \eqref{eq:rad_patches}, we obtain the following approximation (noting that $M$ should be an integer)
\begin{align}
M = \lceil a N/ n \rceil.
\label{eq:num_patches}
\end{align}
Figure \ref{fig:patches} illustrates the patches for $N=4096$, $n=100$, and $a=2.5$.

Having defined these quantities, we can now estimate the computational cost of the RBF-PU interpolation scheme (per time-step).  Similar to the local RBF scheme, the matrices in the $M$ linear systems \eqref{eq:rbf_poly_linsys} for determining $s_{\ell}$, $\ell=1,\ldots,M$ do not change per time step and can thus be $LU$ decomposed as a preprocessing step.  Since each linear system contains approximately $n$ nodes and there are $M$ linear systems, this cost is $O(n^3 M)$, which can be estimated as $O(a n^2 N)$ from \eqref{eq:num_patches}.  Using the $LU$ factorizations, the $M$ linear systems \eqref{eq:rbf_poly_linsys} can be solved in $O(a n N)$ operations per time step.  The evaluation of any of the interpolants $s_{\ell}$ in \eqref{eq:PU_interp} at a point on $\Sph$ requires $O(n^2)$ operations.  Each point of $\Xi^m$ will, on average, belong to $a$ patches, so that the sum in \eqref{eq:PU_interp} will, on average, only involve $a$ non-zero terms.  Thus, the computational cost of evaluating the RBF-PU interpolant \eqref{eq:PU_interp} at $\Xi^m$ is $O(a n^2 N)$.  To determine which patches a given point in $\Xi^m$ belongs to, we use a kd-tree of the patch centers, which only has to be constructed once at a cost of $O(M\log M)=O(a N/n \log (a N/n))$ operations and can then be searched each time-step at a cost of $O(N \log (a N/n))$.  Since $n$ and $a$ are chosen independent of $N$ and typically $a << n << N$, the method has an asymptotic computational cost of $O(N \log N)$ per time-step, with a one time initial cost of $O(N\log N)$.
%\begin{rmk}
%An important special case of SL advection is when the problem is restricted to the real line and forward Euler for the trajectory calculation.  If we restrict $\Delta t \sim h$ where $h \propto \frac{1}{N}$ is the node spacing on the real line, and use linear interpolation in conjunction with a forward Euler backtrace to find $\vq_{\Xi}^m$, the SL advection scheme reduces to the classical first-order upwind scheme. Even if one lifts the CFL restriction on $\Delta t$ in this setting, the SL advection scheme with linear interpolation is unconditionally stable~\cite{FalconeFerrettiBook}. Unfortunately, it is unclear how to prove stability for either high-order polynomial interpolation or RBF interpolation. However, in practice, all of our RBF methods are stable for values of $\Delta t$ much larger than the CFL condition dictates. 
%\end{rmk}

%\input{Measures}
\section{Results}
\label{sec:results}

We investigate various properties of the three new methods proposed in this article on three standard test cases for transport on the sphere from the literature: solid-body rotation of a cosine bell from~\cite{Wil92} and deformational flow with two different initial conditions from~\cite{NairLauritzen2010}. We test our methods using quasi-uniformly distributed node sets $X$ of various sizes $N$.  For the global RBF method, we use maximum determinant node sets~\cite{WomersleySloan:2001} of sizes $N=3136$, 4096, 5041, 6084, 7744, 9025, 10000, 11881, 13689, and 15129.  For the local RBF and RBF-PU methods, we use equidistributed icosahedral nodes~\cite{MICHAELS2017} of sizes $N=2562$, 5762, 10242, 23042, 40962, and 92162.  All these node sets are available from~\cite{SpherePts}.  In the results that follow, we estimate the convergence of our method verses $\sqrt{N}$ since this is roughly inversely proportional to the average spacing between the nodes (see Section \ref{sec:eulerian_nodes} for a discussion).  For the local RBF and RBF-PU methods we present results for $n=17$, $31$, $49$, and $84$.  These are common values used in RBF-FD methods for advective PDEs on the sphere~\cite{FoL11,FlyerLehto2012}, and thus provide easy comparisons with those methods.  

The main focus of the investigation is on accuracy of the methods, which we measure using the relative $\ell_2$ and $\ell_{\infty}$ norms.  We also give results on the dissipation and dispersion errors using a visual inspection of the errors and the \emph{a posteriori} quantitative measures derived in~\cite{Takacs:1985} (see also~\cite{RestelliBonaventura:2006}).  These measures are derived from a decomposition of the mean square error of the numerical solutions and are given as follows:
\begin{align}
\text{Dissipation error:}\; & \left[\sigma(q)-\sigma(q_{X})\right]^2 + \left[\bar{q}-\bar{q_{X}}\right]^2, \label{eq:dissipation}\\
\text{Dispersion error:}\; & 2\left[\sigma(q)\sigma(q_{X}) - \frac{1}{4\pi} \int_{\Sph} (q-\bar{q})(q_X-\bar{q_X}) dS\right] \label{eq:dispersion},
\end{align}
where $q$ is the exact solution and $q_X$ is the approximate solution, and bars on these variables denote the mean, while $\sigma$ denotes their standard deviation over $\Sph$.  In the results, we divide these errors by the mean square error of the solution to give a relative measure.  Finally, we present results on the conservation properties of the methods by plotting the absolute errors in the total mass: $\lf | \frac{1}{4\pi}\int_{\Sph} (q - q_X) dS\rt|$.   All of the above quantities are computed using discrete approximations to the continuous operators defining them on $\Sph$, which requires computing numerical approximations to surface integrals over $\Sph$ using the nodes $X$.  To compute these surface integrals, we use the sixth-order kernel-based meshfree quadrature method from~\cite{FHNWW2013b}.  

While not reported here, we also experimented with other quasi-uniform node sets for the global RBF method, specifically the the ME and spherical $t$-design nodes.  In~\cite{FornbergPiret:2008} (see also~\cite{FornbergMartel2014}), it was found that the MD vastly outperform the ME nodes when using flat RBFs (close to the spherical harmonic limit~\cite{FornbergPiret:2007}).  For this study, we used shape parameters, $\ep$, that are not in the flat regime and found the the resulting errors for all the test cases considered below were roughly comparable amongst the different types of node sets, with the ME nodes giving slightly smaller errors overall.  We chose the MD nodes since there are many more of these node sets available for the sizes we were interested in testing.  Finally, it should be noted that we did not try and optimize the value of $\ep$ for each test problem and each node set, so there is potential room for improvement in some of the results presented for the global RBF method.  However, optimization of $\ep$ is rarely possible in real applications.

\subsection{Solid body rotation of a cosine bell}
As a first test problem, we consider the standard Test Case 1 from Williamson et~al.~\cite{Wil92}. The components of the steady velocity field for this test case in spherical coordinates ($-\pi \leq \lambda \leq \pi$, $-\pi/2 \leq \theta \leq \pi/2$) are given by
\begin{align*}
u(\lambda, \theta) =  \sin(\theta)\sin(\lambda)\sin(\alpha)-\cos(\theta)\cos(\alpha),\quad v(\lambda, \theta) =  \cos(\lambda)\sin(\alpha).
\end{align*}
This velocity field results in solid body rotation at an angle of $\alpha$ with respect to the equator.  In all of our tests, we use one of the standard choices of $\alpha = \pi/2$, which corresponds to advection of the initial condition directly over the north and south poles. Since our method works purely in Cartesian coordinates, we use a change of basis to obtain the same velocity field in these coordinates.  The initial condition is taken as a compactly supported cosine bell centered at $(1,0,0)$:
\begin{align}
q(\vx,0) = 
\begin{cases}
\frac{1}{2}\lf(1 + \cos\lf(\frac{\pi r}{R_{\rm b}} \rt) \rt) & \text{if}\ r < R_{\rm b}, \\
0 & \text{if}\ r \geq R_{\rm b},
\end{cases}
\end{align}
where $\vx = (x,y,z)$, $r = \arccos(x)$, and the support is set as $R_{\rm b}=1/3$. This initial condition has a jump in the second derivative at $r=R_{\rm b}$, which makes the test susceptible to both dispersive and diffusive errors. The test calls for simulating advection of the initial condition to the final time of $T=2\pi$, which corresponds to one full revolution of the bell over the sphere.  For the convergence tests, we set the time step to $\Delta t= \pi/10$, which is necessary so that spatial errors dominate for all the finest node sets. On the finest node set for the global RBF method, this time-step gives a CFL number of approximately 12, while for the finest node set for the local RBF and RBF-PU methods this gives a CFL number of approximately 28.
\begin{figure}[h!]
\centering
\begin{tabular}{c|c}
Local RBF Method & RBF-PU Method \cr
\hline
{
	\includegraphics[width=0.41\textwidth]{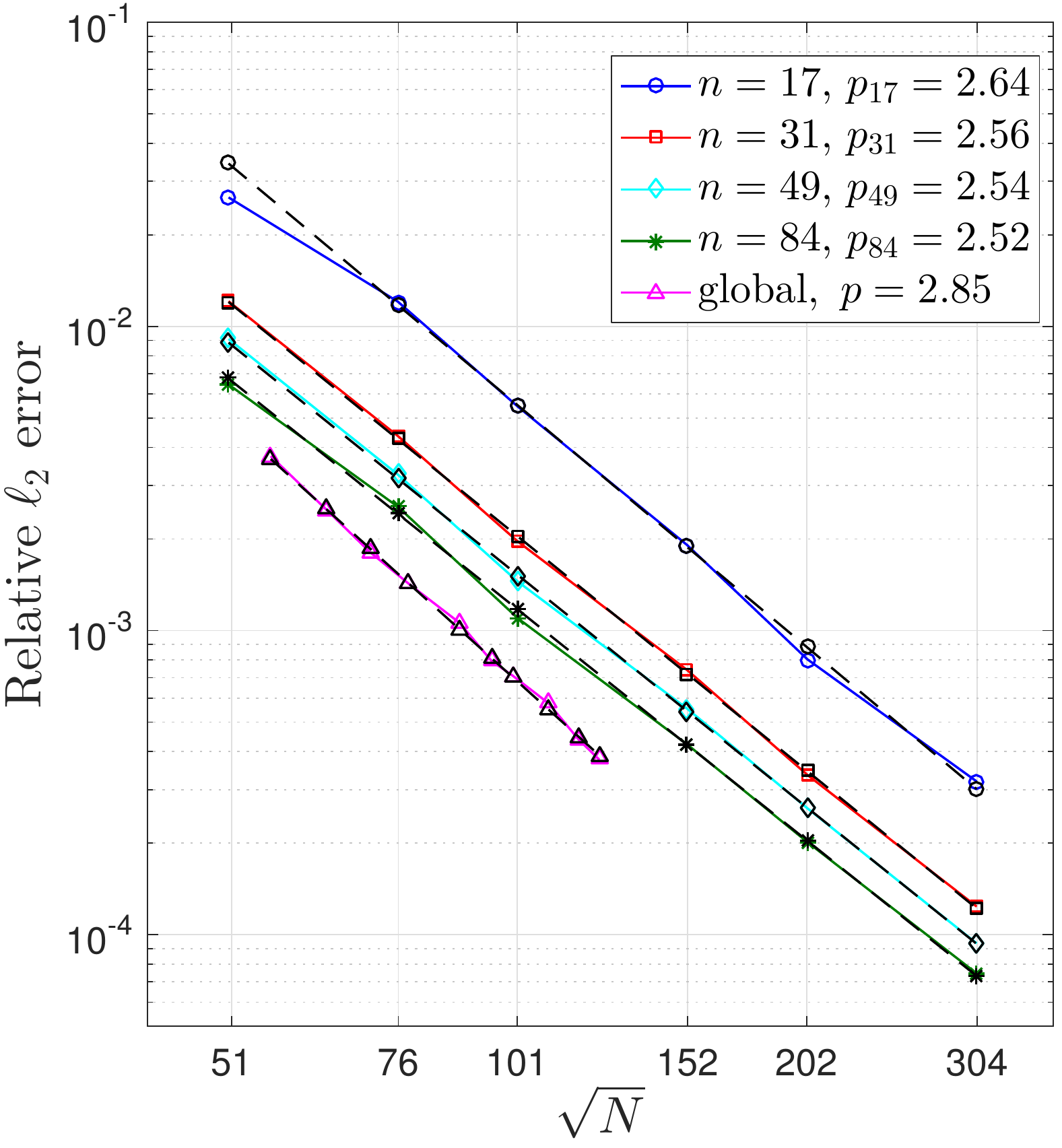} 	
}
&
{
	\includegraphics[width=0.41\textwidth]{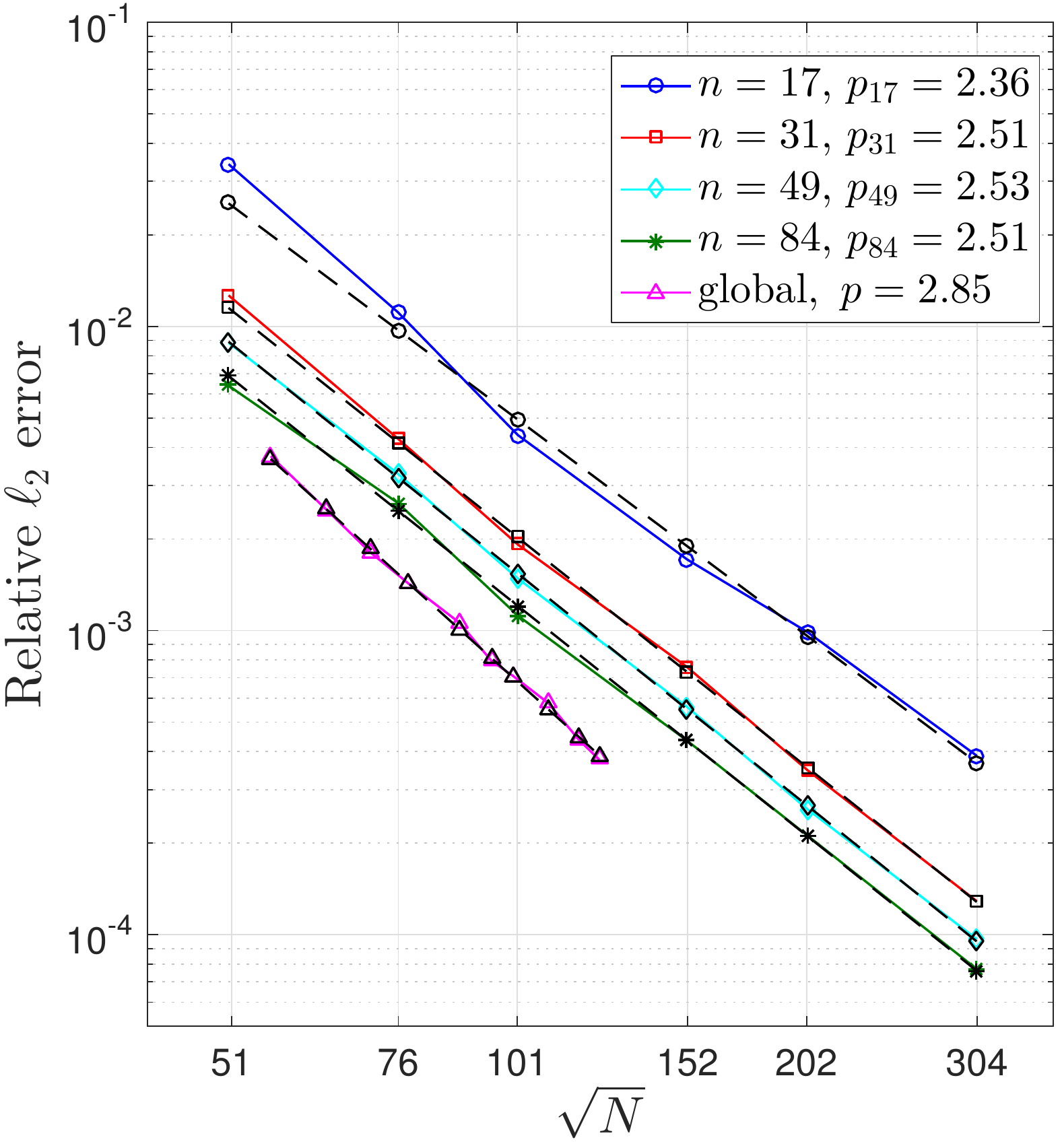} 	
}
\cr
{
	\includegraphics[width=0.41\textwidth]{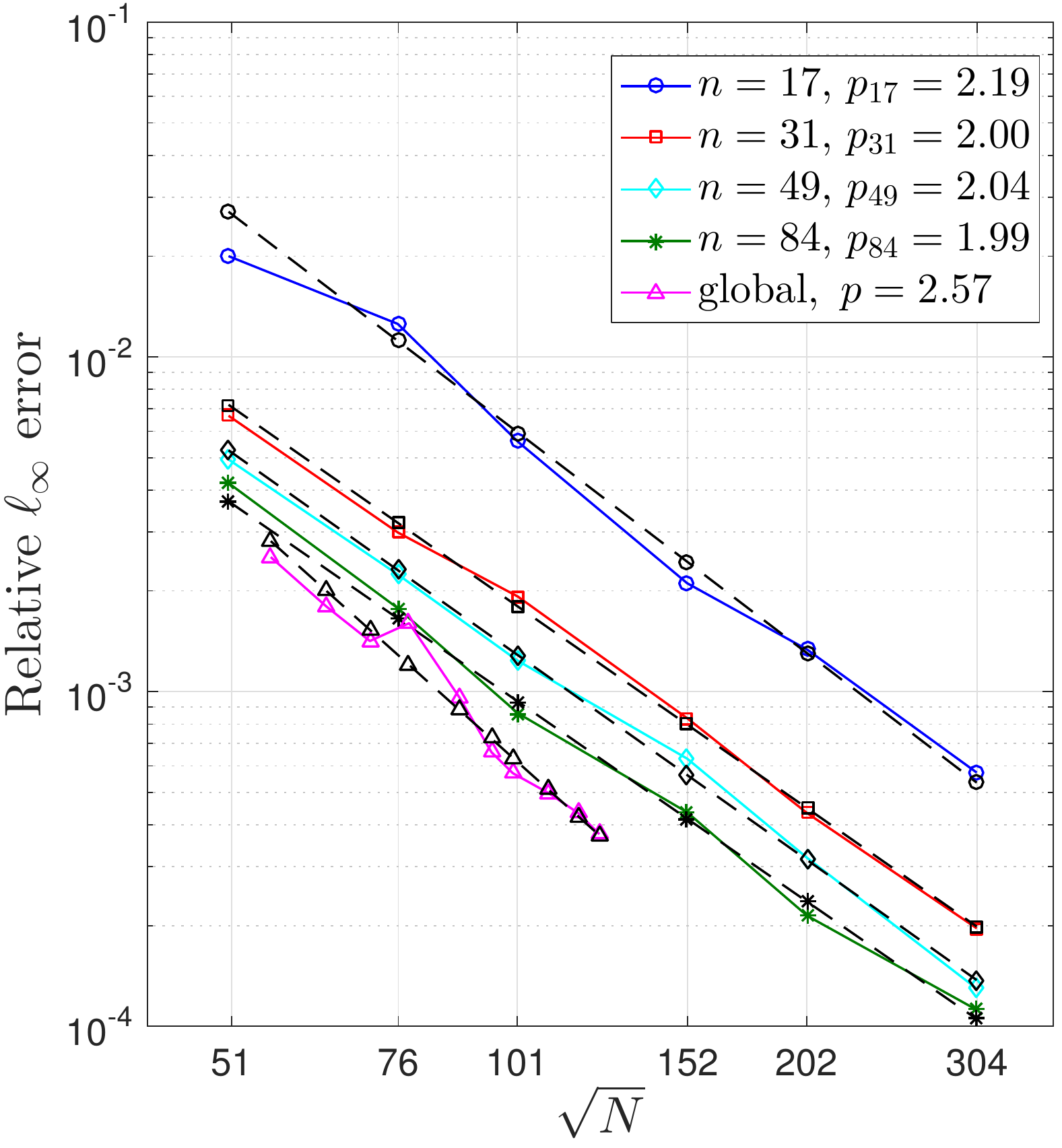} 	
}
&
{
	\includegraphics[width=0.41\textwidth]{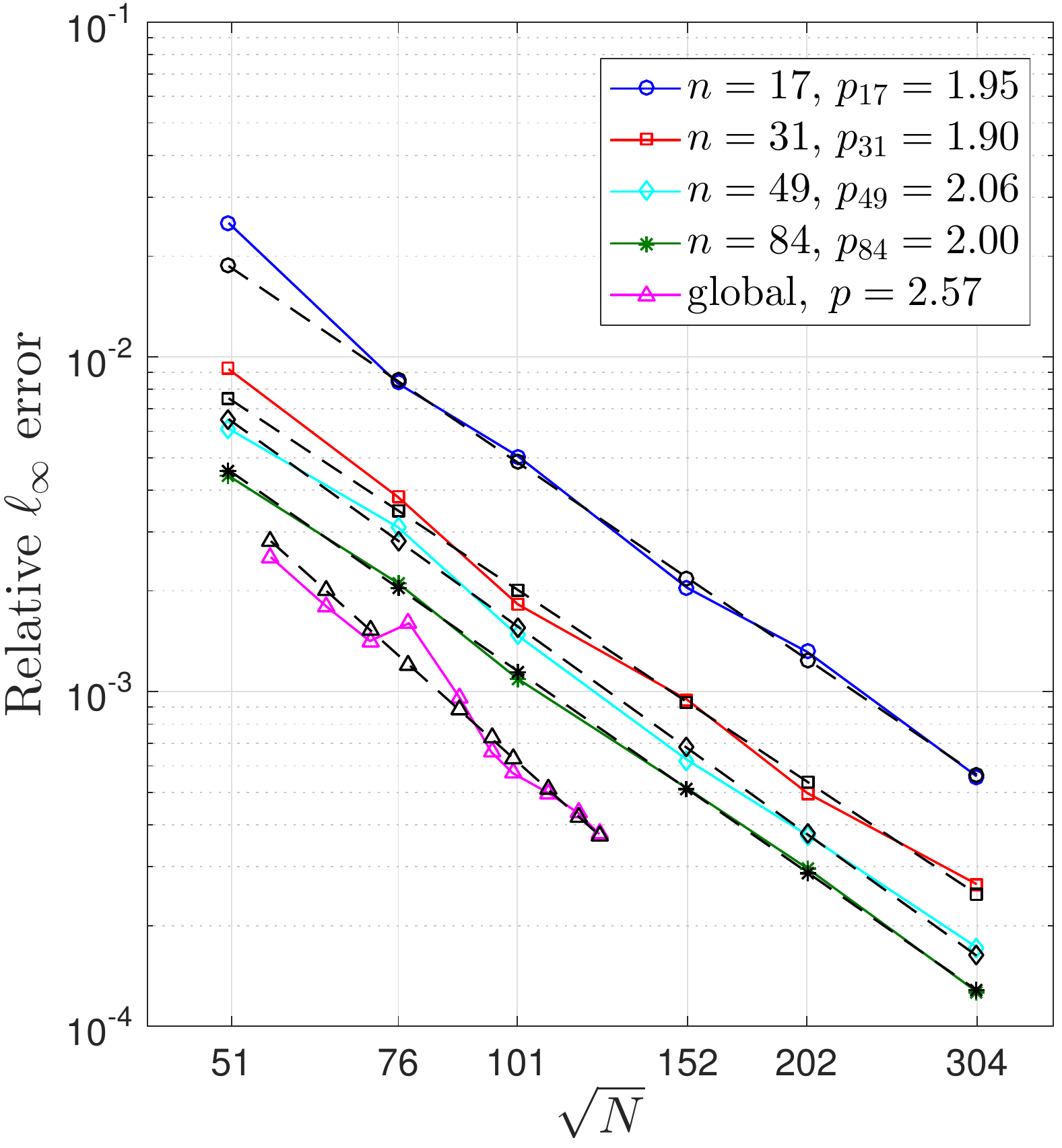} 	
}
\end{tabular}
\caption{Convergence results of the relative $\ell_2$ (top row) and $\ell_{\infty}$ (bottom row) errors for the solid body rotation of a cosine bell test case after one revolution.  The results for the global RBF method are included in all plots for comparison purposes.  The dashed lines are the lines of best fit to the data (without the first point included) of the form $C_n N^{-p_n/2})$.  The values of $p_n$, which estimate the order of accuracy of the methods, are listed in the legend and the subscript on $p$ is dropped for the global RBF method.}
\label{fig:lrbf_sbr2}	
\end{figure}

Convergence results in relative $\ell_2$ and $\ell_{\infty}$ norms for increasing $N$ are shown in Figure \ref{fig:lrbf_sbr2} for the three methods, together with the estimated rates of convergence.  From the top row of this plot, we see that the $\ell_2$ error for the local RBF and RBF-PU methods appears to converges at a rate close to $2.5$ for all $n$, while the convergence for the global RBF method appears closer to $3$.  The $\ell_{\infty}$ error is shown in the bottom row and we see that the convergence rates are now around $2$ for the local RBF and RBF-PU methods and $2.5$ for the global RBF method.  These rates of convergence are dictated by the smoothness of the initial condition, which is only $C^1(\Sph)$, and are consistent with the convergence results for the global RBF collocation and RBF-FD methods for this same problem~\cite{FlyerWright:2007,FoL11}. For the local RBF and RBF-PU methods, we see that increasing $n$ leads to a decrease in the error, but not in the convergence rates (again because of the limited smoothness of the solution).  We also see that the errors for the global method are smaller for similar values of $N$.  

\begin{figure}[h!]
\centering
\begin{tabular}{cc}
{
	\includegraphics[width=0.41\textwidth]{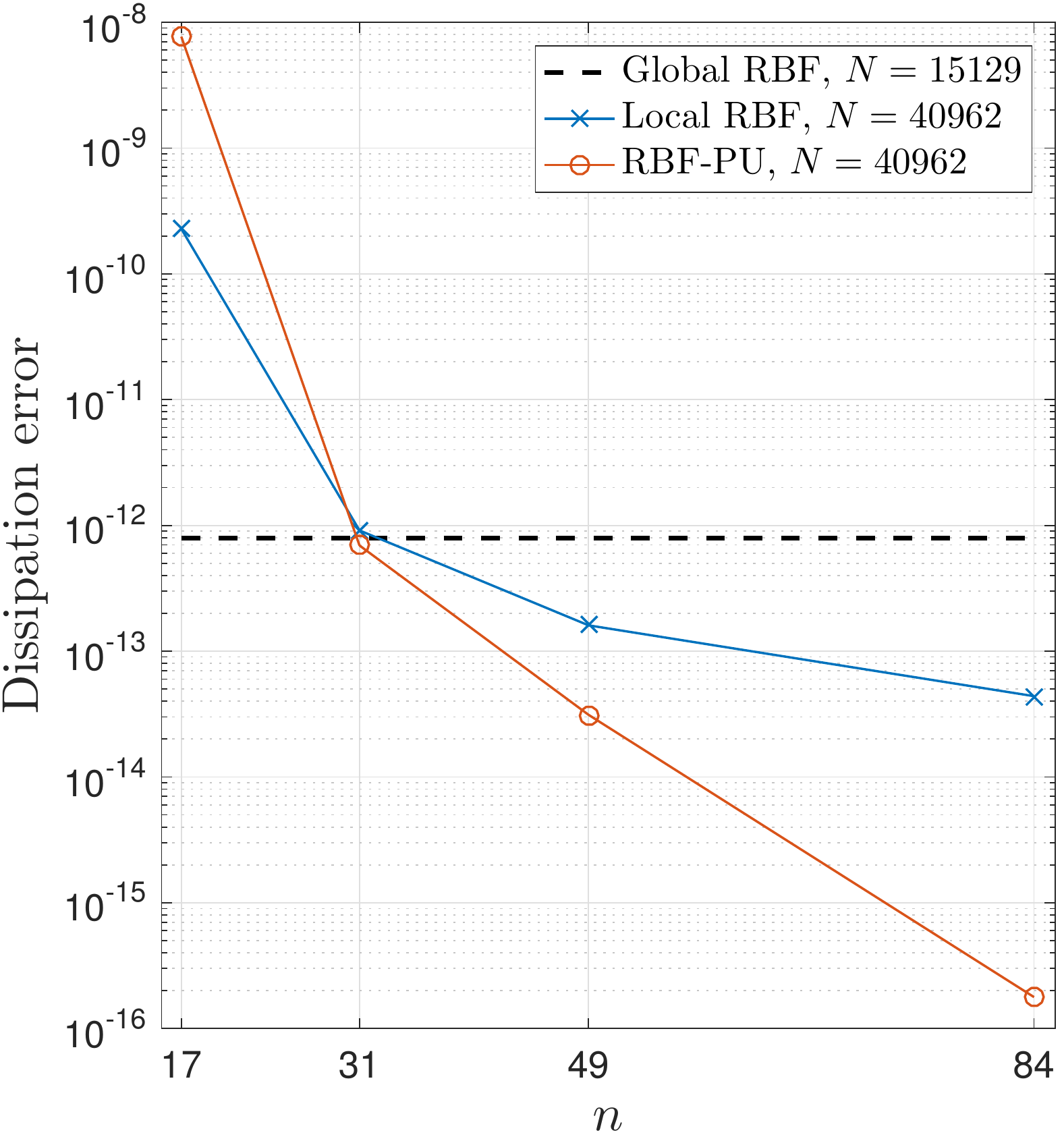} 	
}
&
{
	\includegraphics[width=0.41\textwidth]{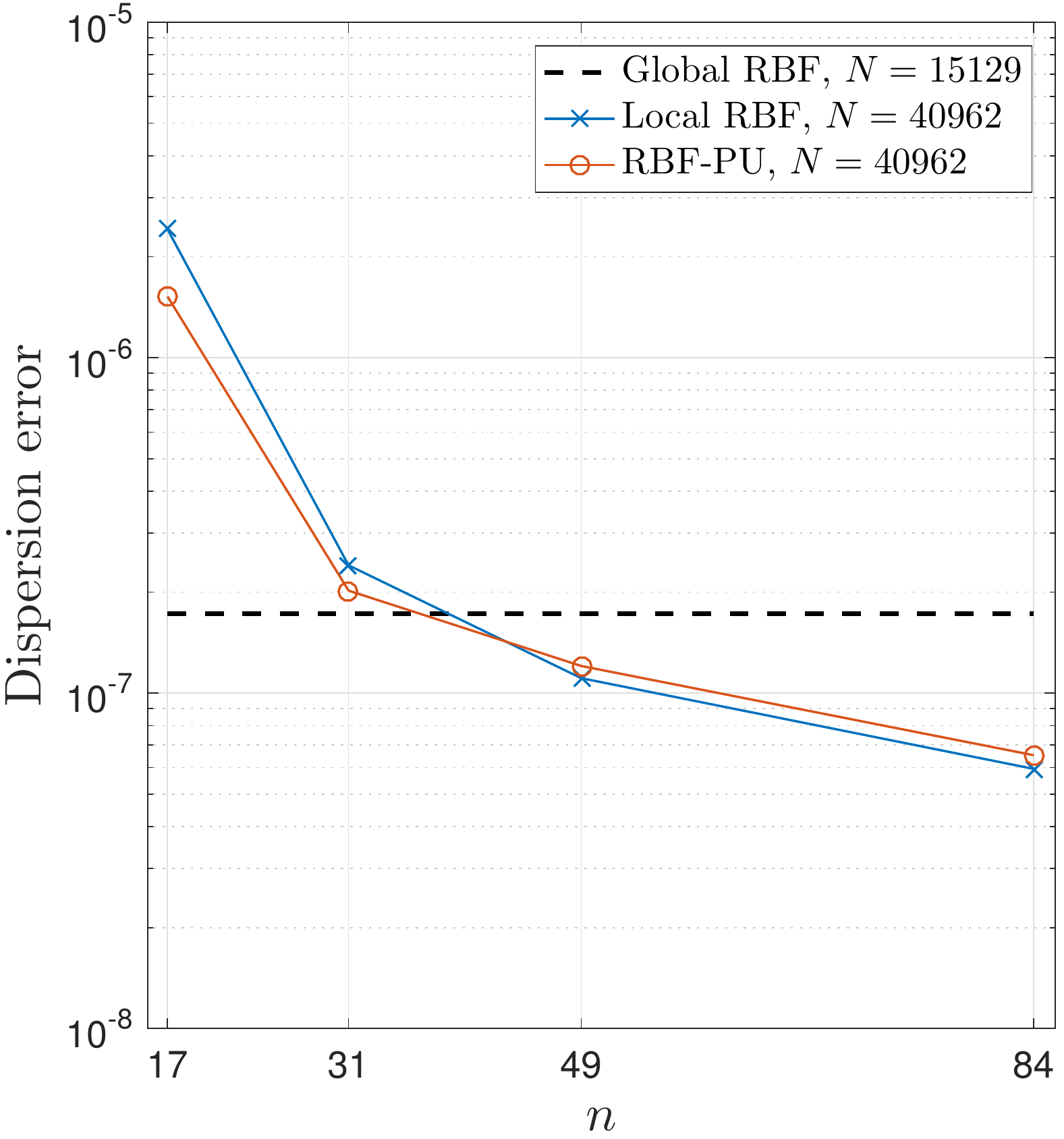} 	
}
\end{tabular}
\caption{Relative dissipation (left) and dispersion (right) errors \eqref{eq:dissipation}--\eqref{eq:dispersion} for the solid body rotation of a cosine bell test case after one revolution as $n$ increases in the local RBF and RBF-PU methods.  The global RBF method does not have a dependence on $n$ and is included as a dashed line for reference.}
\label{fig:dissdisp}	
\end{figure}

Figure \ref{fig:dissdisp} displays the dissipation and dispersion errors \eqref{eq:dissipation}--\eqref{eq:dispersion}.  Here we have fixed $N$ at 40962 for the local RBF and RBF-PU methods and plotted dissipation and dispersion errors against with $n$, the local stencil/patch size.  Also included in the plots are the results for the global RBF method (dashed line) with $N=15129$.  From the plots, we see that dispersion errors dominate the numerical solutions for all the methods, which is expected since the initial condition is only $C^1(\Sph)$.  We also see that for the local RBF and RBF-PU methods increasing $n$ (which increases the order of accuracy of these methods) leads to a decrease in both dissipation and dispersion errors, with the decrease in the former being much more pronounced.  

\begin{figure}[h!]
\centering
\subfloat[Error global RBF: $T = 2\pi$]
{
	\includegraphics[width=0.3\textwidth]{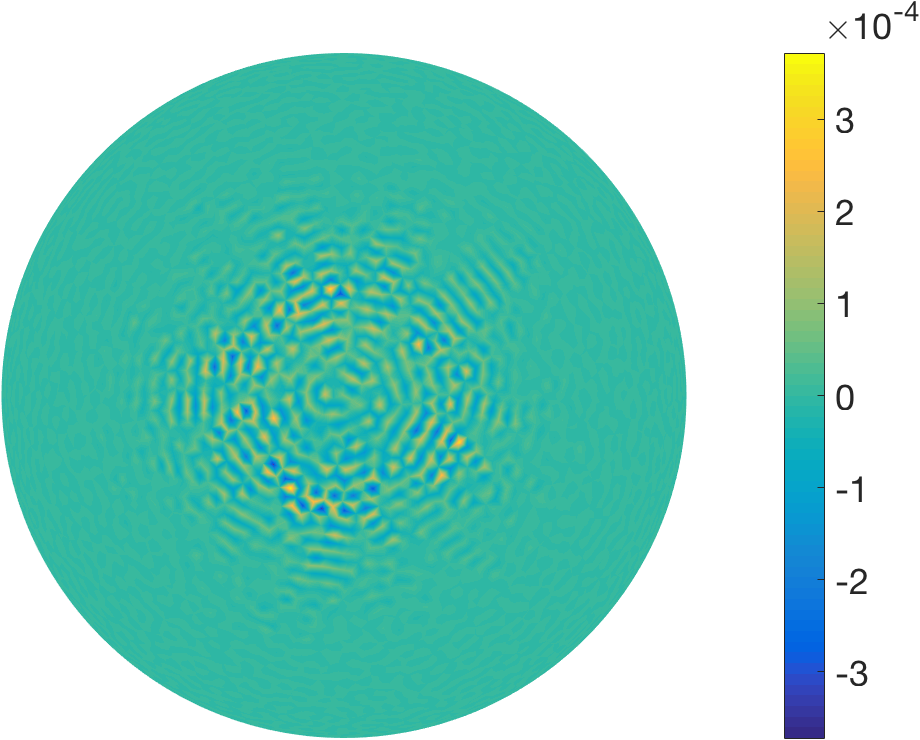} 	
	\label{fig:ldvt1}
}
\subfloat[Error local RBF: $T = 2\pi$]
{
	\includegraphics[width=0.3\textwidth]{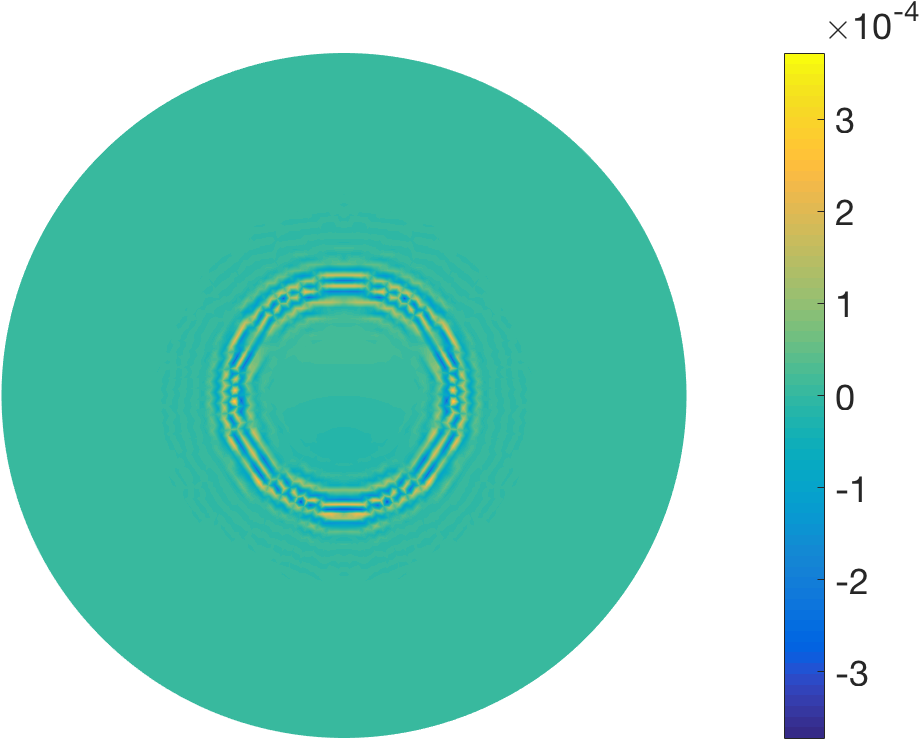} 	
	\label{fig:ldvt2}
}
\subfloat[Error RBF-PU: $T = 2\pi$]
{
	\includegraphics[width=0.3\textwidth]{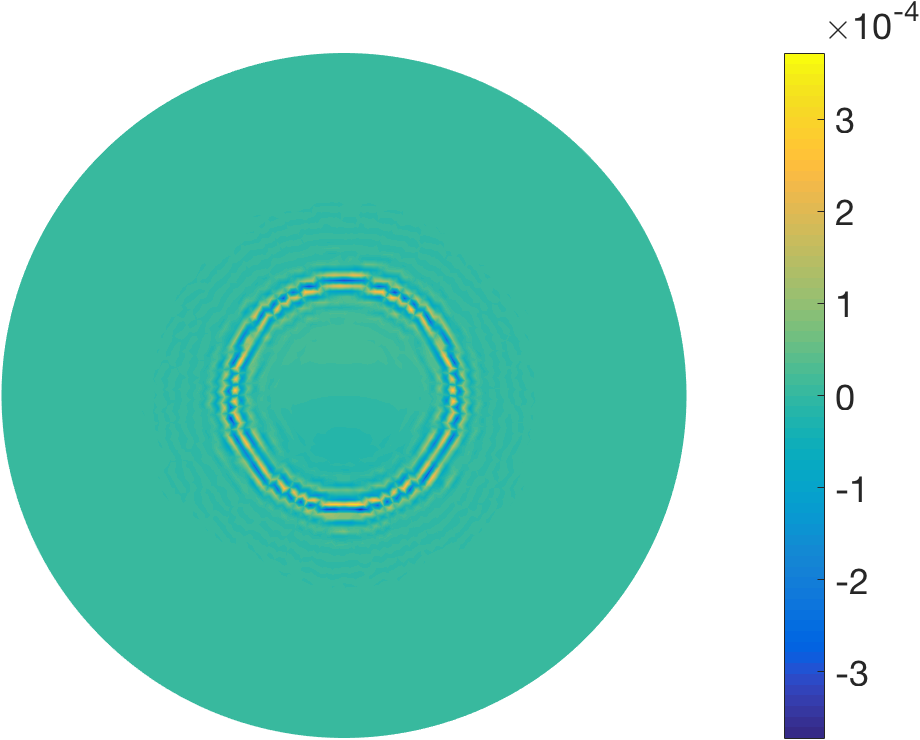}
	\label{fig:pdvt}
}

\subfloat[Error global RBF: $T = 20\pi$]
{
	\includegraphics[width=0.3\textwidth]{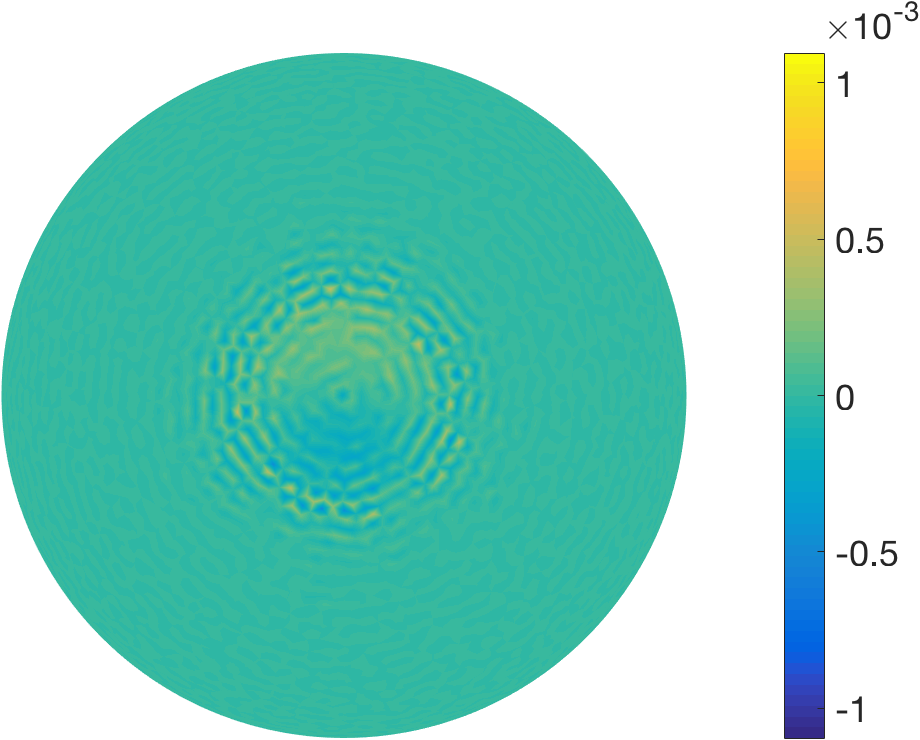} 	
	\label{fig:ldvt3}
}
\subfloat[Error local RBF: $T = 20\pi$]
{
	\includegraphics[width=0.3\textwidth]{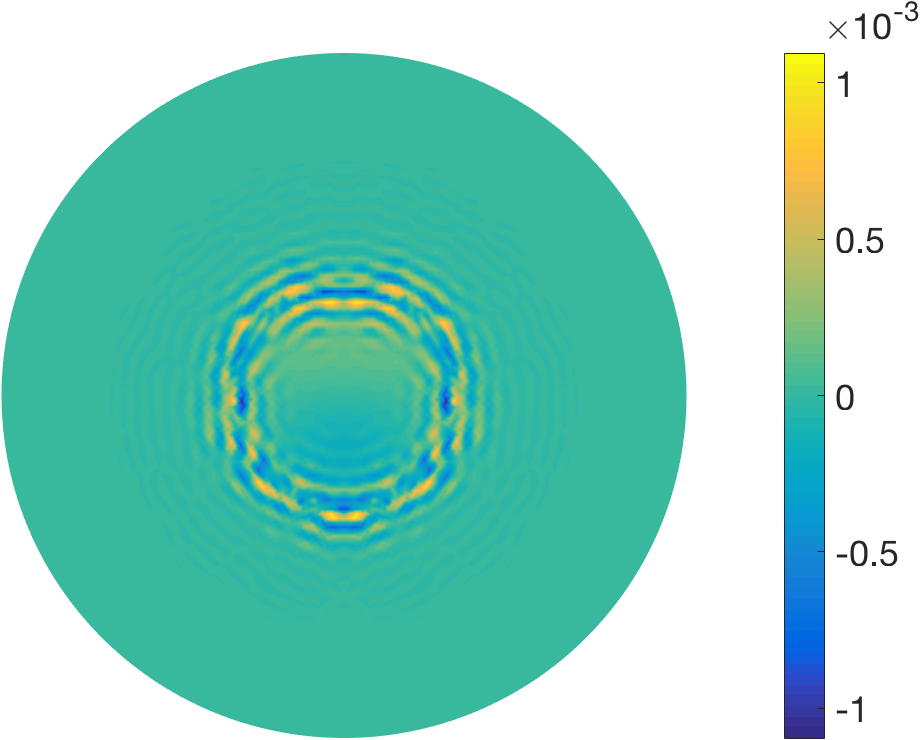} 	
	\label{fig:lcvt}
}
\subfloat[RBF-PU: $T = 20\pi$]
{
	\includegraphics[width=0.3\textwidth]{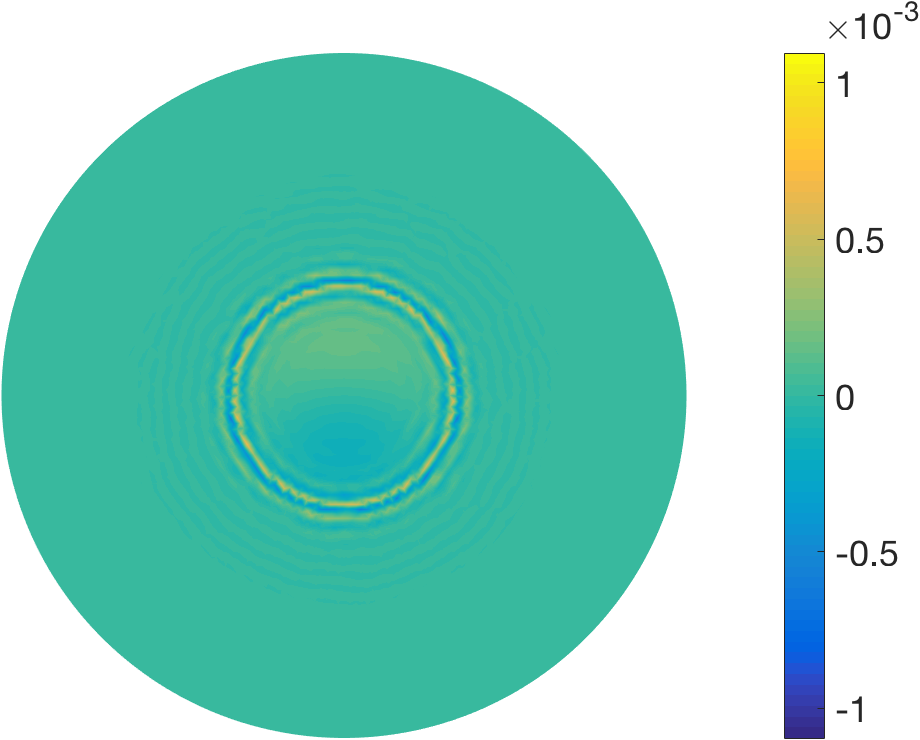}
	\label{fig:pcvt}
}
\caption{Pseudocolor plots of the errors, $q_X-q$, on the sphere for solid body rotation test case of the cosine bell.  The top row shows the errors after one revolution ($T=2\pi$), while the bottom row shows the errors after ten revolutions ($T=20\pi$).  (a) \& (d) show the results for the global RBF method using $N=15129$, (b) \& (d) show the results for the local RBF method using $N=40962$ and $n=49$, and (c) \& (e) show the results for the RBF-PU method using $N=40962$ and $n=49$.   Note the color scales differ for the two values of $T$.  The view for all plots is from the initial center of the bell.}
\label{fig:errloc}	
\end{figure}

We further explore the dissipation/dispersion errors of the methods by plotting the difference in the exact and approximate solutions after one ($T=2\pi$) and ten ($T=20\pi$) revolutions of the cosine bell. This test was suggested in~\cite{FoL11} to examine whether the errors remain well-localized to the support of the bells over time, and to give a visual indication of the dissipation/dispersion errors.  Figure \ref{fig:errloc} displays the results for the global RBF method for the case of N=15129 and the local RBF and RBF-PU methods, both for the case of $N=40962$ and $n=49$.  These values were selected since the errors were of similar magnitude.   We see from the first row of the figure (one revolution) that the errors for the local RBF and RBF-PU methods are localized around the discontinuities in the second derivative of the solution, with the error being lower and more localized for the RBF-PU method.  The errors for the global RBF method are spread more over the support of the entire bell.  None of the methods display a dispersive wave-train over the whole sphere.  After ten revolutions (second row of the figure) the errors for the local RBF and RBF-PU methods increase by about a factor of 4, but still remain localized around the discontinuities in the solution.  The RBF-PU method again displays better localization of the errors, or lower dispersive-type errors.  The dominant errors for the global method also remain restricted to the support of the bell, but only increase by about a factor of $1.5$.   The results for the local RBF and RBF-PU methods are qualitatively similar to those in~\cite{FoL11} for the RBF-FD method with hyperviscosity stabilization.

\begin{figure}[h!]
\centering
	\includegraphics[width=0.5\textwidth]{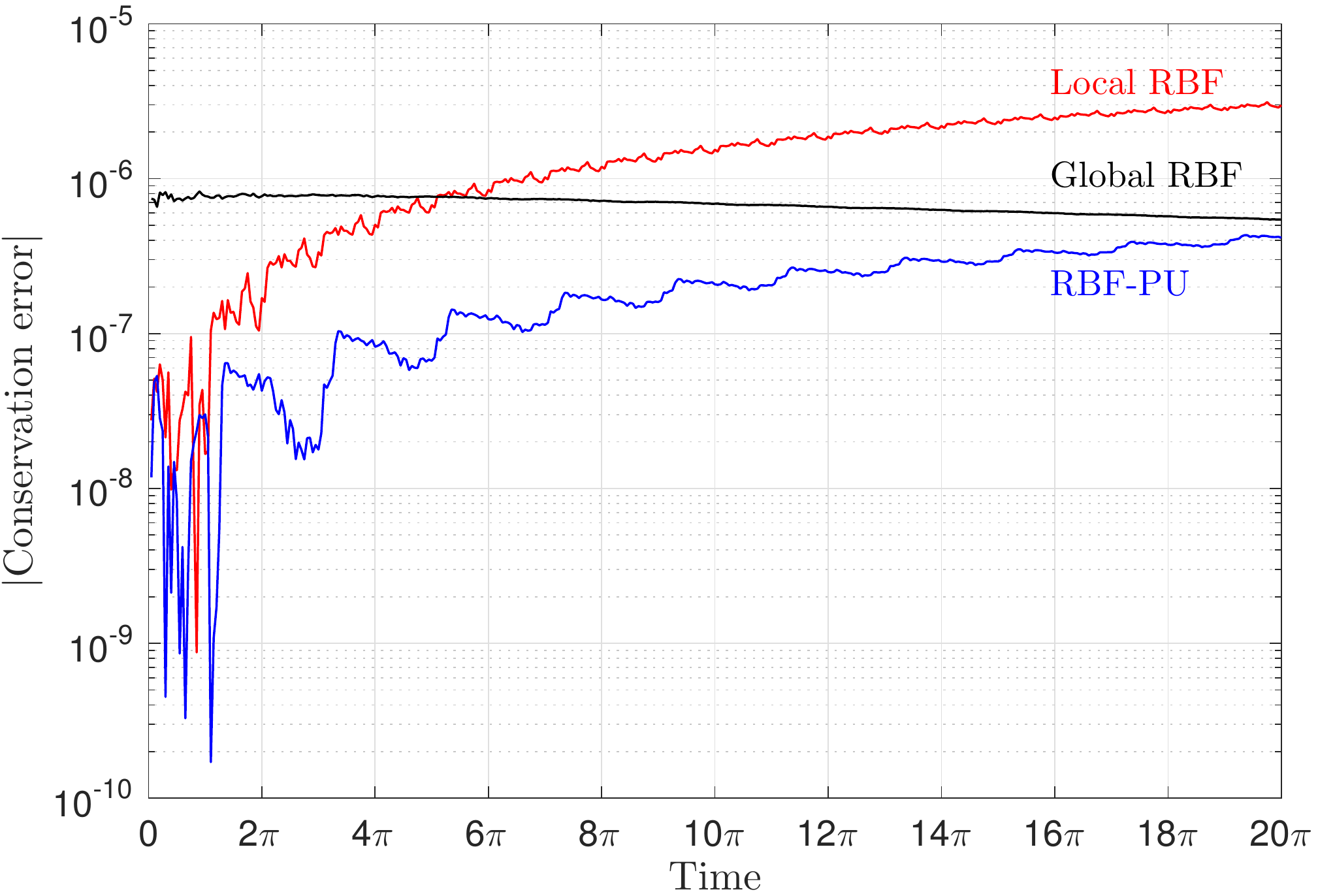}
\caption{Time traces of the mass conservation errors over ten revolutions of the cosine bell in the solid body rotation test case. For the local RBF and RBF-PU methods $N=40962$ and $n=49$, while for the global method $N=15129$.  For this test we used $\Delta t=\pi/20$.}
\label{fig:consSolidBody}	
\end{figure}

Time traces of the absolute value of the mass conservation errors over ten revolutions ($T=20\pi$) of the cosine bell are displayed in Figure \ref{fig:consSolidBody}.  We again use $N=40962$ and $n=49$ for the local RBF and RBF-PU methods and $N=15129$ for the global method. From this figure we see that the conservation errors for the local RBF and RBF-PU methods increase slightly over the integration time, with the growth of the local method being larger, while the errors in the global method remain consistent over the integration interval.  While not presented here, we found that increasing $n$ to 84 for the local RBF and RBF-PU methods had a marginal effect on decreasing the conservation errors.

\subsection{Deformational flow}

For the second test problem, we consider the deformational flow test case from~\cite{RP10}. The components of the velocity field for this test are given in spherical coordinates as
\begin{align}
u(\lambda,\theta,t) &= \frac{10}{T}\cos\lf(\frac{\pi t}{T}\rt)\sin^2\lf(\lambda - \frac{2\pi t}{T}\rt)\sin\lf(2\theta\rt) + \frac{2\pi}{T}\cos\theta, \\
v(\lambda,\theta,t) &= \frac{10}{T}\cos\lf(\frac{\pi t}{T}\rt)\sin\lf(2\lambda - \frac{2\pi t}{T}\rt)\cos(\theta),
\end{align}
where $T=5$. This flow field is designed to deform the initial condition up to time $t=2.5$ and then reverse so that the solution is returned to its initial position and value at time $t=5$. This value serves as the final time for the simulation. As before, we use a simple change of basis to convert this velocity field into Cartesian coordinates. 

The following two initial conditions are considered:
\begin{enumerate}
\item Two cosine bells: $q(\vx,0) = 0.1 + 0.9(q_1(\vx,0) + q_2(\vx,0))$,
where for $j=1,2$
\begin{align*}
q_j(\vx,0) &=
\begin{cases}
\frac{1}{2}\lf(1 + \cos\lf(2\pi \cos^{-1}\lf(\vx \cdot \vp_j\rt) \rt) \rt) & \text{if}\ \cos^{-1}\lf(\vx \cdot \vp_j\rt) < \frac{1}{2}, \\
0 & \text{otherwise}.
\end{cases}
\end{align*}
\item Two Gaussian bells:  $
q(\vx,0) = 0.95\lf( e^{-5\|\vx - \vp_1\|_2^2} + e^{-5\|\vx - \vp_2\|^2_2} \rt)$.
\end{enumerate}
In both cases, $\vp_1 = \lf(\frac{\sqrt{3}}{2},\frac{1}{2},0 \rt)$ and $\vp_2 = \lf(\frac{\sqrt{3}}{2},-\frac{1}{2},0 \rt)$.  Like the previous example, the first initial condition is a good test of the sensitivity of the three methods to dispersion errors as it is only $C^1(\Sph)$.  The second initial condition is $C^{\infty}(\Sph)$ and hence is a good test of the maximum convergence rates of the three methods.  For all experiments involving the local RBF and RBF-PU methods, we set $\Delta t = 1/10$.  This gives a CFL number of approximately 27 on the finest node set.  For the global RBF method we use $\Delta t = 1/10$ for the cosine bells and $\Delta t = 1/40$ for the Gaussian bells, which results in CFL numbers of $12$ and $3$, respectively, on the finest node sets.  All time-steps were chosen so that spatial errors dominate for all values of $N$.

\subsubsection{Results for the cosine bells}

\begin{figure}[h!]
\centering
\begin{tabular}{c|c}
Local RBF Method & RBF-PU Method \cr
\hline
{
	\includegraphics[width=0.41\textwidth]{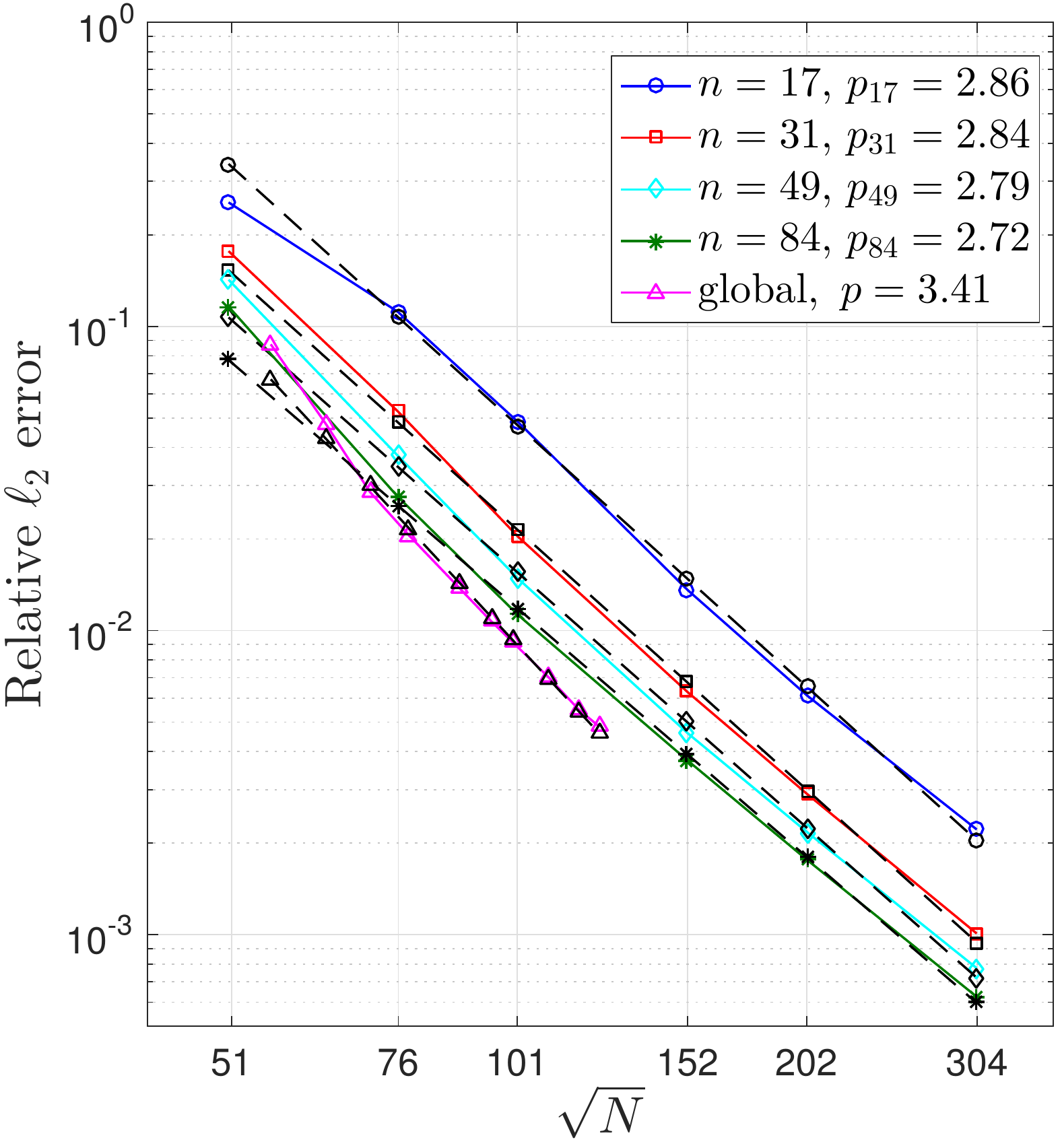} 	
}
&
{
	\includegraphics[width=0.41\textwidth]{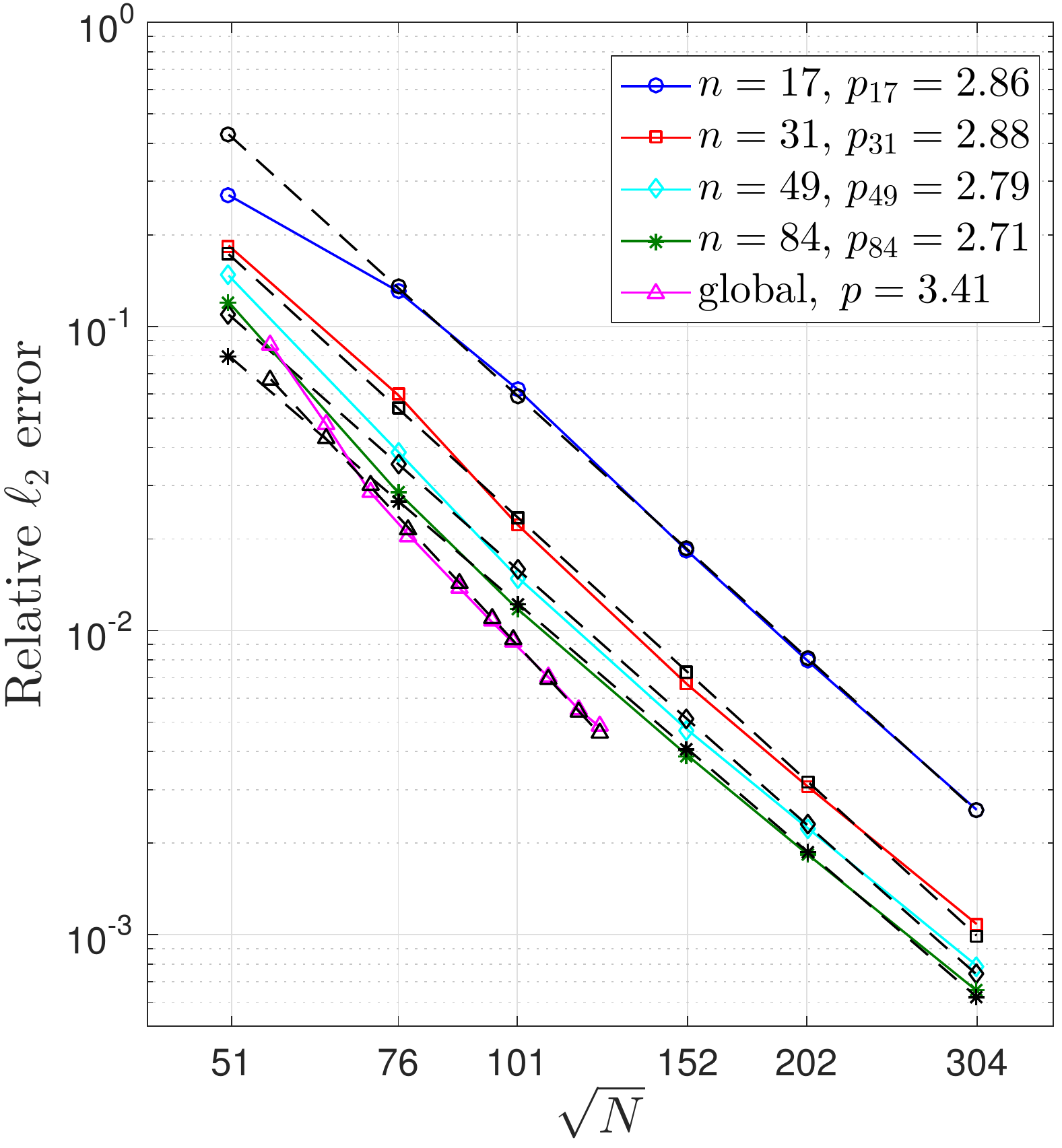} 	
}
\cr
{
	\includegraphics[width=0.41\textwidth]{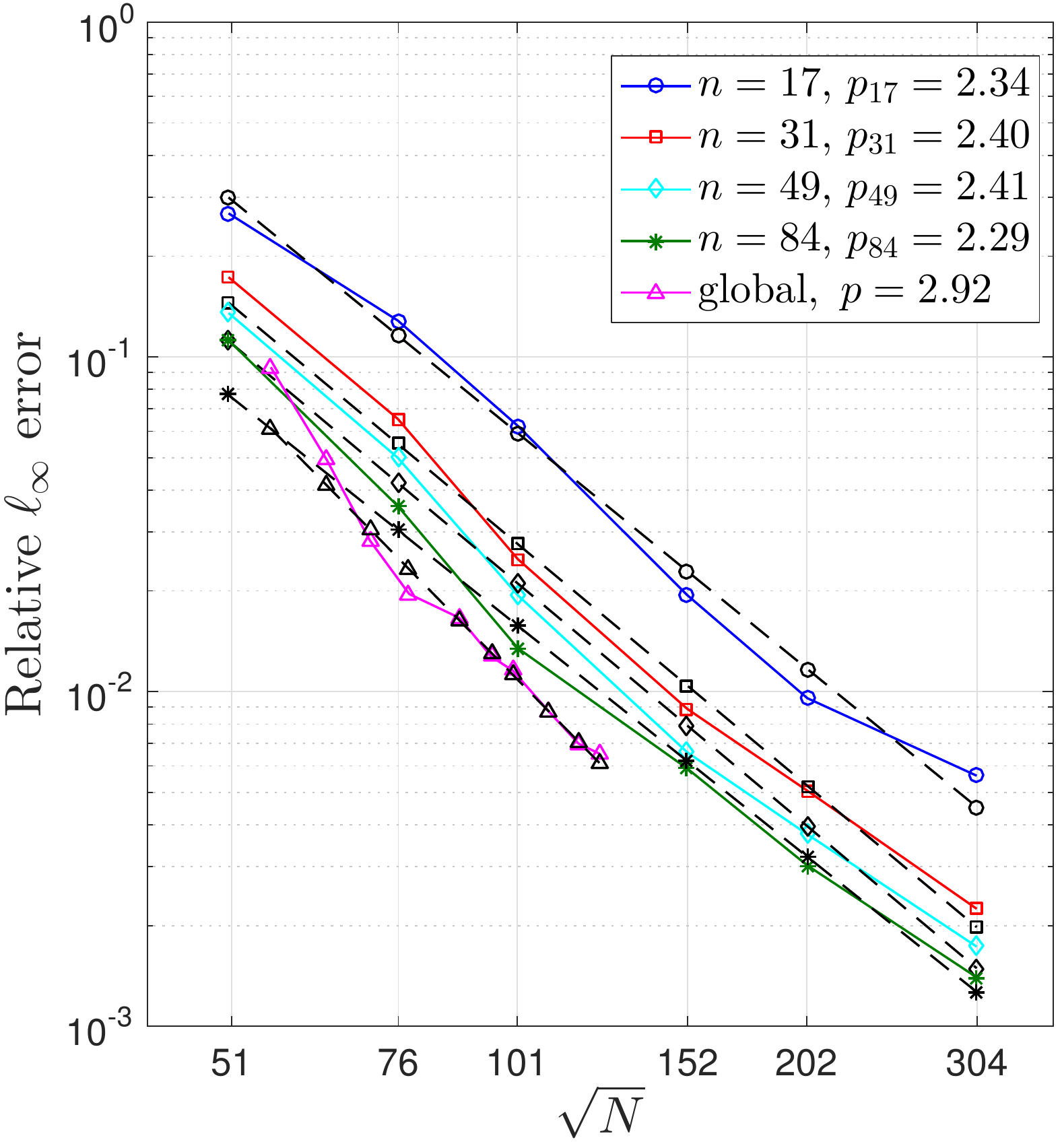} 	
}
&
{
	\includegraphics[width=0.41\textwidth]{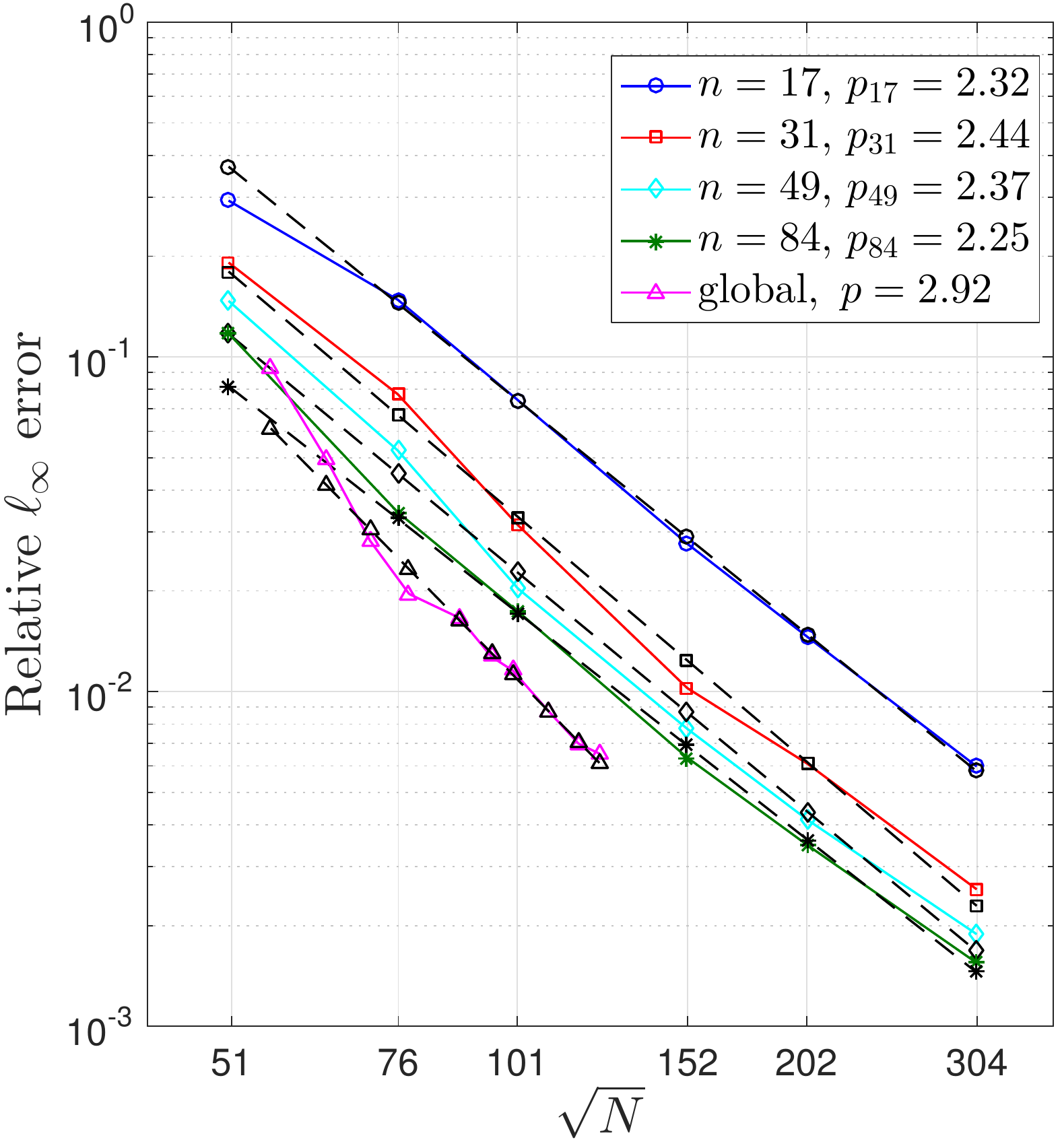} 	
}
\end{tabular}
\caption{Convergence results of the relative $\ell_2$ (top row) and $\ell_{\infty}$ (bottom row) errors the deformational flow of two cosine bells test case.  The results for the global RBF method are included in all plots for comparison purposes.  The dashed lines are the lines of best fit to the data (without the first point included) of the form $C_n N^{-p_n/2})$.  The values of $p_n$, which estimate the order of accuracy of the methods, are listed in the legend and the subscript on $p$ is dropped for the global RBF method.}
\label{fig:lrbf_def2}	
\end{figure}

Figure \ref{fig:lrbf_def2} displays the convergence results in relative $\ell_2$ and $\ell_{\infty}$ norms for increasing $N$ for the three methods, together with the estimated rates of convergence.  We see that the convergence rates for this test case are a bit higher than the solid body rotation test, but that the errors are larger for a given $N$.  The rates of convergence are again limited by the smoothness of the solution.  The larger errors for this test are expected as the solution undergoes much more dramatic changes over the simulation period than the solid body rotation test.  As in the solid body rotation test, the errors for the local RBF and RBF-PU methods are larger than the global method for a given $N$, but because the former methods are more computationally efficient, we can push them to larger $N$ and reach smaller overall errors.  At the end of this section, we compare the accuracy of all three methods to their runtime costs to get a better picture of their overall efficiency.

The relative dissipation and dispersion errors \eqref{eq:dissipation}--\eqref{eq:dispersion} are displayed in Figure \ref{fig:dissdisp}.  As in the previous test, $N$ is fixed for the local RBF and RBF-PU methods and the errors are plotted against $n$.  The global RBF method with $N=15129$ is displayed as a dashed line. The figure shows that dispersion errors again dominate the numerical solutions for all the methods and that increasing $n$ for the local RBF and RBF-PU methods leads to a decrease in these errors.  Visual depictions of these dispersive errors are given in Figure \ref{fig:errloc_deform_cosine}.  Part (a) shows the initial condition and final solution for this test case and parts (b)--(d) show the difference between this solution and the numerical solutions at the final time for the three methods.  The error for the global method (part (b)) are about 1.5 times higher than the local RBF and RBF-PU methods and appears to be much more dispersive in nature.  Also, the error for the RBF-PU method is more localized to the discontinuities of the initial condition than both the global and local RBF methods. 

\begin{figure}[h!]
\centering
\begin{tabular}{cc}
{
	\includegraphics[width=0.41\textwidth]{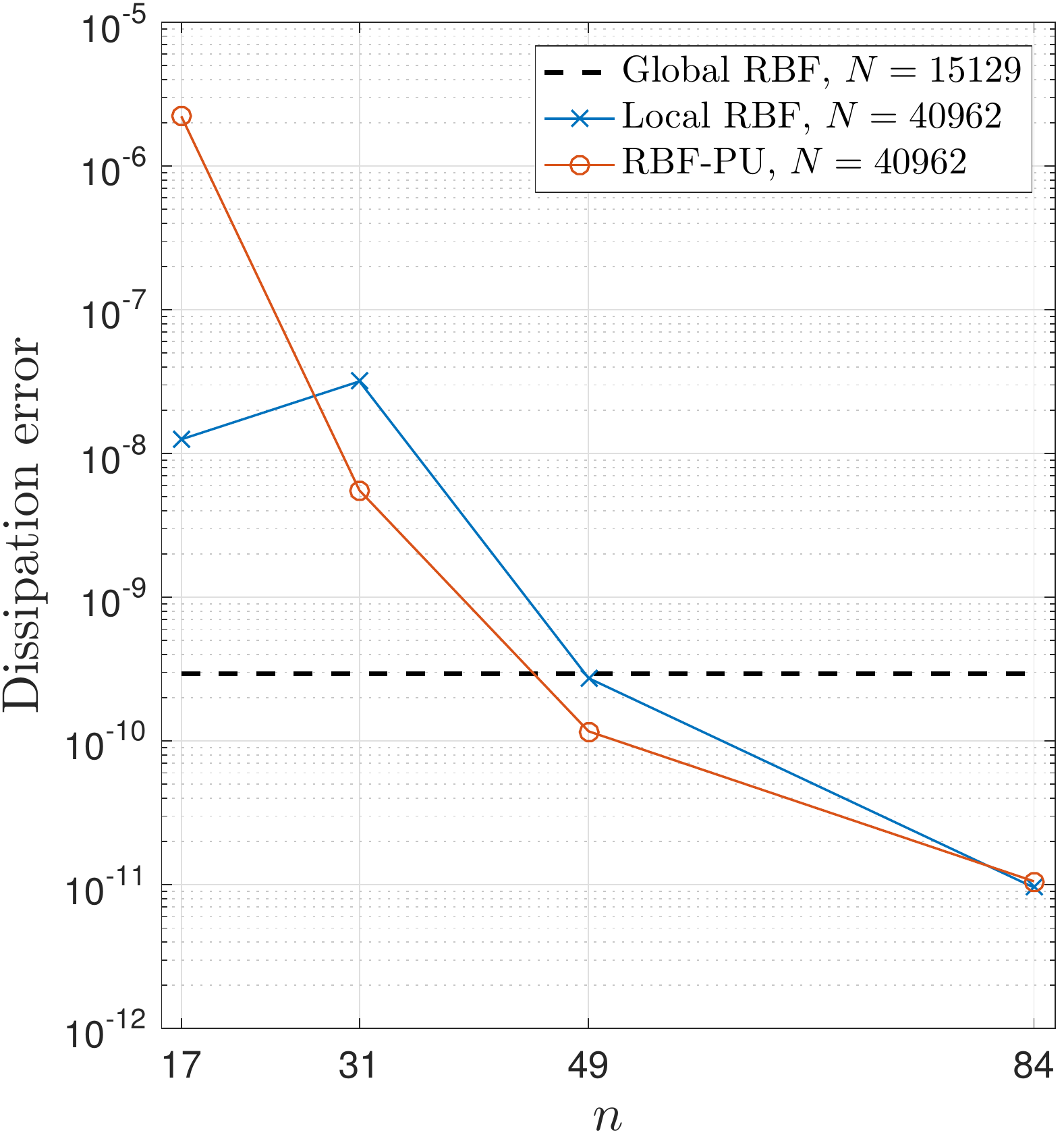} 	
}
&
{
	\includegraphics[width=0.41\textwidth]{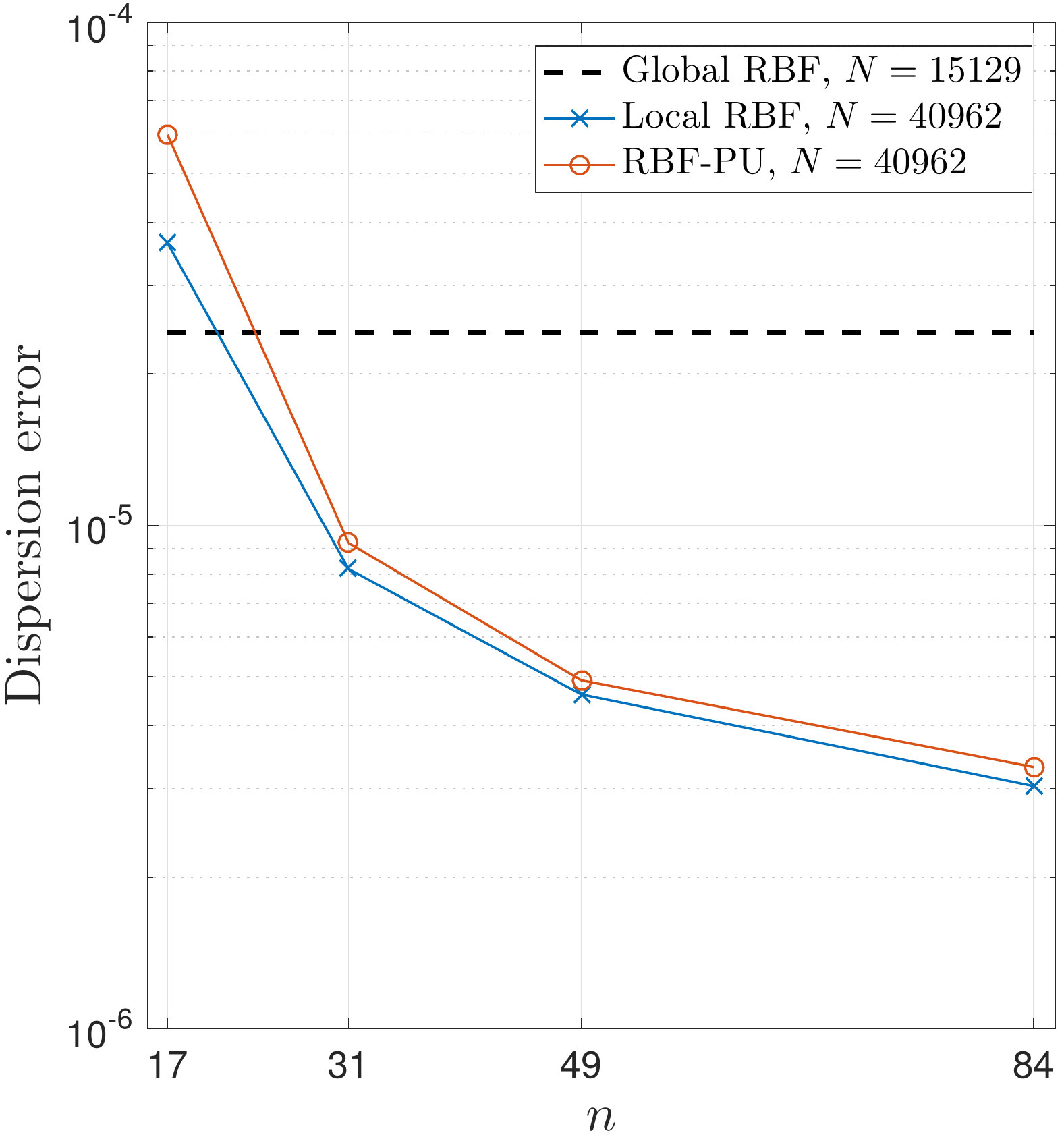} 	
}
\end{tabular}
\caption{Relative dissipation (left) and dispersion (right) errors \eqref{eq:dissipation}--\eqref{eq:dispersion} for the deformational flow test case of two cosine bells after one revolution as $n$ increases in the local RBF and RBF-PU methods.  The global RBF method does not have a dependence on $n$ and is included as a dashed line for reference.}
\label{fig:diss_disp_deform_cosine}	
\end{figure}

Time traces of the mass conservation errors for the three methods over the simulation time are displayed in Figure \ref{fig:consDeformFlow} (a).   For the local RBF and RBF-PU methods, we see that, after a relatively large initial growth, the errors level off around $t=1$ and then start to grow very slowly towards the end of the simulation. The global RBF error is overall larger than the other two methods, but does not exhibit a discernible growth rate.

\begin{figure}[h!]
\centering
\begin{tabular}{cc}
\includegraphics[width=0.35\textwidth]{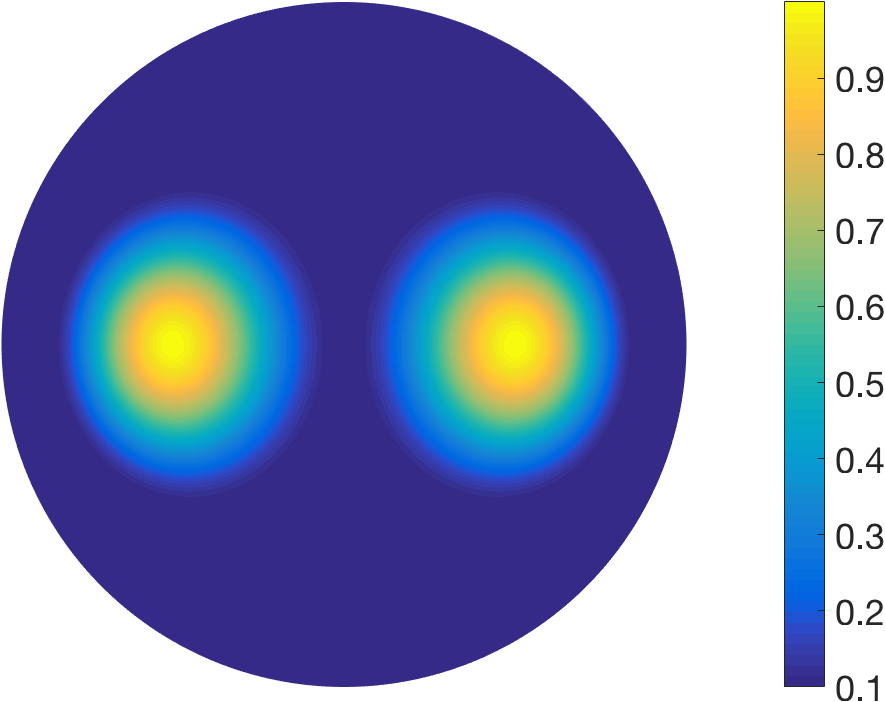}  &
\includegraphics[width=0.35\textwidth]{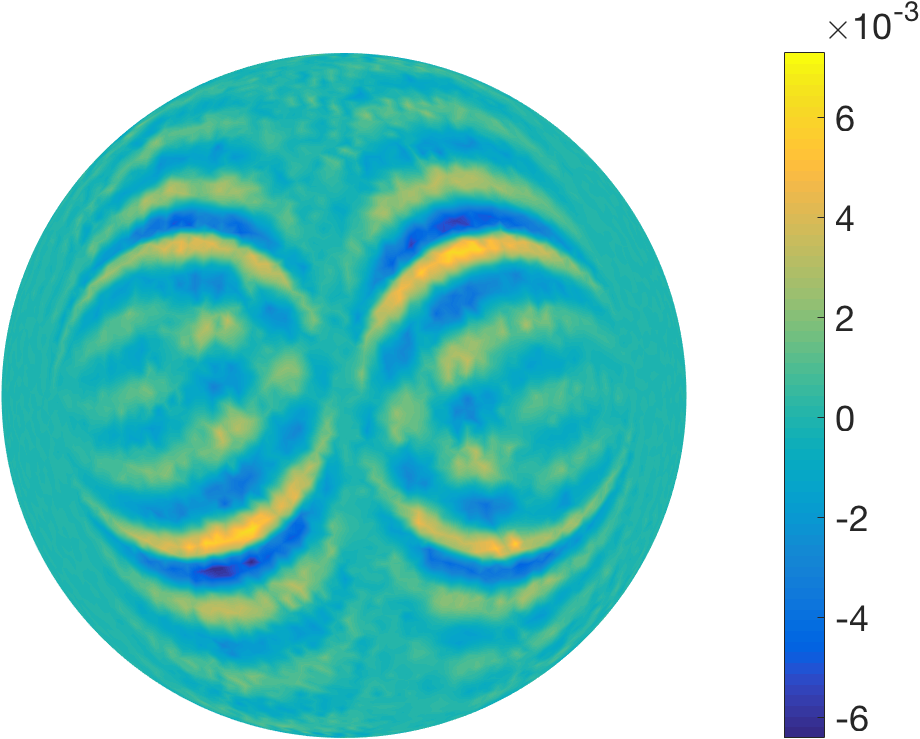} \\
{\footnotesize (a) Exact solution: T=0 \& 5} & (b) {\footnotesize Error global RBF, $T=5$} \\
\includegraphics[width=0.35\textwidth]{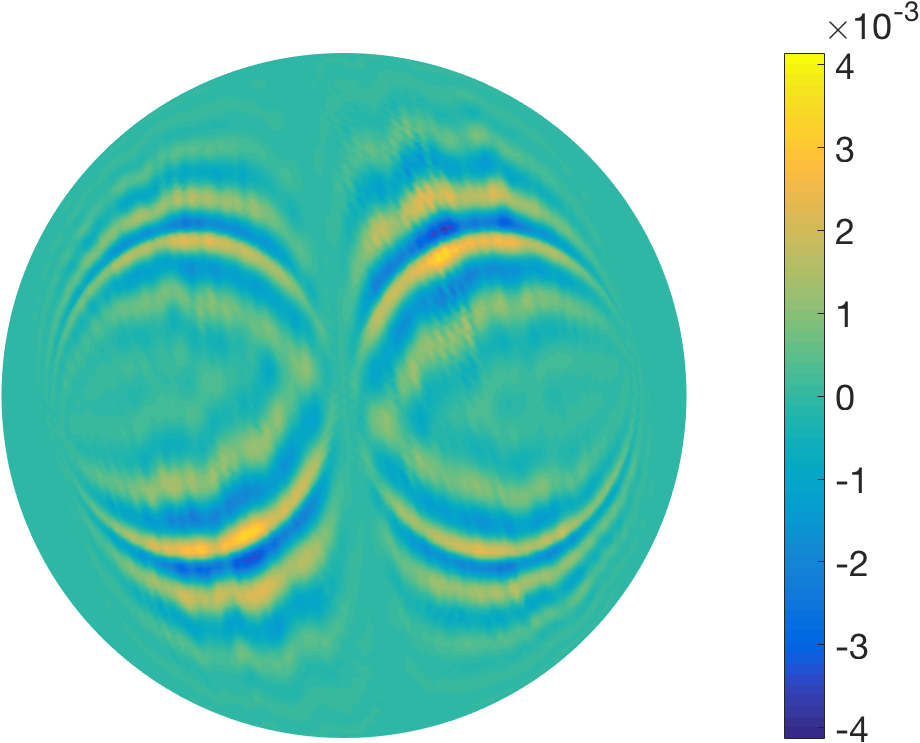} & 
\includegraphics[width=0.35\textwidth]{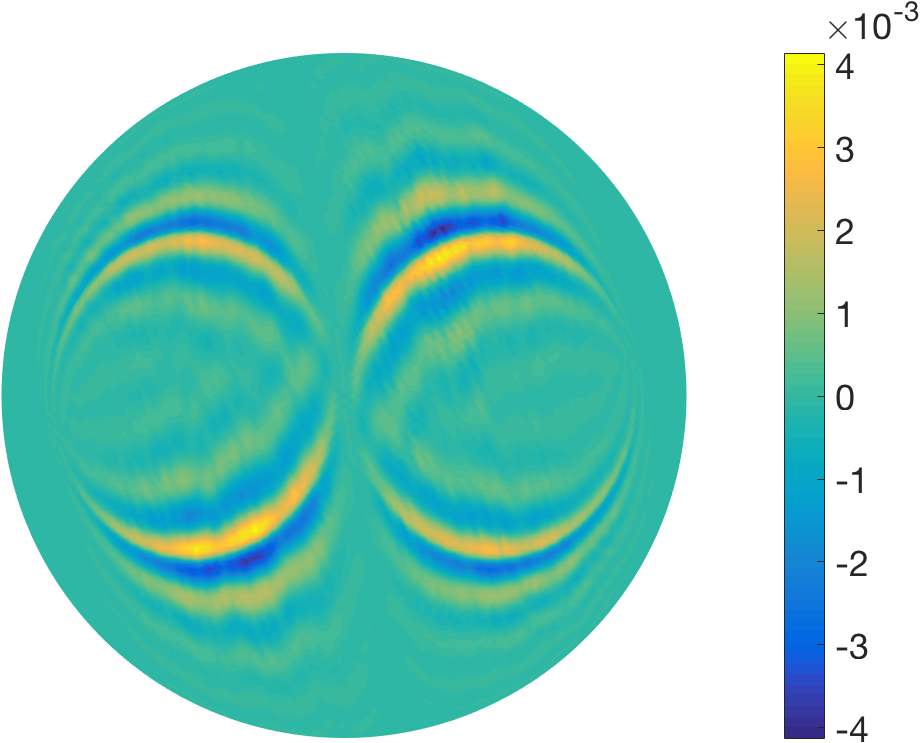} \\
{\footnotesize (c) Error local RBF: $T=5$} & {\footnotesize (d) Error RBF-PU: $T=5$}
\end{tabular}
\caption{(a) Pseudocolor plots on the sphere of the initial and final solutions for the two cosine bells test.  (b)--(c) Show pseudocolor plots of the error, $q_X-q$, at $T=5$ for the global RBF, local RBF, and RBF-PU methods, respectively.  For the global RBF method, $N=15129$, while for the local and RBF-PU methods, $N=40962$ and $n=49$.  Note that the color scale is about 1.5 times larger for the global RBF method.  The view for all plots is from the midpoint between the centers of the two bells.}
\label{fig:errloc_deform_cosine}	
\end{figure}

\begin{figure}[h!]
\centering
\begin{tabular}{cc}
	\includegraphics[width=0.48\textwidth]{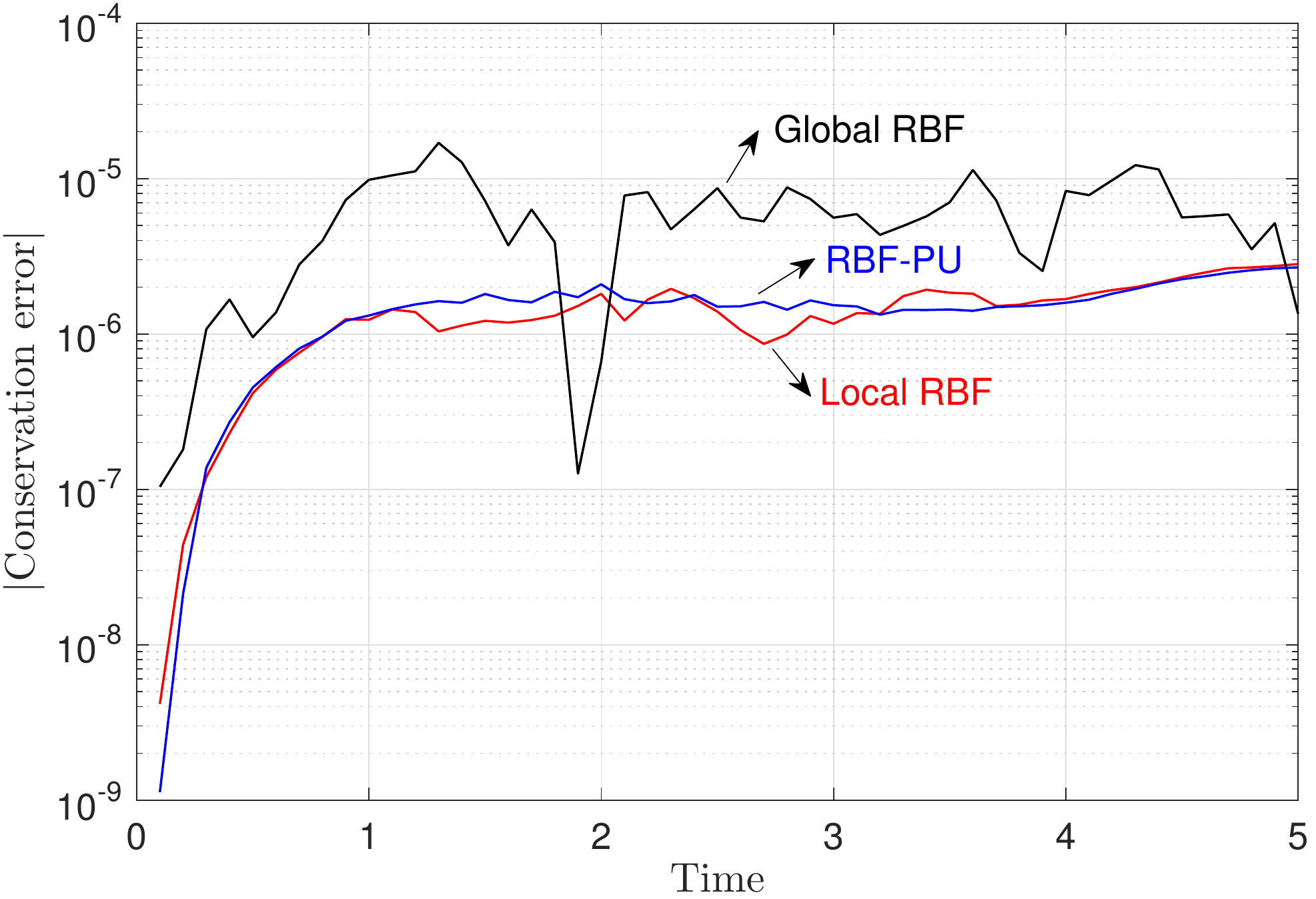} &
	\includegraphics[width=0.48\textwidth]{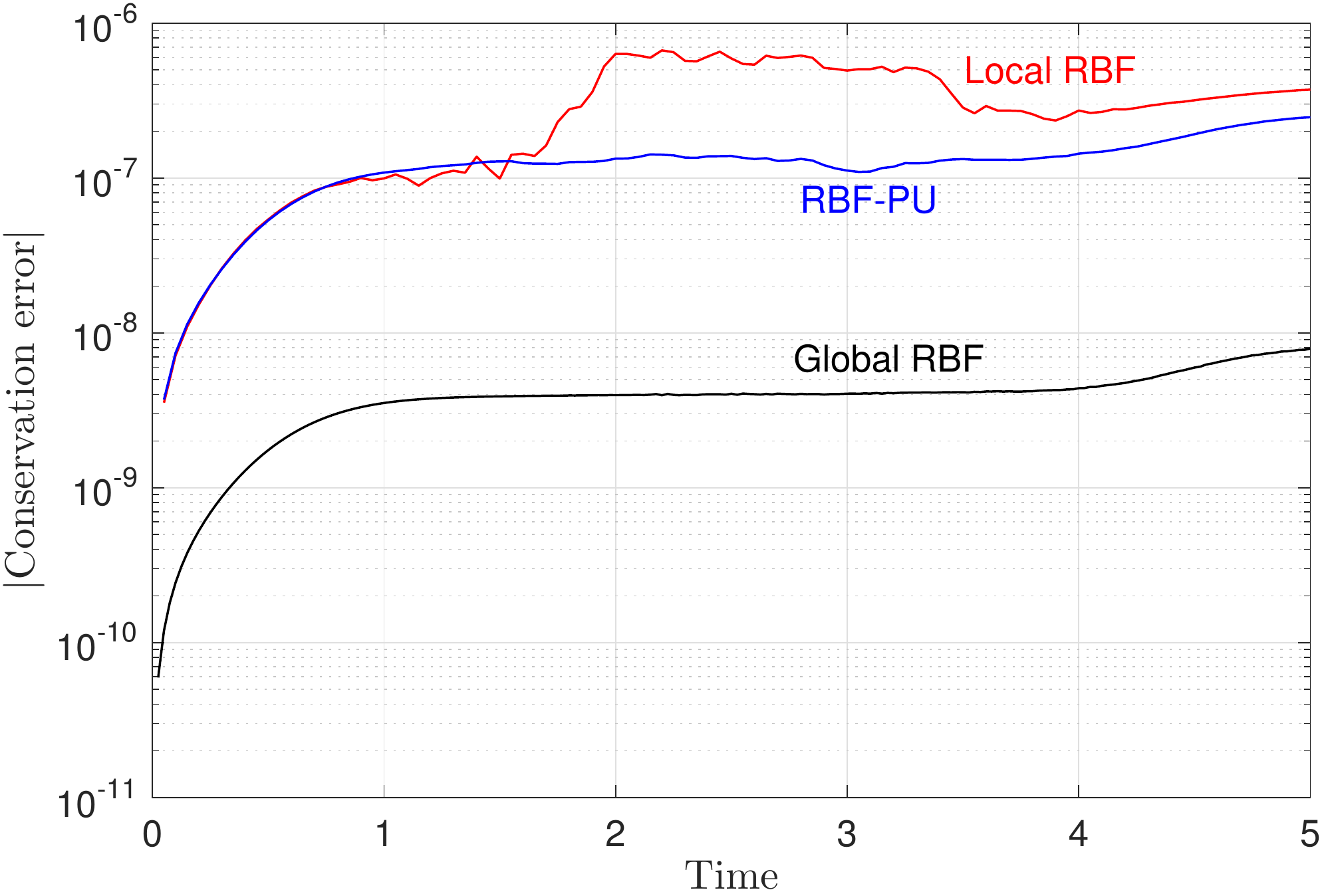} \\
	(a) Cosine bells & (b) Gaussian bells
\end{tabular}
\caption{Time traces of the mass conservation errors over the simulation time for the deformational flow tests.  For the local RBF and RBF-PU methods $N=92162$ and $n=49$ and $\Delta t = 1/10$ for both part (a) and (b), while for the global method $N=15129$ and $\Delta t = 1/10$ for part (a) and $N=15129$ and $\Delta t = 1/40$ for part (b).}
\label{fig:consDeformFlow}	
\end{figure}

\subsubsection{Results for the Gaussian bells}

\begin{figure}[h!]
\centering
\begin{tabular}{c|c}
Local RBF Method & RBF-PU Method \cr
\hline
{
	\includegraphics[width=0.46\textwidth]{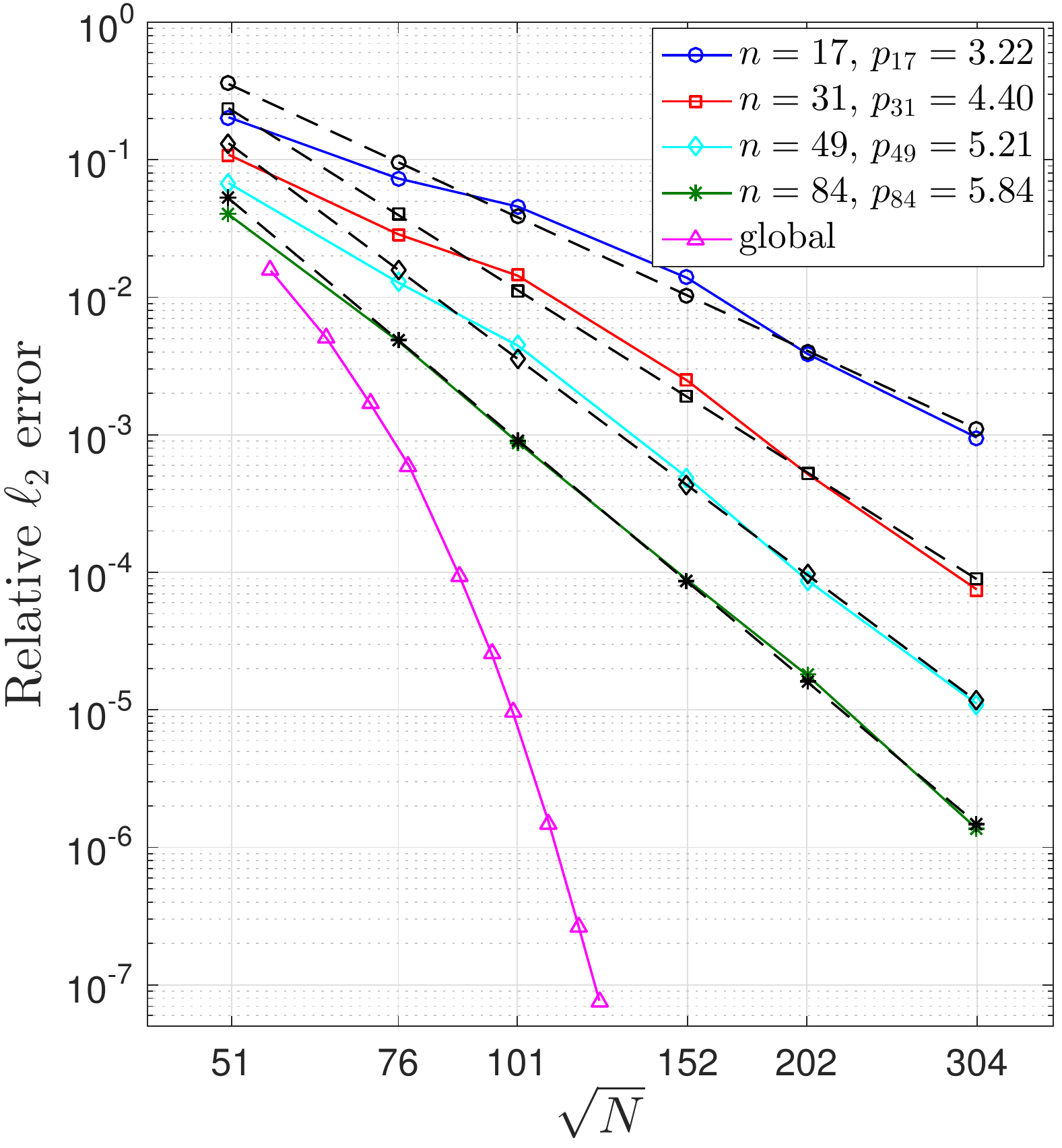} 	
}
&
{
	\includegraphics[width=0.46\textwidth]{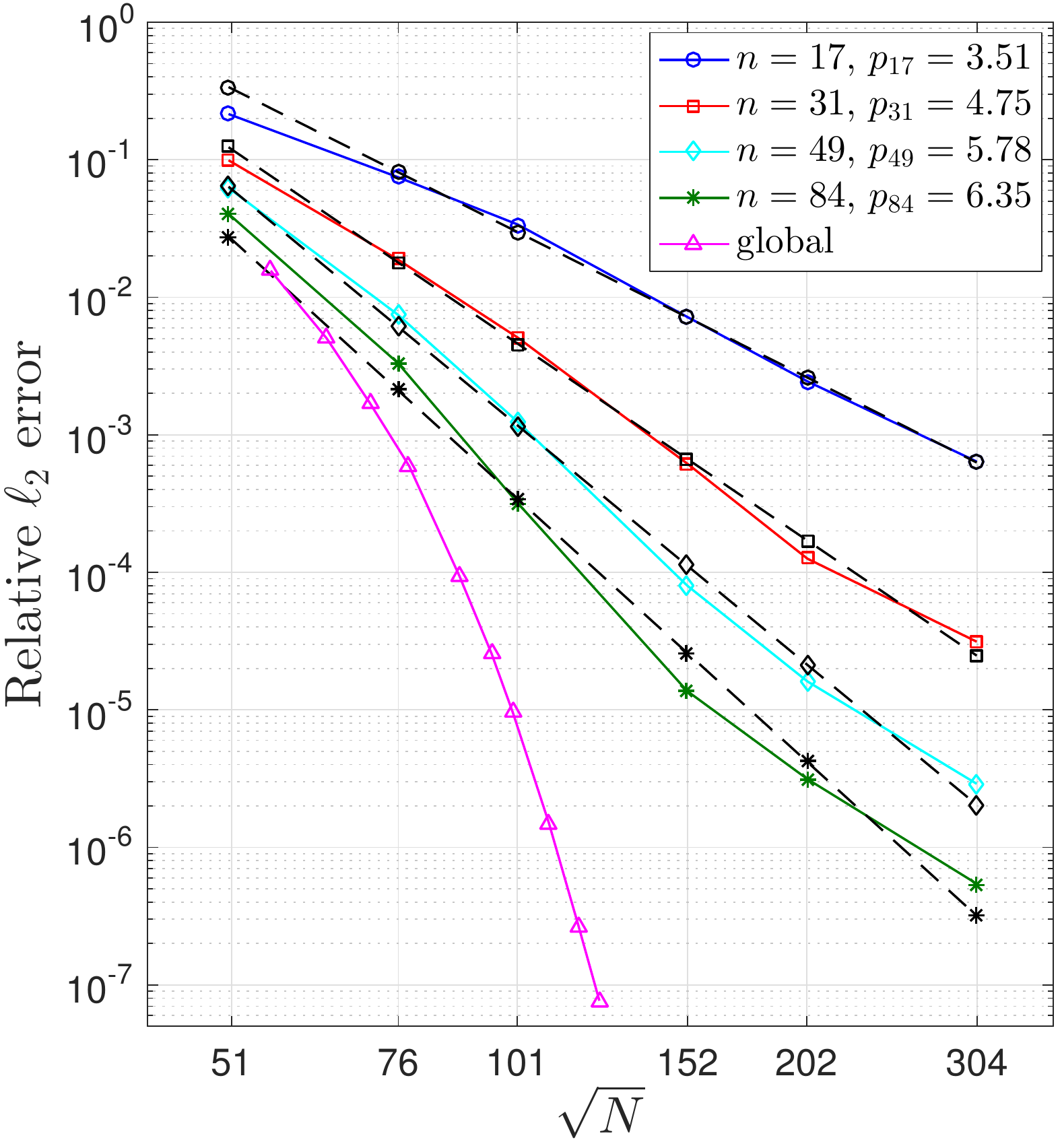} 	
	\label{fig:prbf_drel1}
}
\cr
{
	\includegraphics[width=0.46\textwidth]{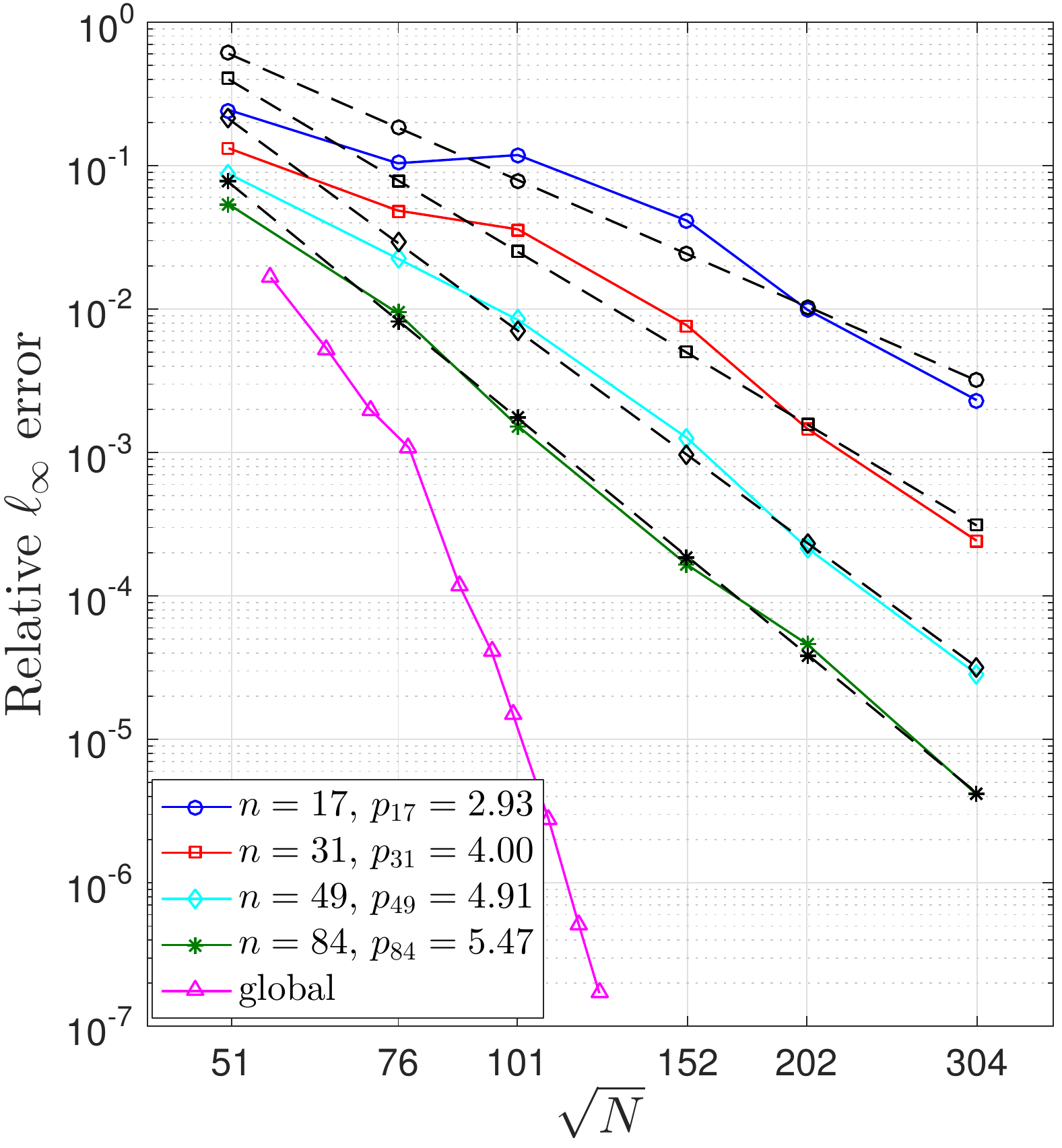} 	
}
&
{
	\includegraphics[width=0.46\textwidth]{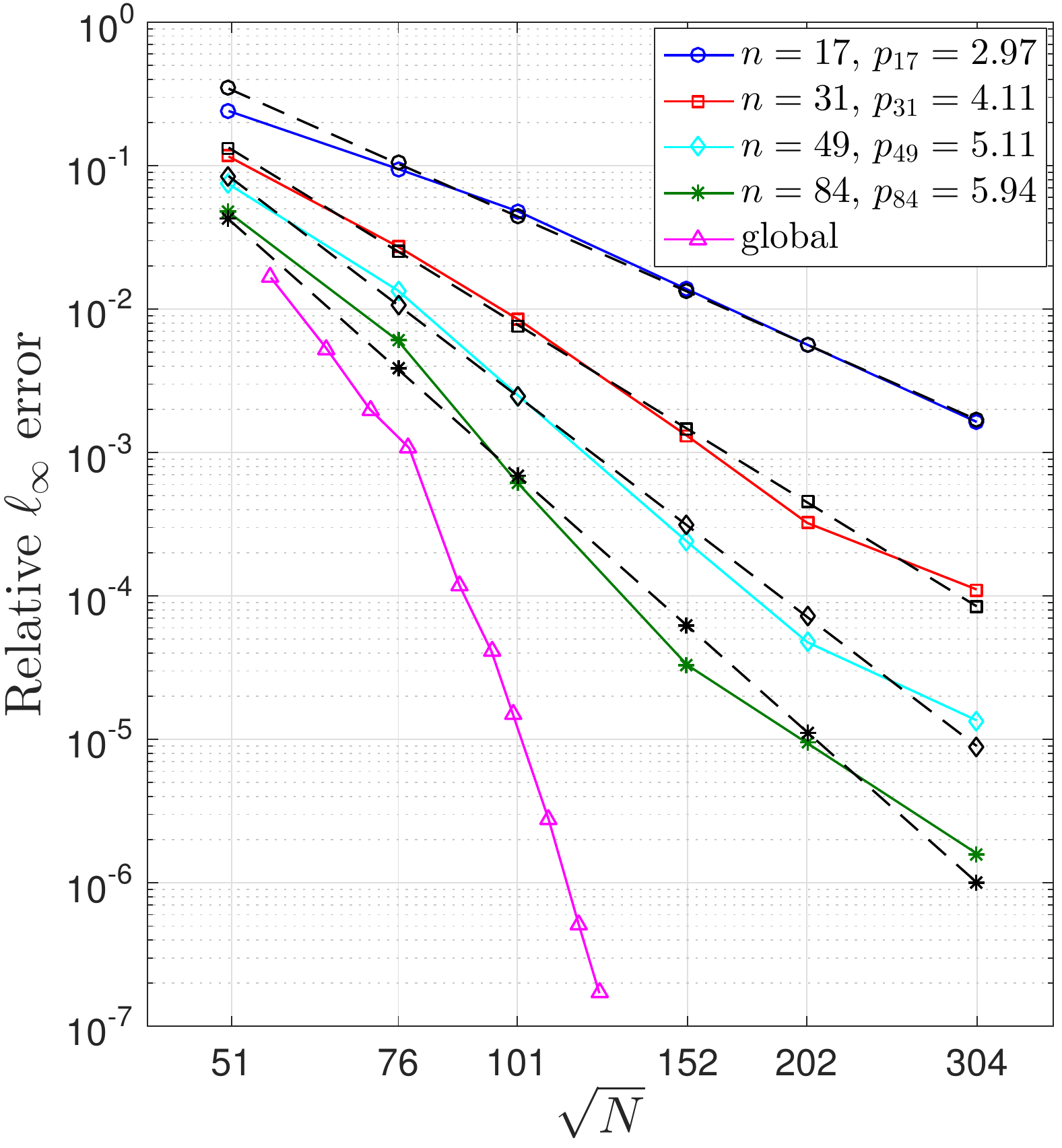} 	
	\label{fig:prbf_drel2}
}
\end{tabular}
\caption{Convergence results of the relative $\ell_2$ (top row) and $\ell_{\infty}$ (bottom row) errors the deformational flow of two Gaussian bells test case.  The results for the global RBF method are included in all plots for comparison purposes.  The dashed lines are the lines of best fit to the data (without the first point included) of the form $C_n N^{-p_n/2})$.  The values of $p_n$, which estimate the order of accuracy of the methods, are listed in the legend.  Lines of best fit are omitted for the global case since the convergence appears to be faster than any polynomial order.}
\label{fig:lrbf_def1}	
\end{figure}

The convergence results for the relative $\ell_2$ and $\ell_{\infty}$ errors for this test are displayed in Figure \ref{fig:lrbf_def1}.  We see from these plots that the norm of the errors for global RBF method appear to converge faster than any polynomial rate, which is expected since the solution is $C^{\infty}(\Sph)$ (see Section \ref{sec:global_rbfs}). The convergence rates of the errors for the local RBF and RBF-PU methods are also higher for this smooth test case.  Unlike the two previous tests, we see that these convergence rates also increase as $n$ increases.  The RBF-PU method appears to have a higher convergence rate for the same $n$ than the local RBF method, and the errors for the RBF-PU method are lower for each corresponding $n$ and $N$ value. This is likely due to the global smoothness of the RBF-PU interpolant.  

\begin{figure}[htb]
\centering
\begin{tabular}{cc}
{
	\includegraphics[width=0.41\textwidth]{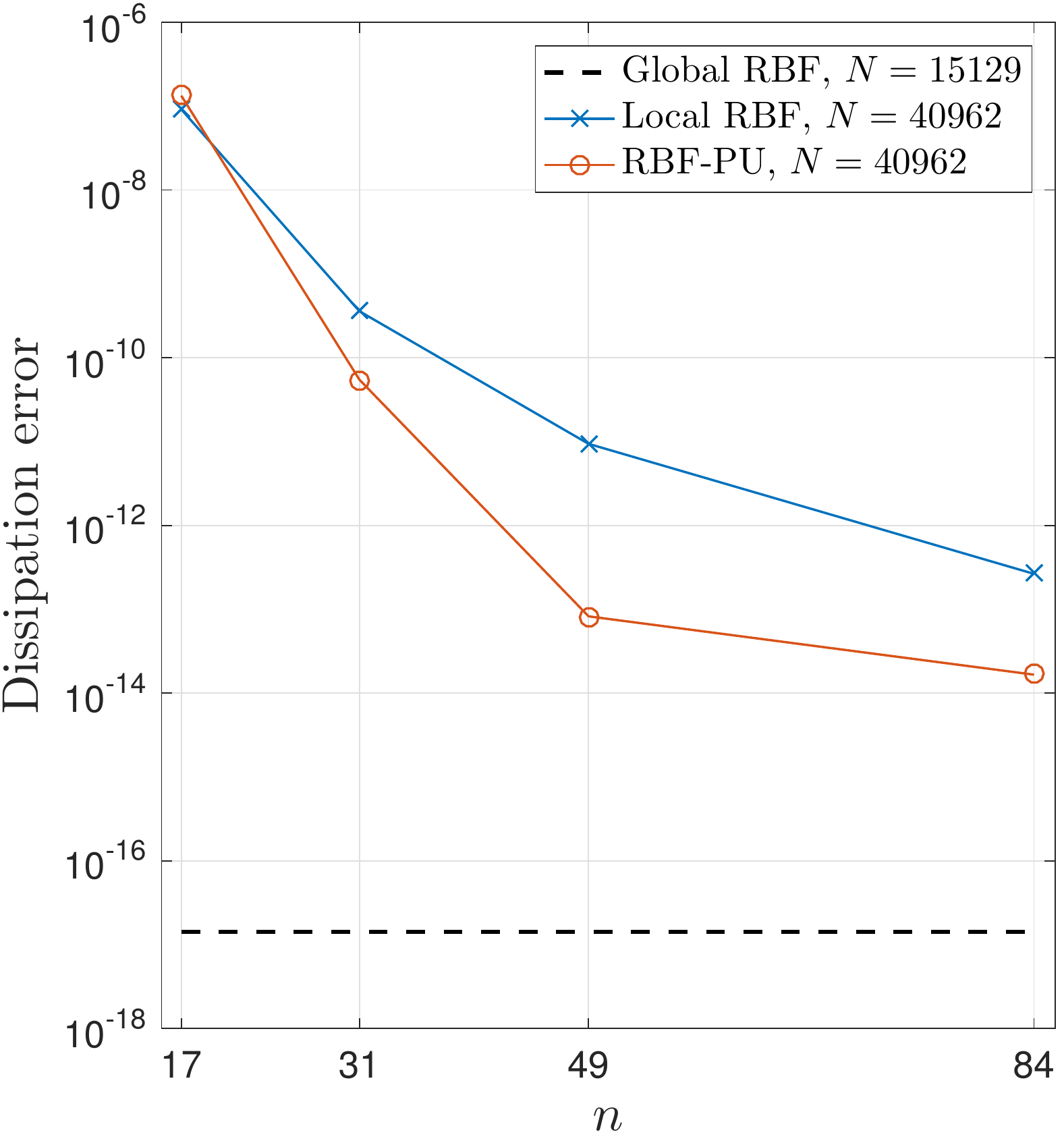} 	
}
&
{
	\includegraphics[width=0.41\textwidth]{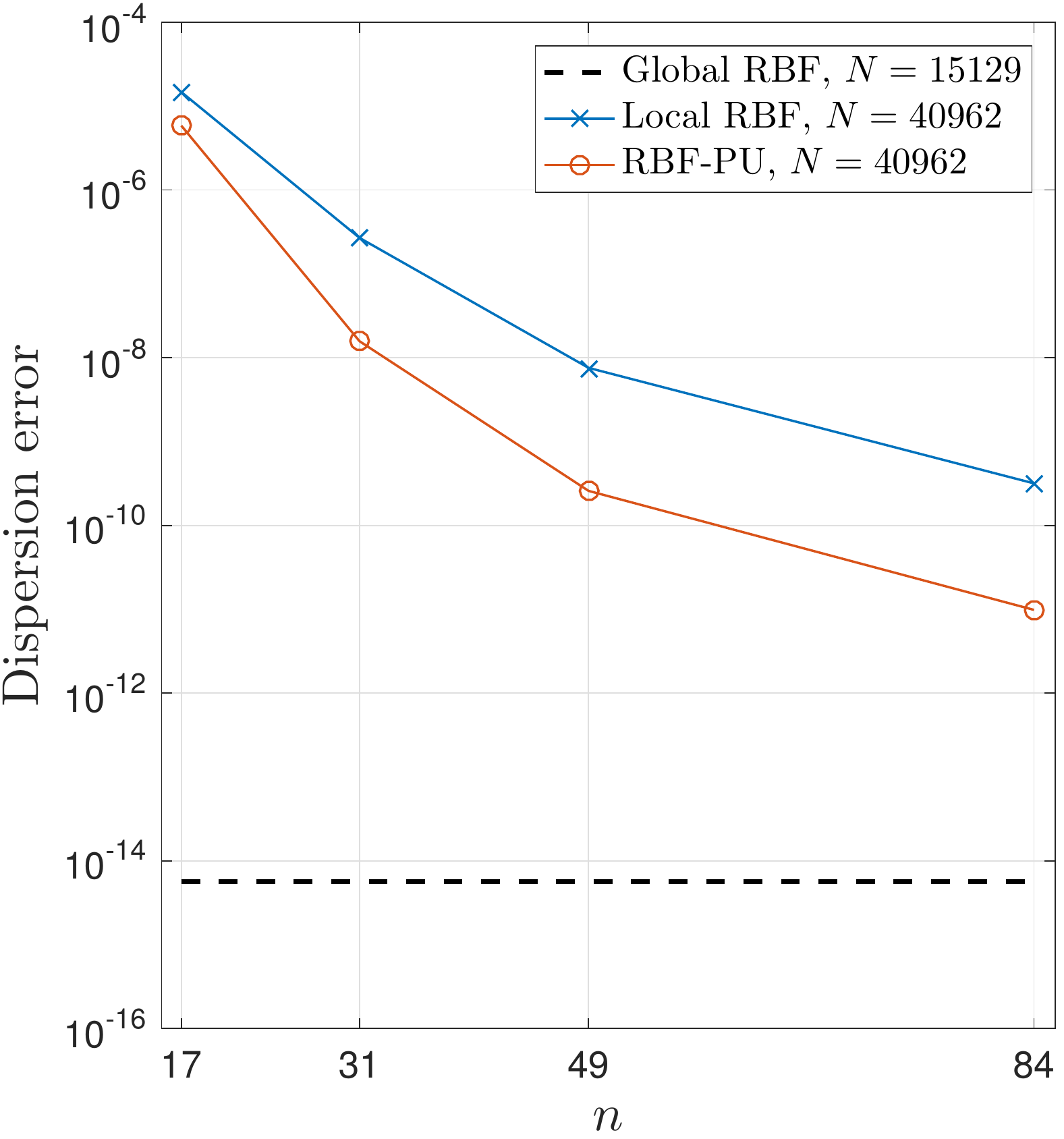} 	
}
\end{tabular}
\caption{Relative dissipation (left) and dispersion (right) errors \eqref{eq:dissipation}--\eqref{eq:dispersion} for the deformational flow test case of two Gaussian bells after one revolution as $n$ increases in the local RBF and RBF-PU methods.  The global RBF method does not have a dependence on $n$ and is included as a dashed line for reference.}
\label{fig:diss_disp_deform_gaussian}	
\end{figure}

Figure \ref{fig:diss_disp_deform_gaussian} displays the relative dissipation and dispersion errors just like the other test cases.  We see that the dissipation and dispersion errors are smallest for the global method and largest for the local RBF method.  We also see that increasing $n$ for local and RBF-PU methods leads to a much larger decrease in both the dissipation and dispersion errors than the previous two test cases.   We omit plots of the errors over the sphere since the dispersion errors are not a real issue for this test case.  

Time traces of the mass conservation errors are displayed in Figure \ref{fig:consDeformFlow} (b).   As with the dissipation/dispersion errors, the global method shows the best results followed by the RBF-PU method and then the local RBF method.  All of these methods show a slight growth in the conservation errors towards the end of the simulation time.   While not reported here, increasing $n$ in the local and RBF-PU methods does lead to a decrease in the conservation errors.

%The convergence rates of the errors for the local RBF and RBF-PU methods are also higher for this smooth test case, as shown Figure \ref{fig:lrbf_def1}.  Unlike the previous tests, we see that these convergence rates also increase as $n$ increases.  Also, the RBF-PU method has slightly higher convergence rates for the same $n$ than the local RBF method, and the errors for the RBF-PU method are lower (up to an order of magnitude) for each corresponding $n$ and $N$ value.  This is likely due to the global smoothness of the RBF-PU interpolant.  The second and third rows of plots in Figure \ref{fig:lrbf_def1} show similar trends in the dissipation and dispersion errors.  These results seem to indicate that for a smooth initial condition, having global smoothness in the interpolant can help decrease the overall error, in addition to the dissipation and dispersion errors.  However, as we show next, this global smoothness comes with a higher computational cost.

\subsection{Cost versus accuracy}\label{ssec:costvsaccuracy}
\begin{figure}[h!]
\centering
\begin{tabular}{cc}
Local RBF Method & RBF-PU Method \\
{
	\includegraphics[height=0.35\textheight]{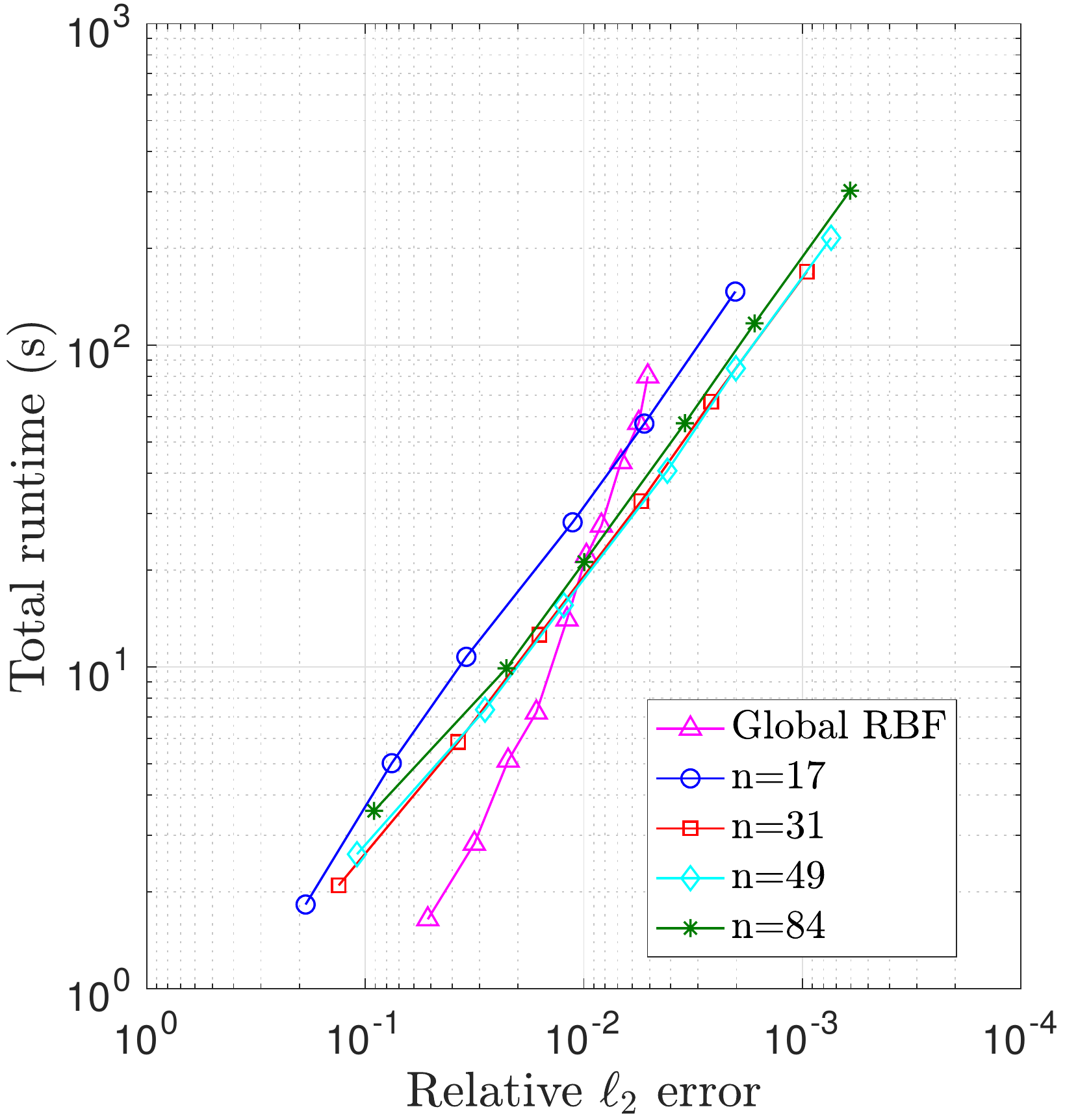} 	
} 
&
{
	\includegraphics[height=0.35\textheight]{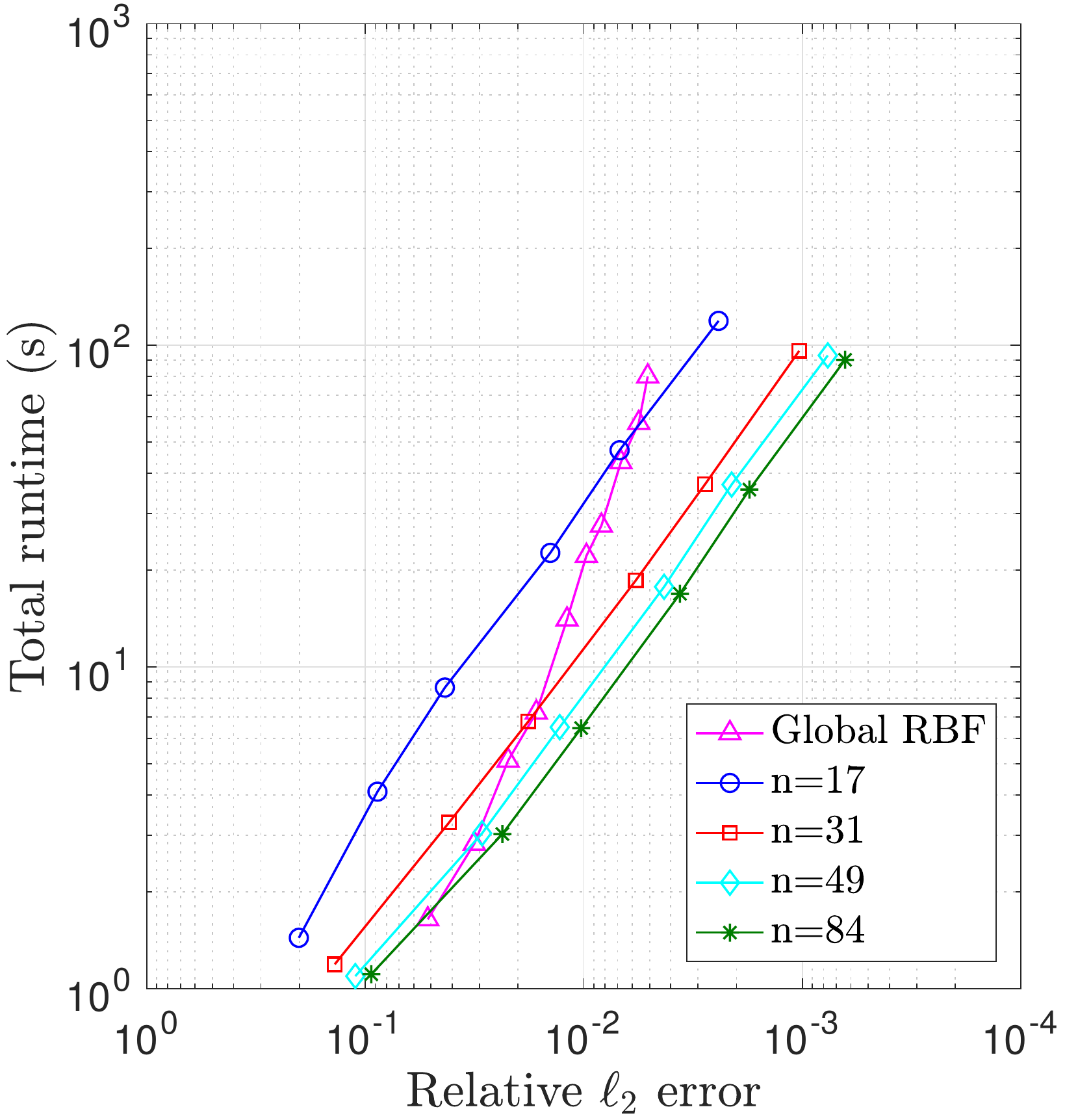}
%	\label{fig:gvpuc}
} \\
\multicolumn{2}{c}{(a) Deformational flow: Cosine Bells}\\
\phantom{space} \\
{
	\includegraphics[height=0.35\textheight]{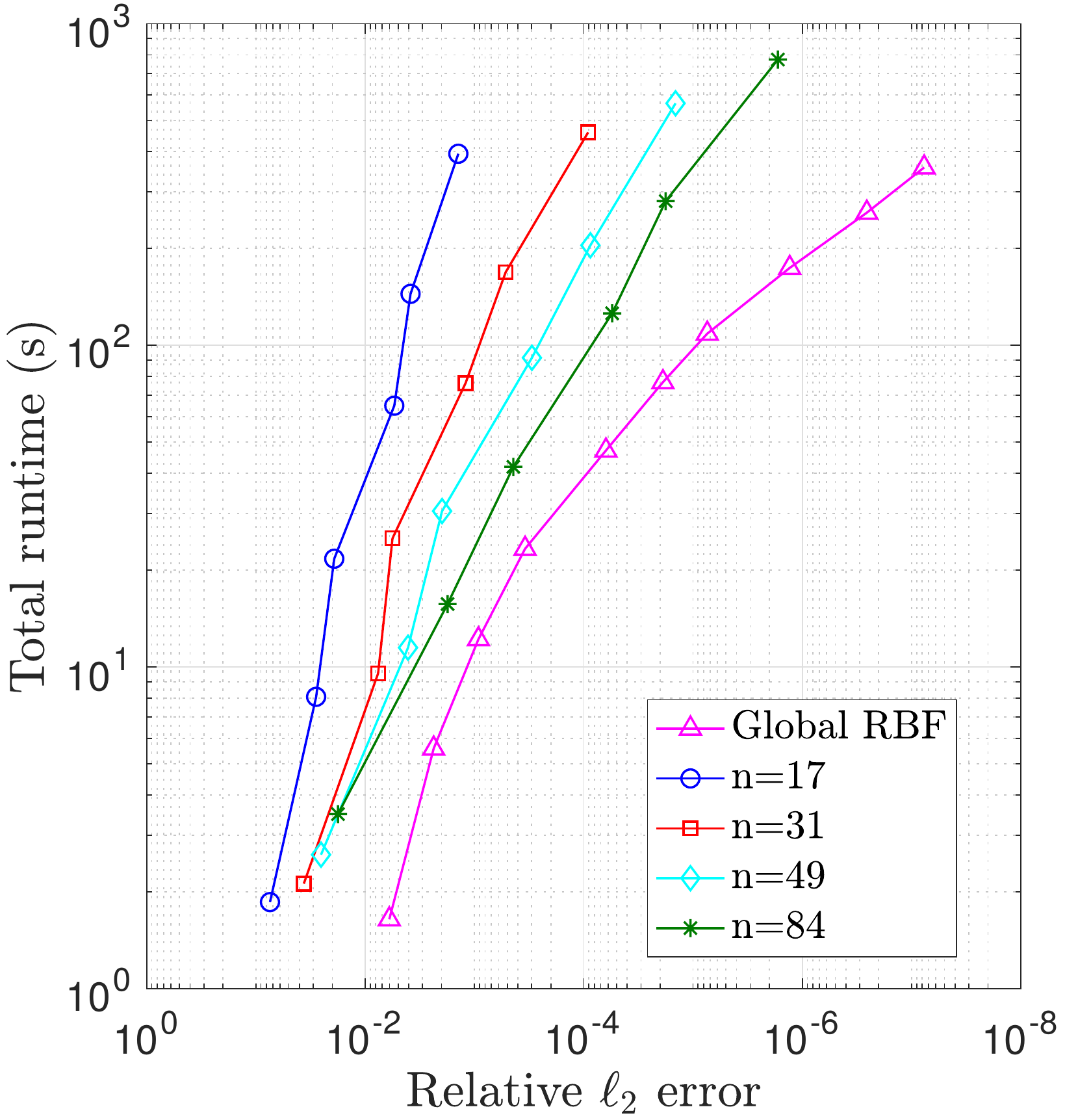} 	
}
&
{
	\includegraphics[height=0.35\textheight]{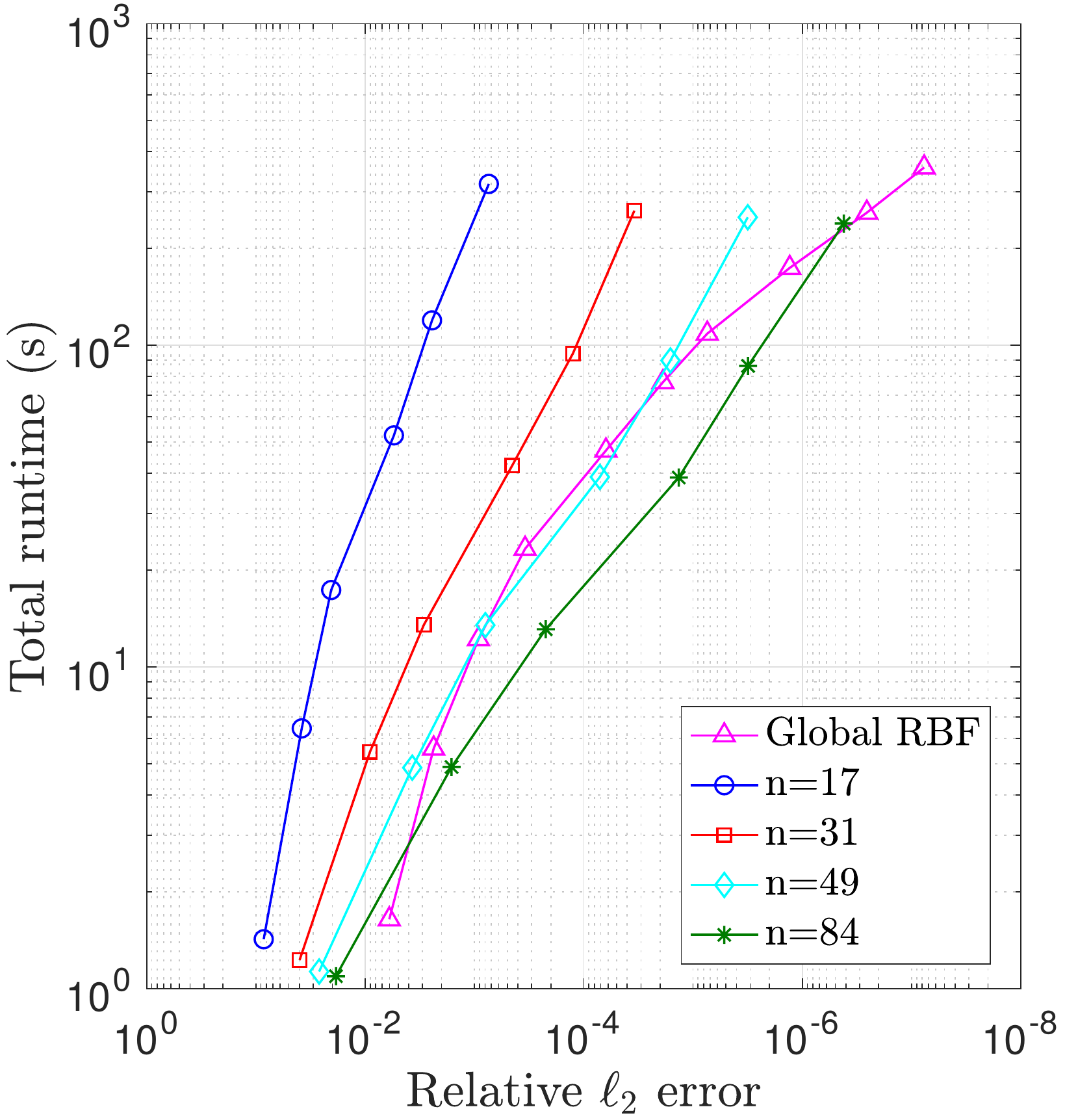}
} \\
\multicolumn{2}{c}{(b) Deformational flow: Gaussian Bells}\\
\end{tabular}
\caption{Computational cost (runtime) verses accuracy for the deformational flow test case using (a) Cosine Bells and (b) Gaussian bells.  The left column shows the results for the local RBF method, while the right is for the RBF-PU method.  The global RBF results are repeated in all plots for comparison purposes.  The runtime is for the total simulation and does not include any pre-computational costs.}
\label{fig:accvcost}	
\end{figure}

\begin{table}
\centering
\begin{tabular}{|c||cccccccccc|}
\hline
\rowcolor{Gray}\multicolumn{11}{|c|}{Local RBF and RBF-PU} \\
\hline
$N$ & 2562 & 5762 & 10242 & 23042 & 40962 & 92162 & & & & \\
\hline
CB, $\Delta t$ & 5/20 & 5/25 & 5/30 & 5/35 & 5/40 & 5/45 & & & &  \\
GB, $\Delta t$ & 5/20 & 5/40 & 5/60 & 5/80 & 5/100 & 5/120 & & & &  \\
\hline
\hline
\rowcolor{Gray}\multicolumn{11}{|c|}{Global RBF} \\
\hline
$N$ & 3136 & 4096 & 5041 & 6084 & 7744 & 9025 & 10000 & 11881 & 13689 & 15129 \\
\hline
CB, $\Delta t$ & 5/20 & 5/20 & 5/25 & 5/25 & 5/30 & 5/35 & 5/35 & 5/40 & 5/40 & 5/45 \\
GB, $\Delta t$ & 5/20 & 5/40 & 5/60 & 5/80 & 5/100 & 5/120 & 5/140 & 5/160 & 5/180 & 5/200 \\
\hline
\hline
\end{tabular}
\caption{\revtwo{Time-steps used for the numerical results in Sections \ref{ssec:costvsaccuracy} and \ref{ssec:morecompare} for the deformational flow test cases of the cosine bells (CB) and Gaussian bells (GB).\label{tbl:timesteps}}}
\end{table}

To properly compare all three methods, it is important to look at their cost (measured in wall-clock time) versus accuracy.  In this section we give such a comparison for the deformational flow test case using both the non-smooth and smooth initial conditions.  We examine the simulation time (ignoring preprocessing costs) and accuracy of each method (in the relative $\ell_2$ norm), and how these change with $N$ and $n$. All tests were run on a Linux workstation with a 3 GHz Intel Core i7-3930K (12 logical cores) and 32 GB of RAM.  All the codes were written and executed in MATLAB (version 2016a) in standard double precision, but without any explicit parallelization.  \revtwo{Unlike the previous sections, we here vary the time-steps with the number of nodes $N$ to approximately optimize the overall efficiency of the methods.  The time-steps used for the various node sizes are listed in Table \ref{tbl:timesteps}.  As mentioned previously, these time-steps were selected so that spatial errors in all the methods dominate.}

The results for the non-smooth cosine bells are displayed in Figure \ref{fig:accvcost} (a).  \revtwo{We see that for low accuracy, the global RBF method has a lower overall cost than both the local RBF method and for the RBF-PU method with $n=17$ and $n=31$.}  For higher accuracies, the cost of the global method increases much more rapidly than the other two methods because of its $O(N^2)$ complexity and, for high enough accuracies, both the local RBF and RBF-PU methods are more efficient overall for all values of $n$.  Comparing the  local RBF and RBF-PU methods, we see that, for a fixed $n$, the latter has a lower computational cost for a given accuracy.  Also, for this non-smooth test case, it does not appear to offer much benefit in terms of cost to use large $n$ with the local method, whereas the RBF-PU method shows slightly better efficiency with increasing $n$.

The results for the smooth Gaussian bells test case are given in Figure \ref{fig:accvcost} (b).  \revtwo{These plots clearly illustrate the advantage of the global RBF method with a smooth kernel over the local RBF method and the RBF-PU method with smaller $n$ when used on a problem with a smooth solution.  Extrapolating out, we see that global method is able to reach a much smaller $\ell_2$ error for same runtime for both the local and RBF-PU methods.  However, the global method will exhaust the memory resources of the machine much more quickly than the local and RBF-PU methods, and its results depend more heavily on memory latency.} The figure also shows that for a smooth initial condition, it pays to use larger values of $n$, with this being much more beneficial for the RBF-PU method.  Comparing the local and RBF-PU methods for this smooth initial condition, we see that overall, in our implementation, the RBF-PU method again gives higher accuracy for the same cost.  

It is important to note that these are serial implementations of all the methods. The local RBF method is easily parallelized on SIMD architectures, allowing for speedups of 2-8 times over a serial implementation, as shown for the related RBF-FD method in~\cite{Bollig12,Tillenius2015406}. The RBF-PU method also has promising parallelization properties using its patch-based structure. In contrast, the global RBF method with infinitely-smooth kernels requires a domain-decomposition style approach to parallelization.  We thus expect parallel versions of both the local and RBF-PU methods to give much better results in cost versus accuracy studies than presented here. We also note that, in all our tests, the global and RBF-PU methods were able to utilize MATLAB's automatically multithreaded BLAS kernels more efficiently than the local RBF method, since these two methods inherently use fewer, denser matrix-matrix multiplications (BLAS-3) operations. 

\revtwo{
\subsection{Comparison with other methods}\label{ssec:morecompare}

\begin{figure}[htb]
\centering
\includegraphics[height=0.35\textheight]{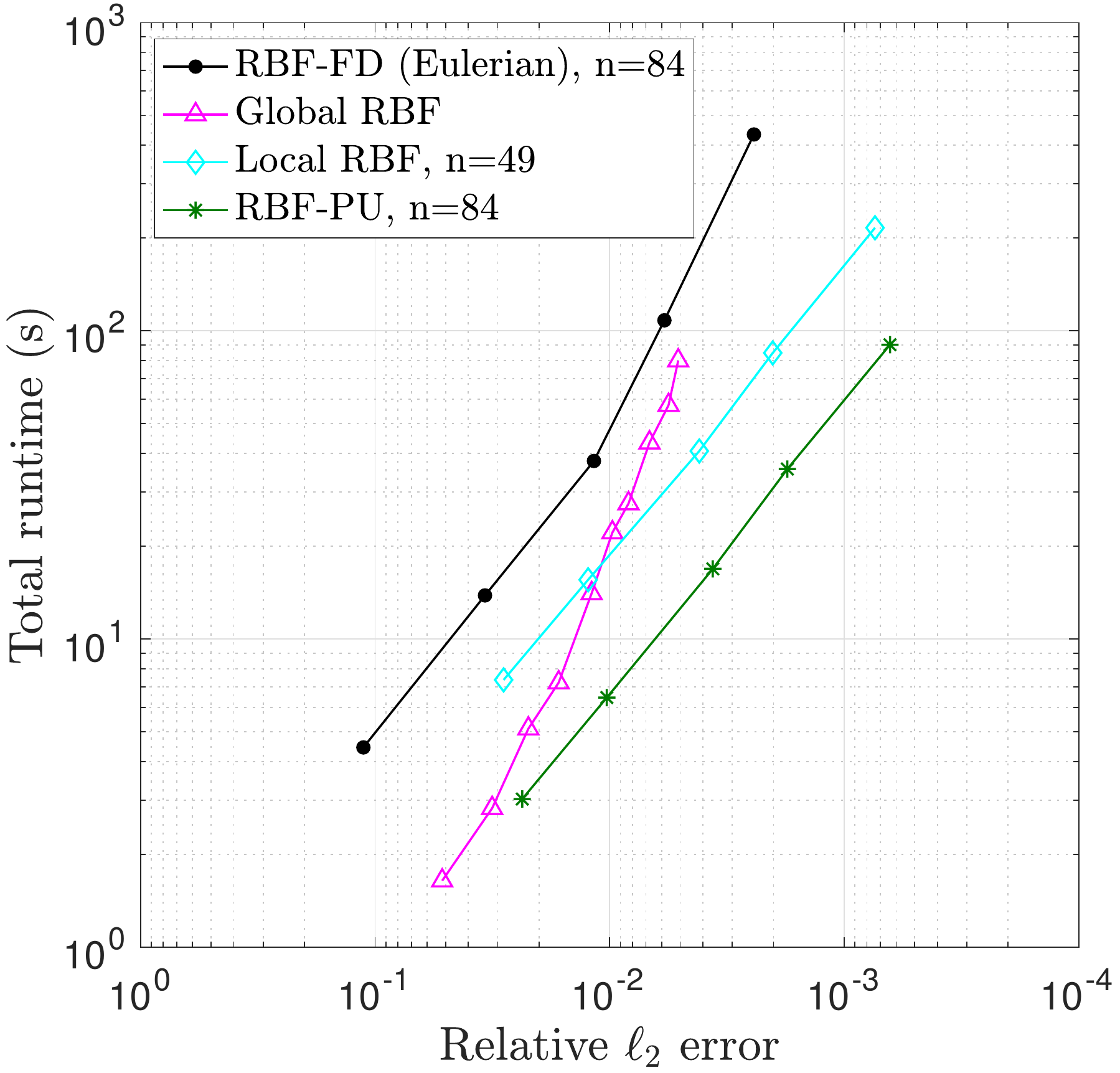} 	
\caption{\revtwo{Comparison of the three new SL RBF schemes and the Eulerian RBF-FD scheme from~\cite{FoL11}.  The plot shows the computational cost (runtime) verses relative $\ell_2$ error for the Cosine Bells deformational flow test case.  The stencil sizes for the global RBF, local RBF, RBF-PU, and Eulerian RBF-FD methods are shown in the legend.  For the latter three methods, the same node sets were used of sizes $N=5762$, 10242, 23042, 40962, and 92162.\label{fig:EulerianRBFFD}}}
\end{figure}

\begin{table}
\centering
\begin{tabular}{|c||c|c|c|c|}
\hline
             & Resolution & Degrees of  & Time-step ($\Delta t$), & Relative $\ell_2$ \\
Method & (in degrees) & freedom ($N$) & non-dimensional & error \\
\hline
\hline
\rowcolor{Gray}\multicolumn{5}{|c|}{Cosine Bells}\\
\hline
%CSLAM~\cite{NairLauritzen2010} & $1.5\degree$ & 21600 & 5/120 & $3.28 \cdot 10^{-2}$ \\
CSLAM~\cite{LauritzenEtAl2012} & $0.75\degree$ & 86400 & 5/240 & $\approx 6 \cdot 10^{-3}$ \\
DG, $p=3$~\cite{NairLauritzen2010} & $1.5\degree$ & 38400 & 5/2400 & $1.39 \cdot 10^{-2}$ \\
RBF-FD, $n=84$~\cite{FoL11} & $1.5\degree$ & 23042 & 5/900 & $1.17 \cdot 10^{-2}$\\
%RBF-FD, $n=84$~\cite{FoL11} & $0.75\degree$ & 92162 & 5/2800 & $2.42 \cdot 10^{-3}$\\
Local RBF, $n=84$ & $1.5\degree$ & 23042 & 5/35 & $3.45 \cdot 10^{-3}$\\
RBF-PU, $n=84$ & $1.5\degree$ & 23042 & 5/35 & $3.63 \cdot 10^{-3}$\\
Global RBF & $1.64\degree$ & 15129 &  5/45 & $5.1 \cdot 10^{-3}$\\
%Global RBF 
\hline
\hline
\rowcolor{Gray}\multicolumn{5}{|c|}{Gaussian Bells}\\
\hline
%CSLAM~\cite{NairLauritzen2010} & $1.5\degree$ & 21600 & 5/120 & $2.46 \cdot 10^{-2}$ \\
CSLAM~\cite{LauritzenEtAl2012} & $0.75\degree$ & 86400 & 5/240 & $\approx 5 \cdot 10^{-4}$ \\
HOMME, $p=6$~\cite{LauritzenEtAl2014} & $1.5\degree$ & 29400 & 5/4800 & $\approx 3 \cdot 10^{-3}$ \\
%HOMME, $p=6$~\cite{LauritzenEtAl2014} & $0.75\degree$ & 117600 & 5/9600 & $6 \cdot 10^{-5}$ \\
RBF-FD, $n=84$~\cite{FoL11} & $1.5\degree$ & 23042 & 5/900 & $3.18 \cdot 10^{-4}$\\
%RBF-FD, $n=84$~\cite{FoL11} & $0.75\degree$ & 92162 & 5/2800 & $2.42 \cdot 10^{-3}$\\
Local RBF, $n=84$ & $1.5\degree$ & 23042 & 5/80 & $5.50 \cdot 10^{-5}$\\
RBF-PU, $n=84$ & $1.5\degree$ & 23042 & 5/80 & $1.35 \cdot 10^{-5}$\\
Global RBF & $1.64\degree$ & 15129 &  5/200 & $7.68 \cdot 10^{-8}$\\
\hline
\end{tabular}
\caption{\revtwo{Comparison of various methods from the literature for the two deformational flow test cases.  CSLAM is the conservative semi-Lagrangian multi-tracer transport scheme used in~\cite{LauritzenEtAl2012}, which is based on a cubed sphere grid.  DG is the discontinuous Galerkin scheme used in~\cite{NairLauritzen2010}.  This is an Eulerian scheme that uses $p=3$ degree polynomials (fourth-order accurate) and the cubed sphere grid.  HOMME is the High-Order Methods Modeling Environment scheme used in~\cite{LauritzenEtAl2014}.  This is also an Eulerian scheme using a cubed-sphere grid and the results given are for a continuous Galerkin formulation using $p=6$ degree polynomials.  RBF-FD is the mesh-free Eulerian scheme from~\cite{FoL11}  using $n=84$ point stencils and the Gaussian RBF with the shape-parameter selected using the formulas from~\cite{FlyerLehto2012}.  Values appearing with a $\approx$ were obtained from plots of the relative $\ell_2$ errors given in the corresponding referenced papers, as no exact values are reported.  For the RBF-FD method the results were generated by the present authors using the code from~\cite{FlyerLehto2012}.  The results given for the CSLAM, DG, and HOMME methods are for the non-filtered (or non shape-preserving) versions, which give the lowest errors for these test cases. Resolution (in degrees) is an approximate measure of the spacing of the nodes (or grid points) around the equator.\label{tbl:compare}}}
\end{table}

In this section, we compare the three new SL RBF methods against other methods for transport on the sphere found in the literature.  First, we compare the overall efficiency (wall-clock time vs.\ accuracy) of the SL RBF methods against the standard Eulerian RBF-FD method developed in~\cite{FoL11} (see also~\cite{FlyerLehto2012,Bollig12,Tillenius2015406}) on the deformational flow of two cosine bells test case. For the results presented here, we used the Gaussian RBF, similar to~\cite{FoL11,FlyerLehto2012}, with an $n=84$ point stencil.  For the shape parameter $\ep$, we used the approach suggested in~\cite{FlyerLehto2012} and selected $\ep$ to scale with $N$ as  $\ep = 0.063\sqrt{N}$.  As mentioned in the introduction, this method requires a hyperviscosity stabilization term of the form $(-1)^{k+1}\gamma_{N} \Delta^{k}$, where $\Delta$ is the Laplacian, $k$ is a positive integer, and $\gamma_{N}$ is a constant that decreases with $N$.  With the choices of RBF, $\ep$, and $n$ above, we found that selecting the hyperviscosity parameters as $k=6$ and $\gamma_{N} = 15*2^{-8}N^{-5}$ produced stable results with the time-steps used, without overdamping the solutions.  As in~\cite{FoL11,FlyerLehto2012}, we used RK4 as the time-integrator and node sets of size $N=5762$, 10242, 23042, 40962, and 92162, with corresponding time-steps $\Delta t = 5/400$, $5/700$, $5/900$, $5/1400$, and $5/2800$.  These are approximately the largest time-steps that could be used so that the solutions were stable and spatial errors dominate.  For the SL RBF methods, we used the time-steps listed in Table \ref{tbl:timesteps}.  The same machine described in the previous section was used for all the simulations.

Figure \ref{fig:EulerianRBFFD} shows the results of the comparison between the SL RBF methods and the Eulerian RBF-FD method.  We see that the local RBF and RBF-PU methods have a lower error for the same degrees of freedom and that they require less computational time than the Eulerian RBF-FD method.  This increase in efficiency comes from the fact that much larger time-steps can be used in the SL simulations.

For the second comparison, we also use the deformational flow test, but now for both the cosine bells and the Gaussian bells.  The focus here is on the errors for similar resolutions, degrees of freedom, and time-steps.  We also include three commonly used methods from the geosciences (CSLAM, DG, and HOMME) in the comparison.  Table \ref{tbl:compare} displays the results and the caption contains descriptions of the methods used.   We can see from the table that all three SL RBF methods compare quite favorably for these test cases to existing methods.  It should be noted, however, that the CSLAM, DG, and HOMME methods are all mass-conserving, whereas, the present version of the RBF SL methods presented here are not. 
}

\section{Summary}
\label{sec:summary}

In this article, we presented three new SL methods for transport on a sphere based on interpolation with global RBFs, local RBFs, and RBF-PU. The RBF framework allowed us to obtain either high-order convergence rates for smooth problems using only scattered nodes and Cartesian coordinates.  Using scale-free RBFs with the addition of spherical harmonics in the local and RBF-PU methods removes the need to choose a shape parameter and also avoids stagnation errors observed in other applications of these methods. Additionally, the SL framework appears to lend our methods intrinsic stability without the need for hand-tuned artificial hyperviscosity.  We summarize the features of our methods below:
\begin{itemize}
\item \revone{The global RBF method with smooth kernels appears to have the best cost-accuracy profile for problems with smooth solutions that can be resolved using a relatively small number of nodes. However, the gap between the global method and local/RBF-PU methods can be decreased by increasing the stencil/patch sizes, $n$. Furthermore, we expect that the local/RBF-PU methods will be more competitive for problems with smooth solutions with fine structures that need to be resolved with large numbers of points.}
\item \revtwo{For non-smooth solutions, the local and RBF-PU methods both out perform the global method for moderate to high accuracies.}
\item The RBF-PU method gives smaller errors, lower dispersion, and better conservation properties than the local RBF method for a comparable number of degrees of freedom, regardless of the smoothness of the solutions.
\item The local RBF method is easier to implement and may be more amenable to highly efficient implementations on SIMD architectures.
\item \revtwo{All three methods compare favorably to commonly used methods for transport on the sphere for the problems considered here.}
\end{itemize}
Our algorithms were designed for transport in an incompressible velocity field. However, compressible fields arise naturally in many biological and geophysical problems. We are currently working on extending our method to handle this case. In addition, there is a need for quasi-monotone (non-oscillatory) local RBF methods to ameliorate the numerical dispersion seen therein. We plan to address this issue in future work. Furthermore, the local RBF method can obtain very large speedups even in a serial implementation if used with the overlapped RBF-FD framework~\cite{Shankar2016}; we plan to explore this in future work as well. Finally, the RBF framework allows for a straightforward generalization of our methods to more general manifolds than  spheres~\cite{FuselierWright2014,SWFKJSC2014,LSW2016}.  \revtwo{While the methods discussed here are limited to transport on the sphere, we anticipate that they may be extended to the full non-linear shallow water equations using a similar approach as~\cite{LaytonSpotz2003}.} 
%~\cite{MagiThesis} \comment{GBW: This reference is obviously broken, but also seems a bit out of place.  The model is not just about transport it also involves moving surfaces which we haven't talked about.  Also, it should refer to the SIAM Applied Math paper that just appeared}.

\section*{Acknowledgements}
The authors thank the two anonymous referees and associate editor for their comments, which helped to improve the paper.  VS acknowledges support for this project under NSF-CCF 1714844, NSF-DMS 1160432, and NSF-DMS 1521748. GBW acknowledges funding support for this project under grants NSF-CCF 1717556 and NSF-ACI 1440638.

%\appendix
%\input{Appendix}

\section*{References}
\bibliography{article_refs_mod}

\end{document}